\documentclass[12pt]{amsart}
\usepackage[margin=1in,centering]{geometry}
\usepackage{amssymb}
\usepackage{graphicx}
\usepackage[usenames,dvipsnames]{color}
\usepackage[implicit=true, 
  colorlinks=true,linkcolor=Blue,citecolor=PineGreen,urlcolor=BrickRed]{hyperref}

\newcommand{\bb}[1]{\mathbb{#1}}

\newtheorem{thm}{Theorem}[section]
\newtheorem{cor}[thm]{Corollary}
\newtheorem{lem}[thm]{Lemma}
\newtheorem{prop}[thm]{Proposition}
\newtheorem{prob}{Problem}
\newtheorem{conj}{Conjecture}
\theoremstyle{definition}
\newtheorem{defn}[thm]{Definition}
\theoremstyle{remark}
\newtheorem{rem}[thm]{Remark}

\numberwithin{equation}{section}
\numberwithin{figure}{section}

\newif\ifStandalone \Standalonetrue
\providecommand{\SeeColorFig}[1]{}
\providecommand{\SeeColorFigs}[2]{}
\newif\ifColorSection \ColorSectionfalse
\providecommand\windex[1]{#1\index{#1}}

\renewcommand{\:}{\colon\! }
\providecommand{\StretchyLines}{\relax }
\providecommand{\url}[1]{\texttt{#1}}

\begin{document}

\ifStandalone
\title{Milnor's Problem on the Growth of Groups and its Consequences}
\author{Rostislav Grigorchuk}
\address{Department of Mathematics;
Mailstop 3368; Texas A\&M University; College Station, TX 77843-3368,
USA }

\email{grigorch@math.tamu.edu}

\else
\address{Rostislav Grigorchuk,
Department of Mathematics;
Mailstop 3368; Texas A\&M University; College Station, TX 77843-3368,
USA }
\email{grigorch@math.tamu.edu}
\fi

\thanks{The author is supported by NSF grant DMS--1207699, ERS grant GA 257110 ``RaWG'',  and by the Simons Foundation}

\dedicatory{Dedicated to John Milnor on the occasion of his  80th birthday.}
\begin{abstract}
We present a survey of results related to Milnor's problem on group growth.  We discuss the cases of
polynomial growth and exponential but not uniformly exponential growth; the main part of the article is
devoted to the intermediate (between polynomial and exponential) growth case. A number of related topics
(growth of manifolds, amenability, asymptotic behavior of random walks) are considered, and  a number of open
problems are suggested.
\end{abstract}
\maketitle

\section{Introduction}   The notion of the growth of a finitely generated group was introduced by A.S. Schwarz
(also spelled Schvarts and \v{S}varc)) \cite{schvarts:55}  and independently by Milnor
\cite{milnor:problem,milnor:note68}. Particular studies of group growth and their use in various situations
have appeared in the works of Krause \cite{krause:growth53},  Adelson-Velskii and Shreider
\cite{adelson:banach57}, Dixmier \cite{dixmier:growth60}, Dye \cite{dye:59,dye:groups63}, Arnold and Krylov
\cite{arnold_krylov:63}, Kirillov \cite{kirillov:dynamic67}, Avez \cite{avez:riemann70}, Guivarc'h
\cite{guivarch:70,guivarch:growth71,guivarch:73}, Hartley, Margulis, Tempelman and other researchers.
The note of Schwarz  did not attract a lot of attention in the mathematical community, and was essentially
unknown to mathematicians both in the USSR and the West (the same happened with papers of Adelson-Velskii,
Dixmier and of some other mathematicians).  By contrast, the note of Milnor \cite{milnor:note68}, and especially the
problem raised by him in \cite{milnor:problem}, initiated a lot of activity and opened  new directions in
group theory and areas of its applications.

The motivation for  Schwarz  and Milnor's studies on the growth of groups were of geometric character. For
instance, it was observed by Schwarz  that the rate of  volume growth of the universal cover $\tilde M$ of a
compact Riemannian manifold $M$ coincides with the rate of growth of the fundamental group $\pi_1(M)$
\cite{schvarts:55}. At the same time Milnor and Wolf demonstrated that  growth type of the fundamental group
gives some important information about the curvature of the manifold. A relation between the growth of a group
and its amenability was discovered by Adel'son-Vel'skii and Shreider (spelled also as
\v{S}re\u{i}der) \cite{adelson:banach57}; specifically,
subexponential growth implies the amenability of the group,  i.e., the existence of an invariant mean).

\medskip
The problem posed by Milnor focuses on  two main questions:
\begin{enumerate}
\item Are there groups of intermediate growth between polynomial and exponential?
\item What are the groups with polynomial growth?
\end{enumerate}

Moreover, Milnor formulated a remarkable conjecture about the coincidence of the class of groups with polynomial
growth and the class of groups  containing  a nilpotent subgroup of finite index (i.e.\ with the class of
virtually nilpotent groups) which was later proved by Gromov \cite{gromov:poly_growth}.

The first part of  Milnor's question was formulated originally by him in the following form: ``Is it true that
the growth function of every finitely generated group is necessarily equivalent to a polynomial or to the
function $2^n$?'' This was answered in the negative by the author in 1983 \cite{grigorch:milnor,grigorch:degrees}.
Despite the negative character of the answer,  the
existence of groups of intermediate growth  made  group theory and  the areas of its applications much richer.
Eventually it led to the appearance   of  new directions in group  theory:  self-similar groups
\cite{gns00:automata,nekrash:self-similar},  branch groups
\cite{grigorch:branch,grigorch:jibranch,bar_gs:branch}, and   iterated monodromy groups
\cite{bartholdi_gn:fractal,nekrash:self-similar}.
 Moreover, completely new  methods, used in the study of groups of intermediate growth, stimulated
intensive studies of groups generated by finite automata - a direction on the border between computer science
and algebra initiated by V.M.~Glushkov in the beginning of the 1960s
~\cite{glushkov:automata,aleshin-book,gns00:automata,bondarenko_gkmnss:full_clas32}.   Groups of intermediate
growth were used and continue to be used in many different situations: in the study of amenable
groups~\cite{grigorch:degrees,grigorch:example}, in topology~\cite{freedman_teich:subexp95}, in  theory of
random walks on groups~\cite{erschler:ICM11,kaiman:munchhausen,erschler:boundary04}, in  theory of operator
algebras~\cite{grigorch:solved,grigorchuk-n:schur,grigorch:expanders11}, in dynamical
systems~\cite{bartholdi_nekr:quandratic,bux_p:iter_monodromy}, in
percolation~\cite{muchnik_pak:percolation01}, in the study of cellular automata~\cite{machi_mignosi:garden93},
in the theory of Riemannian surfaces~\cite{grigorch:topological89}, etc.

 As was already mentioned, the study of the growth of groups
was introduced by  Schwarz and Milnor.  However, some preliminary cases  were considered by H.U.~Krause in his
Ph.D thesis, defended in Zurich in 1953 \cite{krause:growth53}, and geometric growth was considered around 1953 by
V.A.~Efremovich \cite{efremovic:growth53}. Later, the question of group growth  was considered in certain situations by
Adelson-Velskii and Shreider, Dixmier, Dye, Arnold and Krylov, Guivarc'h, Kirillov, Margulis, Tempelman, and
other researchers; however, as has been already mentioned, the systematic study of group growth only started
with Milnor's work in 1965.

\medskip

The first period of studies concerning group growth (in the 1960s and 1970s) was dedicated to the study of
growth of nilpotent and solvable groups.  During a short period of time, Milnor, Wolf, Hartley, Guivarc'h and
Bass discovered that nilpotent groups have  polynomial growth of integer degree, given by a number
(\ref{polyn5}) expressed in terms of the lower central series. The converse  fact, namely that polynomial
growth implies the virtual nilpotence of the group  was proven by M.~Gromov in his remarkable paper
\cite{gromov:poly_growth}, which stimulated a lot of activity in different areas of mathematics. Gromov's
proof uses the idea of the limit of  a sequence of metric spaces, as well as  Montgomery and Zippin's
solution of Hilbert's 5th problem \cite{montgomery_zip:topolog55}. Van den Dries and Wilkie
\cite{dries_wilkie:gromov84} used methods of nonstandard analysis to explore  Gromov's idea in order to
slightly  improve his result (in particular, to outline a broader  approach to the notion of a cone of a
group), while B.~Kleiner, using the ideas of  Colding-Minicozzi \cite{colding_minic:harmonic97} and applying
the techniques of harmonic functions,  gave a proof of Gromov's theorem which doesn't rely on the techniques
in the Montgomery-Zipin solution to Hilbert's 5th problem. This approach was explored  by Y.~Shalom and
T.~Tao to get effective and quantitative results about the polynomial growth case.

The author's paper~\cite{grigorch:hilbert}, followed by the ICM Kyoto paper~\cite{grigorch:ICM90},
raised a very interesting question which we discuss in detail in this article (section~\ref{gapconj}) and
call the \emph{Gap Conjecture}. The conjecture states that if the growth degree of a finitely generated group is
strictly smaller than $e^{\sqrt n}$ then it is polynomial and the group is virtually nilpotent.  If proven,
this conjecture would give a far reaching generalization of Gromov's polynomial growth theorem. There are
several results, some presented here, supporting the conjecture.

It is interesting that the $p$-adic analogue of Hilbert's 5th Problem,  solved  by M.~Lazard
\cite{lazard:analytic65}, was used by the author~\cite{grigorch:hilbert,grigorch:ICM90} to obtain results
related to the Gap Conjecture;  these are  discussed in
section~\ref{gapconj}.  See also the work of A.~Lubotzky, A.~Mann and
D.~Segal~\cite{lubotzky_mann:polyn91,lubotzky_man_seg:93} for a description of the groups with polynomial
subgroup growth.

The notion of growth can be defined for many algebraic and combinatorial objects, in particular for
semigroups, associative and Lie algebras, graphs and discrete metric spaces as discussed, for instance, in
\cite{harpe_cg:paradoxical}. There is much more freedom in the asymptotic behavior of growth in each of these
classes of algebraic objects, and there is no hope to obtain a general result characterizing  polynomial
growth similar to Gromov's theorem. But in~\cite{grigorch:semigroups_with_cancel88}, the author was able to
extend the theorem of Gromov to  the case of cancellative semigroups (using the notion of nilpotent semigroup
introduced by A.I.~Malcev~\cite{malcev:nilpotent53}).

Let us say a few words about the exponential growth case.  Perhaps typically, a finitely generated group has
exponential growth. The fact that any solvable group which is not virtually nilpotent has exponential growth was
established by Milnor and Wolf  in 1968 \cite{milnor:solv68,wolf:growth68}. A direct consequence of the
\windex{Tits Alternative} \cite{tits:alternative72} is  that a finitely generated linear group  has  exponential growth if
it is not virtually nilpotent.

The exponential growth of a group immediately follows from  the existence of a free subgroup on two
generators, or even a free sub-semigroup on two generators. There are many classes of groups (i.e.\
solvable non-virtually nilpotent groups, non-elementary Gromov hyperbolic groups, linear groups without free
subgroup on two generators, etc.) which are known to contain such objects.  But many infinite torsion groups
also have exponential growth, as was first proved by S.I.~Adian in the case of free Burnside groups of large
odd exponent \cite{adian:b-burnside}.

\medskip
The second period of studies of group growth begins in the 1980s and splits
into three directions:
the study of
 analytic properties of growth series  $\Gamma(z)=\sum_{n=0}^{\infty}\gamma(n)z^n$,
the study  of groups of intermediate growth,
and the study  around  Gromov's problem  on  the existence of groups of exponential but not uniformly
exponential growth \cite{gromov_panlaf:book81}.
This paper describes  the main  developments of the second and
third  directions in detail, while the direction of study of growth series is
only briefly mentioned towards the end of Section~\ref{prelim}.

It was discovered that for many groups and even for classes of groups, the function represented by the growth series
$\Gamma(z)$ is rational.  For instance, this holds for virtually abelian groups
\cite{klarner:rational1,klarner:rational2,benson:growth83}, for Gromov hyperbolic groups
(with any finite system of generators) \cite{gromov:hyperbolic,cdp90}, and for Coxeter groups with the canonical
system of generators \cite{paris:growth91}. But the growth series can be irrational even in the case of
 nilpotent groups of nilpotency degree 2 (for instance, for some higher rank Heisenberg groups),
and it can be rational for one system of generators and
 irrational for another \cite{stoll:trans96}. A very interesting approach to the study of growth functions,
 together with applications, is suggested by  Kyoji Saito \cite{saito:limit10,saito:series11}.

All known examples of groups of intermediate growth are groups of branch type and are self-similar groups, or have
 some self-similarity features  (for example, groups from the class $\mathcal{G}_{\omega}, \omega \in \Omega$ discussed
 in   section~\ref{inter1}), or arise from constructions  based on the use of branch groups of self-similar
 type. These groups act on spherically homogeneous rooted trees, and can be studied using techniques which
have been developed during the last three decades.

 Until recently, there were no examples of groups of intermediate growth with a precise estimate on the growth rate
 (there were only upper and lower bounds on the growth as discussed in  detail in section~\ref{inter2}).
 Results about oscillatory behavior of the growth function were presented in \cite{grigorch:degrees} and \cite{grigorch:habil} in
 connection with the existence of groups with incomparable growth.
Recent articles of L.~Bartholdi and A.~Erschler \cite{bartholdi_Ers:growth10,bartholdi_ersch:given(2)11}
provide examples of groups with explicitly computed intermediate behavior of the type $\exp{n^{\alpha}}, \alpha <1$ and many  other
 types of growth, as given in Theorem~\ref{ersch}.

   During the last three decades, many remarkable properties of groups of intermediate growth
  were found, and some of these properties are listed below.  But unfortunately,  we do not know yet if there are finitely presented
  groups of intermediate growth, and this is the main open problem in the field.

\medskip

 The paper is organized as follows.  After providing some background about  group growth in section~\ref{prelim},
 and formulation of Milnor's Problem in section~\ref{problem},  we
 discuss  geometric motivations for the group growth and  the relation of growth with amenability (sections~\ref{geom} and \ref{amenab}).
Section~\ref{1981} describes results obtained before 1981  and section~\ref{polynom0}
 contains an  account of results about the polynomial growth case.

  We describe the main construction of groups of intermediate growth,  methods of
 obtaining upper and lower bounds, and the Bartholdi-Erschler construction in sections~\ref{inter1}, \ref{inter2} and~\ref{inter3}).
The Gap Conjecture is the subject of section ~\ref{gapconj}. In section~\ref{probabil}, we discuss the
relation between the growth of groups and asymptotic characteristics of random walks on them. The final
section~\ref{miscellaneous} contains a short discussion concerning  Gromov's Problem on groups of uniformly
exponential growth,
 on the oscillation phenomenon that holds in the intermediate growth case, and on the role of just-infinite
 groups in the study of group growth.

A number of open problems are included in the text. In our opinion they are among the most important problems
in the field of  group growth. We hope that these problems  will stimulate  further studies of group growth
and related topics.


\medskip
During the last two decades, a number of nice expositions of various topics related to the growth of groups
and their application have been published.  These include the following sources
\cite{wagon:banach93,navas:book,grigorch_h:on_problems,ceccher_machi_scar:grigorch,grigorch-pak:intermediate_growth}.
We especially recommend the book of P. de~la~Harpe \cite{harpe:topics},  which is a comprehensive source of
information about finitely generated groups in general and, in particular, about their growth, as well as the
recent books of T.~Ceccherini-Silberstein and M.~Coornaert \cite{ceccherini_coorn:book10} and A.~Mann
\cite{mann:growthbook} which contains additional material.

The theory of growth of groups is a part of a bigger area of
mathematics  that studies coarse asymptotic properties of various
algebraic and geometric objects.
Some views in this direction can be found in  the following sources:
\cite{lubotzky:discrete_groups-book94,bridson_haefl:book99,roe:coarse03,bekka_harpe:book08,nowak_yu:book12}.

The above sources are recommended to the reader for an introduction to the subject,  to learn more details
about some of the topics considered in this article, or to get information about some other directions of
research involving the growth of groups.

\section{Acknowledgments} The author would like to express his thanks to L.~Bartholdi, A.~Bonifant,
T. Chec\-cher\-ini-Silberstein, T.~Delzant, A.~Erschler, V.~Kaimanovich, A.~Mann, T.~Nagnibeda, P.~de~la~Harpe,
I.~Pak, M.~Sapir, and J.S.~Wilson. My special thanks to S.~Sutherland and Z.~Sunic for tremendous help with
preparation of this paper and for valuable remarks, comments and suggestions. Parts of this work were
completed during a visit to the Institut Mittag-Leffler (Djursholm, Sweden) during the program Geometric and
Analytic Aspects of Group Theory.

\section{Preliminary facts}\label{prelim}

Let $G$ be a finitely generated group with a system of generators $A\!=\!\{a_1,a_2,\dots,a_m\}$ (throughout the
paper we consider only infinite finitely generated groups and only finite systems of generators). The
\emph{length} $|g|=|g|_A$ of the element $g\in G$ with respect to $A$ is the length $n$ of the shortest
presentation of $g$ in the form
\[g=a^{\pm 1}_{i_1}a^{\pm 1}_{i_2}\dots a^{\pm 1}_{i_n},\]
where $a_{i_j}$ are elements in $A$. This depends on the set of generators, but for any two systems of generators $A$
and $B$ there is a constant $ C \in \mathbb{N}$ such that the inequalities
\begin{equation} \label{length} |g|_A \leq C|g|_B,  \qquad |g|_B \leq C|g|_A.
\end{equation}
hold. To justify this it is enough   to express each $a$-generator as a word in $b$-generators and vice
versa: $a_i=A_i(b_{\mu}), b_j=B_j(a_{\nu})$. Then $C=max_{i,j}\{|A|_i,|B|_j\}$, where $|W|$ denotes the
length of the word $W$.  In addition to the length, one can also introduce a \emph{word metric}
$d_A(g,h)=|g^{-1}h|, g,h\in G$ on a group  which is left invariant. More general type of left-invariant
metrics and length functions on groups can be considered and studied as well. For instance, one can assign
different positive weights to generators and define the length and metric according to these weights.

The \emph{growth function} of a  group $G$ with respect to the generating set $A$ is the function
\[\gamma_G^A(n)=\bigl|\{g\in G : |g|_A\leq n\}\bigr|,\]
where $|E|$ denotes the cardinality of the set $E$, and $n$ is a natural number.

If $\Gamma=\Gamma(G,A)$ is the Cayley graph of a group $G$ with respect to a generating set $A$, then $|g|$
is the combinatorial distance from vertex $g$ to  vertex $e$ (represented by  identity element $e\in G$), and
$\gamma_G^A(n)$ counts the  number of vertices at combinatorial distance $\leq n$
 from  $e$ (i.e., it counts the number of elements in the
ball of radius $n$ with  center at the identity element).

It follows from (\ref{length}) that the growth functions $\gamma_G^A(n), \gamma_G^B(n)$ satisfy  the
inequalities
\begin{equation}\label{growth} \gamma_G^A(n) \leq  \gamma_G^B(Cn), \qquad
\gamma_G^B(n) \leq \gamma_G^A(Cn).
\end{equation}

The dependence of the growth function on the  generating set is an inconvenience and it is customary to avoid
it by using the following trick. Following J.~Milnor,  two functions on the naturals $\gamma_1(n)$ and $\gamma_2(n)$
are called \emph{equivalent} (written $\gamma_1(n)\sim \gamma_2(n)$)
 if there is a constant
$C\in \mathbb{N}$ such that $ \gamma_1(n) \leq  \gamma_2(Cn)$, $\gamma_2(n) \leq  \gamma_1(Cn)$ for all
$n\geq 1.$ Then according to (\ref{growth}), the growth functions constructed with respect to two different
systems of generators are equivalent. The class of equivalence $[\gamma_G^A(n)]$ of the growth function  is
called  the \emph{degree of growth}  (the \emph{growth degree}), or  the \emph{rate of growth} of a group
$G$. It is an invariant of a group not only up to   isomorphism but also  up to a weaker equivalence relation
\emph{quasi-isometry}.

\smallskip
Recall that two finitely generated groups $G$ and $H$   with  generating sets $A$ and $B$ respectively
 are quasi-isometric if the metric spaces
$(G,d_A)$ , $(H,d_B)$ are quasi-isometric. Here $d_A$  and $d_B$ are the \emph{word metrics} on $G$ and $H$
defined as $d_A(f,g)=|f^{-1}g|_A, d_B(h,l)=|h^{-1}l|_B$, with $f,g\in G$ and $h,l\in H$. Two metric spaces $(X,d_1),
(Y,d_2)$ are quasi-isometric if there is a map $\phi\: X\rightarrow Y$ and constants $C\geq 1, D \geq 0$ such
that
\[ \frac{1}{C}d_1(x_1,x_2)-D \leq d_2(\phi (x_1), \phi(x_2))\leq C d_1(x_1,x_2)+D \]
for all $x_1,x_2 \in X$; a further requirement is that there is a constant $L >0$ so that for any point $y \in Y$, there
is a point $x\in X$ with
\[d_2(y, \phi(x)) \leq L.\]
This concept is due to M.~Gromov \cite{gromov:hyperbolic} and is one of the most important notions in geometric
group theory, allowing the study of groups from coarse point of view.

\medskip
It is easy to see that the growth of a group coincides with the growth of a subgroup of finite index, and
that the growth of a group is not smaller than the growth of a finitely generated subgroup or a factor group.
We will say that a group is \emph{virtually nilpotent} (resp. virtually solvable) if it contains a nilpotent
(solvable) subgroup of finite index.

We will also consider a preoder $\preceq  $ on the set of growth functions:
\begin{equation}\label{preoder}\gamma_1(n) \preceq \gamma_2(n)
\end{equation}
if there is an integer $C>1$ such that $\gamma_1(n) \leq \gamma_2(Cn)$ for all $n\geq 1$. This makes a  set
$\mathcal{W}$ of growth degrees of finitely generated groups  a partially ordered set. The notation $\prec$ will
be used in this article  to indicate a strict inequality.

Observe that Schwarz  in his note \cite{schvarts:55}  used  formally  weaker equivalence relation $\sim_1$
given by
 inequalities
\begin{equation}\label{growth0} \gamma_G^A(n) \leq  C\gamma_G^B(Cn), \qquad
\gamma_G^B(n) \leq C\gamma_G^A(Cn).
\end{equation}
But indeed both equivalence relations coincide when restricted to the set of growth functions of infinite finitely
generated groups,  as was observed in \cite[Proposition3.1]{grigorch:degrees}. Therefore we will use
either of them depending on the situation.

Because of the independence of the growth rate  on a generating set, we will usually omit subscripts in the notation.

\bigskip
Let us list  the main  examples of  growth rates that will be used later.
\begin{itemize}
\item  The power functions $n^{\alpha}$ belong to different equivalence classes for different $\alpha
    \geq 0$.
\item  The polynomial function $P_d(n)=c_dn^d+ \dots +c_1n+c_0,$ where $c_d \neq  0$ is equivalent to the power function $n^d$.
\item All exponential  functions $\lambda^n, \lambda > 1$ are equivalent and belong to the class $[2^n]$
    (or to  $[e^n].$)
\item All functions of \emph{intermediate type}  $e^{n^{\alpha}}, 0< \alpha <1$ belong to different equivalence classes.

\end{itemize}
Observe that this is not a complete list of rates of growth that a group may have.

\medskip
A free group $F_m$ of rank $m$ has $2m(2m-1)^{n-1}$ elements of length $n$ with respect to any free system
of generators $A$  and
\[\gamma_{F_m}^A(n)=1+2m+2m(2m-1)+\dots +2m(2m-1)^{n-1} \sim 2^n.\]
Since a group with $m$ generators can be presented as a quotient group of a free group of rank $m$, the growth
of a finitely generated group cannot be faster than exponential  (i.e., it can not be superexponential) .
Therefore we can split the growth types into three classes:
\begin{itemize}
\item \emph{Polynomial} growth. A group $G$ has a \emph{polynomial} growth is there are constants $C>0$ and $d> 0$ such that
$\gamma(n) < Cn^d$ for all $ n\geq 1.$ This is equivalent to
\begin{equation}\label{polynom}
\overline{lim}_{n \to \infty}\frac{\log \gamma(n)}{\log n} < \infty.
\end{equation}
A  group $G$  has   \emph{weakly polynomial} growth if for some $d\in (0,\infty)$

\begin{equation}\label{polynom1}
\underline{lim}_{n \to \infty}\frac{\log \gamma(n)}{\log n} < \infty.
\end{equation}
(as we will see later (\ref{polynom}) is equivalent to  (\ref{polynom1})).

\item \emph{Intermediate} growth.   A group $G$ has  \emph{intermediate}  growth if $\gamma(n)$ grows
    faster than  any polynomial but slower than  any exponent function $\lambda^n, \lambda >1$ (i.e.
    $\gamma(n) \prec e^n$).
\item \emph{Exponential growth}.  A group $G$ has  \emph{exponential} growth if  $\gamma(n)$ is
    equivalent to $2^n$.
\end{itemize}

The case of exponential growth can be redefined in the following way.   Because of the obvious
semi-multiplicativity
\[\gamma(m+n)\leq \gamma(m)\gamma(n)\]
 the limit
\begin{equation}\label{growthexp}\lim_{n\to \infty} \sqrt[n]{\gamma(n)}=\kappa
\end{equation}
exists. If $ \kappa > 1$ then the growth is exponential. If $\kappa =1$, then the growth is subexponential and
therefore is either  polynomial or intermediate.

\begin{rem}  In some situations it is reasonable to extend the domain of growth function to all nonnegative real numbers.
 This can be done, for instance, by  setting $\gamma(x)=\gamma([x]),$  where $[x]$ is the integer part
of $x$.

\medskip
Milnor's  equivalence relation~$\sim$  on the set of growth functions is a coarse approach to the study of growth.
Sometimes, in order to study the growth, weaker equivalence relations on the set of monotone functions are
used, for instance when the factor $C$ appears not only in the argument of the growth function but also in
front of it (as in the case of the Schwarz equivalence relation). Additionally, in some situations, an
additive term appears in the form  $C\gamma(Cn)+Cn$ (for instance in the study of Dehn functions of finitely
presented groups, as suggested by Gromov \cite{gromov:hyperbolic}). Instead of the constant $C$ one can use
the equivalence relation with $C=C(n)$, where $C(n)$ is a slowly growing function, for instance a polynomial.
In contrast with the suggested ways of weakening the equivalence relation leading to a loss of some
information about the growth, the more precise evaluation of the rate of growth requires  tools from analysis
based on classical asymptotic methods.

There is a standard way to associate a \emph{growth series} to a growth function, defined by
\[\Gamma(z)= \sum_{n=0}^{\infty}\gamma(n)z^n.\]
The radius of converges of the series is $R=1/\kappa \geq 1$. If the analytic function represented by this
series is  a rational function, or more generally an algebraic function, then  the coefficients $\gamma(n)$
grow either polynomially or exponentially and therefore the intermediate growth is impossible in this case.
For some classes of groups (e.g., abelian groups or Gromov hyperbolic groups), the growth series is always a
rational function \cite{cannon:growth80,gromov:hyperbolic}.  There are groups or classes of groups for which
there is a system of generators with rational or algebraic growth series; for instance,  Coxeter groups have
this property with respect to the canonical system of generators  \cite{paris:growth91}. There are examples
of groups for which the growth series is a rational function for one system of generators but is not rational
for another system of generators \cite{stoll:trans96}. Typically the growth series of a group is a
transcendental function.

\medskip

There is also interest in the study of the \emph{complete} growth function
\[\Gamma_\ast(z)= \sum_{g \in G}gz^{|g|} \in \mathbb{Z}[G][[z]]\]
viewed an element of the ring of formal power series with coefficients in a  group ring $\mathbb{Z}[G]$ of a
group, or an \emph{operator} growth function
\[\Gamma_{\pi}(z)= \sum_{g\in G}\pi(g)z^{|g|},\]
where $\pi(g)$ is a representation of the group in a Hilbert space (or more generally, a Banach space)
\cite{liardet:96,grigorch_nagnib:complete97}. An important case is when $\pi$ is a left regular
representation of the group in $l^2(G)$.
\end{rem}

\begin{rem}  The growth can be defined for many algebraic, combinatorial, geometric, probabilistic and
dynamical objects. For instance, one can speak about the growth of a connected locally finite graph $\Gamma$,
by which we mean the growth of the function $\gamma_{\Gamma,v}(n)$, which counts the number of vertices at
combinatorial distance $\leq n$ from a base vertex $v$. Growth of such graphs can be superexponential, but
if the graph has uniformly bounded degree (for instance is a regular graph), the growth is at most exponential. The
growth does not depend on the choice of $v$.

The growth  of a Riemannian manifold $M$ is, by  definition, the  growth  as $r \to \infty$ of the function
$\gamma(r)=Vol(B_x(r))$ expressing the volume of a ball of radius $r$ with  center at fixed point $x \in M$.
The rate of growth is independent of the choice of $x$.

The definition of  growth of a semigroup is similar to the group case. Also one can define  growth of
finitely generated associative algebras, graded associative algebras, Lie algebras, etc. For instance, if
$\mathcal{A}=\bigoplus_{n=0}^{\infty} \mathcal{A}_n$ is a finitely generated associative graded algebra
defined over a field $\mathbb{F}$ then the growth of dimensions $d_n=\dim_\mathbb{F}\mathcal{A}_n$ determines
the growth of $\mathcal{A}_n$ and the corresponding growth series
\begin{equation}\label{hilbert}
\mathcal{H}(z)=\sum_{n=0}^{\infty}d_nz^n
\end{equation}
usually is called a \emph{ Hilbert-Poincare} series.

\medskip

  Given a countable group $G$ and a probabilistic measure $\mu$ on it  whose support generates $G$, one can consider
a right random walk
  on $G$ which begins at the identity  element $e\in G$ and such that the transitions $g \rightarrow gh$ happen with probability
  $\mu(h)$. One of the main characteristics of such random process  is the probability $P^{(n)}_{e,e}$ of return
  after $n$ steps.  This probability may decay exponentially  (as $r^n, r < 1$) or subexponentially. The value
  \begin{equation}\label{spectral}
  r=\overline{lim}_{n \to \infty}\sqrt[n]{P^{(n)}_{e,e}}
 \end{equation}
is called the \emph{spectral radius} of the random walk (observe that 
 $\lim_{n \to\infty}\sqrt[2n]{P^{(2n)}_{e,e}}$ exists).  It was introduced by H.~Kesten in \cite{kesten:symmetric}. In the
case of a symmetric measure (i.e.\ when  for any $g\in G$, the equality $\mu(g)=\mu(g^{-1})$ holds), the spectral
radius coincides with the norm of the Markov operator of the random walk, and the subexponential decay of
 $P^{(n)}_{e,e}$  (i.e.\ the equality $r=1$)  holds if and only if the group is amenable in the von~Neumann sense
(i.e.\ when $G$ has a
 left invariant mean) \cite{vonNeumann:1929,kesten:amenab}.
We will discuss this topic in more details later in this paper.

\end{rem}

\section{The Problem and the Conjecture of Milnor}\label{problem}

In his note  published in the American Mathematical Monthly \cite{milnor:problem},  Milnor formulated a remarkable
problem concerning the growth of groups, as well as an ingenious conjecture concerning the polynomial growth case.
We reproduce them here, but dividing the  problem into two parts (which was not done in
\cite{milnor:problem}). As before, $\gamma(n)$ denotes the growth function of a finitely generated group (we
shall keep this notation through  the paper).

\bigskip

\theoremstyle{plain}
\newtheorem*{MilnorProblem}{Milnor's Problem}
\begin{MilnorProblem}\quad\\[-\baselineskip]
\begin{itemize}
 \item[(I)] Is the function $\gamma(n)$ necessarily equivalent either to a power of $n$ or to  exponential
function $2^n$?

 \item[(II)] In particular, is the growth exponent
\begin{equation}\label{polynom2} d=\lim_{n \to \infty} \frac{\log \gamma(n)}{\log n}
\end{equation}
always either a well-defined integer or infinity? For which groups is $d < \infty$?
\end{itemize}
\end{MilnorProblem}
\medskip

\theoremstyle{plain}
\newtheorem*{MilnorConjecture}{Milnor's Conjecture}
\begin{MilnorConjecture}
 (A possible conjecture would be that $d < \infty$ if and only if $G$ contains a nilpotent
subgroup of finite index). 
\end{MilnorConjecture}
\bigskip

The first part of Milnor's Problem is a question on the existence of groups of intermediate growth. The
second part and the conjecture are oriented toward the study of the polynomial growth case and the
regularity properties of growth functions.

Clearly the motivation for  suggesting such a  problem and  a conjecture was based on Milnor's background
in the area of group growth as of 1968.  Below we will provide  more information on what was known
about group growth around 1968.

In short, the history of solutions of Milnor's Problem and his Conjecture is the following. The Conjecture was
confirmed by M.~Gromov \cite{gromov:poly_growth}.  This, together with the results of Guivarc'h
\cite{guivarch:70,guivarch:growth71}, B.~Hartley (unpublished, but see ``Added in Proof'' in
\cite{wolf:growth68}), and  H.~Bass \cite{bass:polynom72}, showed that the upper limit (\ref{polynom}) is a
non-negative integer and that the finiteness of the limit (\ref{polynom2}) implies that the group is virtually
nilpotent. The existence of the limit (\ref{polynom2}) also follows from results of the mathematicians quoted above,
giving the complete solution of the second part of Milnor's Problem. In the case of nilpotent groups,
the existence of the limit
\begin{equation}\lim_{n \to \infty}\frac{\gamma(n)}{n^d}
\end{equation}
(where $d$ is the degree of polynomial  growth) was proved  by P.~Pansu \cite{pansu:polynom83} and is an
additional bonus in the study of the polynomial growth case.

The first part of Milnor's Problem was solved in the negative  by the author in 1983
\cite{grigorch:milnor,grigorch:degrees,grigorch:degrees85}. It was shown that there are groups of
intermediate growth, and moreover that there are uncountably
many of them. Also, it was shown that there are pairs of  groups  with incomparable growth in the sense of the
order $\preceq$. In addition, many other results about  growth and algebraic properties of groups of intermediate
growth were obtained around 1984 and later. Despite a negative solution of the first part of Milnor's problem,
the fact that there are groups of intermediate growth made group theory richer and  substantially extended
the area of its applications.

All known examples of groups of intermediate growth are infinitely  presented groups.

\begin{prob} \label{finitepresent}
Is it true that the growth function of a finitely presented group is equivalent either to a
polynomial or to the exponential function $2^n$?
\end{prob}

In  \cite{grigorch-pak:intermediate_growth}, it is conjectured that there are no finitely presented  groups
of intermediate growth. At the same time, the author suggests even a stronger conjecture.

\begin{conj} A finitely presented group either contains a free subsemigroup on two generators
 or is virtually nilpotent.
\end{conj}

Problem~\ref{finitepresent} is the main remaining open problem concerning group growth. As we will see
later, there are many recursively presented groups of intermediate growth, in particular, the groups
$\mathcal G = \mathcal{G}_{\xi}$ and $\mathcal{G}_{\eta}$ described in section~\ref{inter1} are recursively
presented; we note that the group $\mathcal G$ will serve in this text as the main illustrating example.
By Higman's embedding theorem \cite{lyndon_s:cgt_book} such groups embed into finitely presented groups. It
is worth mentioning that for the group $\mathcal{G}_{\xi}$, there is a very precise and nice embedding based
on the use of Lysionok presentation (\ref{lysenok}), as was observed in \cite{grigorch:example} (a similar
claim holds for $\mathcal{G}_{\eta}$). In fact, a similar embedding exists for all groups with finite
$L$-presentation that are defined and studied in \cite{bar_gs:branch,bartholdi:endomorphic_presentations}.
The corresponding finitely presented group (which we denote here $\tilde{\mathcal{G}_{\xi}}$) is an ascending
$HNN-$extension of $\mathcal{G}_{\xi}$. It has a normal subgroup $N$ which is ascending union of conjugates
of a subgroup isomorphic to $\mathcal{G}_{\xi}$ with the quotient $\tilde{\mathcal{G}_{\xi}}/N$ isomorphic to
a infinite cyclic group. Unfortunately (or fortunately), the group $\tilde{\mathcal{G}_{\xi}}$ has
exponential growth but shares the property of amenability with $\mathcal{G}_{\xi}$.  The quotients of this
group are described in \cite{sapir_wise:scending02}. The idea of finding a finitely presented group
containing $\mathcal{G}_{\xi}$ as a normal subgroup (with the hope of thus obtaining a finitely presented
group of intermediate growth) fails, as was observed by M.~Sapir, because $\tilde{\mathcal{G}_{\xi}}$  cannot
serve as a normal subgroup of any finitely presented groups
 (see the argument in
\url{http://mathoverflow.net/questions/73076/higman-embedding-theorem/}).

 Moreover, it was observed by
P. de la Harpe and the author that any finitely presented group $\hat{\mathcal G}$ that can be
homomorphically mapped onto $\mathcal{G}$  contains a free subgroup on two generators and hence is of
exponential growth. In view of these facts it would be interesting to find a finitely presented group with a
normal subgroup of intermediate growth and to find a finitely presented group without a free subgroup on two
generators (or perhaps even a finitely presented amenable group) that can be mapped onto a group of
intermediate growth.
This is discussed in detail in \cite{benli_grigorch_harpe:13}.

For the group $\mathcal{G}_{\eta}$ and many other groups possessing  a presentation of the type
(\ref{lysenok}) (i.e.\ a presentation involving a finite set of relators and their iterations by one or more,
but finitely many, substitutions), embeddings into finitely presented groups similar to the one for
$\mathcal{G}$ also exist~\cite{bar_gs:branch,bartholdi:endomorphic_presentations}.

\section{Relations between group growth and Riemannian  geometry}\label{geom}

One of the first results showing the usefulness of the notion of group growth was the result of
 A.S.~Schwarz,  who  proved the following theorem in 1957.

\begin{thm}\label{schvarts}
Let $\tilde{M}$ be the universal cover of a compact Riemannian manifold $M$.  Then the rate of
growth of $\tilde{M}$ is equal to the rate of growth of the fundamental group $\pi_1(M).$
\end{thm}

In addition to the obvious examples of  groups with polynomial or exponential growth (free abelian groups
and free noncommutative groups respectively), in his note Schwarz produced   examples of solvable (and even
metabelian) groups of exponential growth. These are defined as a semidirect product of $\mathbb{Z}^d, d\geq 2$ and a
cyclic group generated by an  automorphism $\varphi \in SL_d(\mathbb{Z})$ given by a matrix with at least one
eigenvalue off the unit circle.

In  theorem~\ref{schvarts}, the comparison of the growth rates of a manifold and a group is considered with
respect to the equivalence relation
\begin{equation} \gamma_1(n)  \sim_1 \gamma_2(n)
\quad\Leftrightarrow\quad
\exists C > 0  \  \forall n  \ \  \gamma_1(n) \geq
C\gamma_2(Cn) \   \& \  \gamma_2(n) \geq C\gamma_1(Cn);
\end{equation}
the domain of $\gamma_{\pi_1(M)}(n)$ is extended in the natural way  to $\mathbb{R}_+$. Of course, the
equivalence relations $\sim $ and $\sim_1$ defined on the set of monotone functions of natural argument are
different.  But, as was already mentioned in the previous section,  they  coincide on the set of growth
functions of finitely-generated  infinite groups.

The theorem of Schwarz  relates the growth of groups with the volume  growth of the universal cover of a
compact Riemannian manifold.  In fact a more general statement which deals with non-universal and even non-regular
coverings holds; we formulate this below.

\medskip
Milnor's investigation into the relation between  growth and curvature  led  him  to the following
two results \cite{milnor:note68}.

\begin{thm}\label{miln1}
 If $M$ is a complete $d$-dimensional Riemannian manifold whose
mean curvature tensor $R_{ij}$  is everywhere positive semi-definite, then the growth
function $\gamma(n)$  associated with any finitely generated subgroup of the fundamental
group $\pi_1(M)$  must satisfy
$\gamma(n) < Cn^d$ for some positive constant $C$ .
\end{thm}

\begin{thm}\label{miln2}
If $M$ is a compact Riemannian manifold with all sectional curvatures less than zero, then the growth
function of the fundamental group $\pi_1(M)$ is exponential: $\gamma(n) > a^n$ for some constant $a > 1$.

\end{thm}

The next  example  considered by Milnor  was the first step in the direction toward understanding  the growth
of nilpotent groups. Let $G$ be the nilpotent Lie group consisting of all $3 \times 3$ triangular real
matrices with $1$'s on the diagonal, and let $\mathcal{H}_3$ be the subgroup consisting of all integer
matrices of the same form. Then the coset space $G/\mathcal{H}_3$ is a compact $3$-dimensional manifold with
fundamental group $\mathcal{H}_3$.
\begin{lem} The growth function of $\mathcal{H}_3$ is quartic:
\[C_1n^4 < \gamma_{\mathcal{H}_3}(n) < C_2n^4\]
with $0 <C_1 <  C_2$.

\end{lem}

\begin{cor}
No Riemannian metric on $G/\mathcal{H}_3$ can satisfy either the hypothesis
of Theorem~\ref{miln1} or the hypothesis of Theorem~\ref{miln2}.
\end{cor}

The ideas of Milnor were used by A.~Aves \cite{avez:riemann70} to get a partial  answer to a conjecture of
E.~Hopf \cite{hopf:47}: Let M be a compact, connected Riemannian manifold without focal points. Then either
the fundamental group $\pi_1(M)$ has exponential growth, or M is flat.

\medskip
Now let us go back to the theorem of Schwarz and present a more general statement. But before that, we need to
recall some notions from geometric group theory. Let $G$ be a finitely generated  group  with a system of
generators $A$ and let $H < G$ be a subgroup. Let $\Gamma= \Gamma(G,H,A)$ be the Schreier graph determined by the
triple $(G,H,A)$. The vertices of $\Gamma$ are in bijection with cosets $gH, g\in G$, and two vertices $gH$
and $hH$ are joined by oriented edge (labeled by $a \in A$) if $hH=agH.$ This  notion is a generalization of
the notion of a Cayley graph of a group (Cayley graphs correspond to the case when $H=\{e\}$).
$\Gamma(G,H,A)$ is a $2m$-regular graph, where $m$ is the cardinality of the generating set. As was already
defined,  the growth function $\gamma_{\Gamma}(n)$ of a graph counts the number of vertices at combinatorial
distance $\leq n$ from the base vertex $v$  (we remind the reader that the choice of $v$ does not play a
role; it is natural for a Schreier graph to choose $v=1H$).

If a group $G$ with a generating set $A$ acts transitively on a set $X$ then $\Gamma(G,H,A)$ is isomorphic to
the graph of the action, i.e.\ the graph $\Gamma_{\ast}$ with set of vertices $X$ and set of edges consisting
of pairs $(x,a(x)), x\in X, a\in A$. In this case the subgroup $H$  coincides with the stabilizer $st_G(x)$.

\begin{thm} Let $M$ be a compact  Riemannian manifold, let $H$ be a subgroup of the fundamental
group $\pi_1(M)$ and let $\tilde{M}$ be a cover of $M$ corresponding to $H$  supplied by a  Riemannian metric
lifted from $M$. Then the growth function $\gamma_{\tilde{M}}(r)$ is $\sim_1$ equivalent to the growth
function $\gamma_{\Gamma(\pi_1(M),H,A)}(r)$ (with  domain naturally extended to $\mathbb{R}_+$) of the
Schreier graph $\Gamma(\pi_1(M),H,A)$ where  $A$ is a finite  system of generators of $\pi_1(M)$.

\end{thm}
\begin{proof} Triangulate $M$, lift  the triangulation to $\tilde{M}$, and make the comparison of the volume growth of $\tilde{M}$ with the growth
function $\gamma_{\Gamma(\pi_1(M),H,A)}(r)$ using the graph of the action of  $\pi_1(M)$ on the preimage
$p^{-1}(x)$, where $p\: \tilde{M} \rightarrow  M$ is  the canonical projection, and then apply the arguments from \cite{schvarts:55}.
\end{proof}

A particular case of this theorem  mentioned in \cite{grigorch:topological89} is the case of a regular (i.e.\
Galois) cover,  i.e.\ when $H$ is a normal subgroup of $\pi_1(M)$.  In this case, the growth of the covering
manifold $\tilde{M}$ coincides with the growth of the quotient group $\pi_1(M)/H$ and the latter is
isomorphic to the group of deck transformations of the cover. This fact together with the results about
groups of intermediate growth obtained in \cite{grigorch:degrees} allowed the author to construct uncountably
many Riemannian surfaces supplied with  groups of isometries acting on them cocompactly, having a topological
type of oriented surface of infinite genus with one end, and which are not pairwise quasi-isometric
\cite{grigorch:topological89}.

\section{Results about group growth obtained  before 1981}\label{1981}
Among the first publications to use of the notion of group growth were \cite{schvarts:55} and
\cite{adelson:banach57}. The result of the article of Adelson-Velskii  and Shreider  relates  growth with
amenability and will be discussed in the  next section. Then after more than a decade of sporadic appearances
of group growth in various articles (partly listed in the introduction), the papers of Milnor
\cite{milnor:note68,milnor:solv68} and his note \cite{milnor:problem} appeared.

These publications,  followed by Wolf's article \cite{wolf:growth68}, attracted attention of many
researchers and a flurry of activity occurred during a short period. The fact that nilpotent groups  have
polynomial growth was already observed in \cite{milnor:note68} (the example of Heisenberg group) and in full
generality was studied by Wolf \cite{wolf:growth68} and  B.~Hartly  (see ``Added in Proof'' in
\cite{wolf:growth68}), Guivarc'h  \cite{guivarch:70,guivarch:growth71} and by H.~Bass \cite{bass:polynom72},
who  showed that for nilpotent groups, the growth function satisfies the inequalities
\begin{equation}\label{polyn3}
C_1n^d\leq \gamma(n) \leq  C_2n^d.
\end{equation}
Here
\begin{equation}\label{polyn5}
d=\sum_{i}i \cdot rank_{\mathbb{Q}}(\gamma_i(G)/\gamma_{i+1}(G)),
\end{equation}
where $\gamma_i(G)$ is $i$th member of the lower central series of the group, $C_1$ and $C_2$ are  positive
constants, and $rank_{\mathbb{Q}}(A)$ is the torsion free rank of the abelian group $A$.

\smallskip
In \cite{wolf:growth68}, Wolf proved that a polycyclic group either contains a nilpotent subgroup of finite
index and has polynomial growth or the growth is exponential. At the same time, Milnor observed that a solvable
but  not  polyciclic group  has exponential growth. This,  in   combination with Wolf's result, led to the fact
that the growth of a solvable group is exponential except for the case when the group is  virtually nilpotent.

A self-contained proof of the main result about growth of solvable groups was given by Tits in the appendix
to Gromov's paper \cite{gromov:poly_growth}.  In \cite{rosenblatt:supermen74} J.~Rosenblatt showed that groups
of subexponential growth are \emph{superamenable} (the class of superamenable groups is a subclass of the
class of amenable groups, for the definition see  the next section), and indicated that a solvable group of
exponential growth contains a \emph{free subsemigroup} on two generators (perhaps the latter  was known
before, but the author has no corresponding reference).

The remarkable theorem of Tits (usually called the \emph{\windex{Tits Alternative}})  implies that a finitely generated subgroup of a linear
group (i.e.\ a subgroup of
  $GL_n(\mathbb{F}), n\geq 1$,  $\mathbb{F}$ a  field),  either contains a free subgroup on two generators or is virtually
solvable \cite{tits:alternative72}. This, together with the  results about the growth of solvable groups, imply
that the growth of a finitely generated linear group is either polynomial or exponential.

The first part of Milnor's problem on growth  was  included by S.I.~Adian in his monograph
\cite{adian:b-burnside},
  dedicated to the one of the most famous problems in Algebra -  the \emph{Burnside
Problem on periodic groups}. Adian showed that  the free Burnside group
\begin{equation}
B(m,n)=\langle a_1,a_2,\dots,a_m | X^n =1 \rangle,
\end{equation}
 of exponent $n\geq 665$  ($n$ odd ) with $m\geq 2$ has exponential growth. This is  a stronger result than
the result of P.S.~Novikov and
 S.I.~Adian \cite{novikov_ad:burnside2} about the infiniteness of
 $B(m,n)$ in the case $m\geq 2$ and odd $n \geq 4381$.
 A number of results about growth of semigroups was obtained by V.~Trofimov \cite{trofimov:semigr80} (for more recent developments  about  growth
 of semigroups see \cite{bartholdi_rs:intermediate,shneerson:intermediate04,shneerson:types05}).  However, we are not going to
 get much into the details of growth in the semigroup case.

\medskip
A useful fact about groups of intermediate growth  is due to S.~Rosset \cite{rosset:76}.
\begin{thm} If $G$ is a finitely generated group which does not grow exponentially and
 $H$ is a normal subgroup such that $G/H$ is solvable, then $H$ is finitely generated.
\end{thm}

The basic tool for the proof of this theorem is   \emph{Milnor's lemma}.
\begin{lem}If $G$ is a finitely generated  group with subexponential growth, and if $x, y \in G$,
 then the group generated by the set of conjugates $y, xyx^{-1}, x^2yx^{-2}, \dots $ is finitely generated.
\end{lem}
The theorem of Rosset was stated in 1976 for groups of subexponential growth  but  its real application is
to the case of groups of intermediate growth, because by a theorem of Gromov  the groups of polynomial growth are
virtually nilpotent and all subgroups in such groups are finitely generated.

\medskip
The theorem of Rosset can be generalized. In the next statement we use the notion of elementary amenable
group which is defined in the next section. Observe that solvable groups constitute a subclass of the class
of elementary amenable groups.

\begin{thm} Let $G$ be a finitely generated group with no free subsemigroups on two generators and let the
quotient $G/N$ be an elementary amenable group.  Then the kernel $N$ is a finitely generated group.
\end{thm}

The proof of this fact  is based on the use of a version of Milnor's lemma
(Lem.~1 from
\cite{rhemtulla_long:freesemig95}) and transfinite induction on the ``complexity'' of elementary amenable
groups defined below. Observe that this induction  was used by Chou to prove the absence of groups of
intermediate growth, infinite finitely generated torsion groups, and infinite finitely generated simple
groups in the class of elementary amenable groups
 \cite{chou:eg}.  Interesting results about algebraic properties of ``generalized'' elementary groups are
 obtained by D.~Osin~\cite{osin:elementary02}.

In the introduction to his paper published in 1984, \cite{grigorch:degrees} the author wrote: ``In the past
decade, in group theory there appeared a direction that could be called
`Asymptotic Group Theory' ''
(perhaps this was the first  time when the name ``asymptotic group theory'' was used  in the mathematical
literature). Milnor is one of the pioneers of this direction and his ideas and results contributed a lot to
its formation. Asymptotic group theory studies various asymptotic invariants of groups, first of all
asymptotic characteristics of groups (like growth), many of which were defined and studied during the last
three decades. For instance,  the notion of \emph{cogrowth} (or \emph{relative growth}) was introduced by
author in \cite{grigorchuk:symmetrical00}, and later was used by him, Olshanskii and Adian to answer some
questions related to the
\emph{von Neumann Conjecture} on non-amenable groups \cite{grigorch:invariant78,olshanskii:means81,adian:nanf}. The notion of \emph{subgroup growth} was introduced by Grunewald, Segal and
Smith in \cite{grunewald_seg_smith:88} and studied by many mathematicians (see \cite{lubotzky_segalbook} and
the literature cited there).  Dehn functions and their growth were introduced by Gromov and also happen to be
a popular subject for investigation (see \cite{bridson_haefl:book99} and citations therein).   Many
asymptotic invariants of groups were introduced and studied by Gromov in
\cite{gromov:asymptotic_invariant93}, and many other asymptotic invariants has been introduced since then.
The asymptotic methods in the case of algebras (including the topics of self-similarity) are discussed in
survey of E.~Zelmanov \cite{zelmanov:gaeta06,zelmanov:openproblems07} and in \cite{gromov:ggd08}. The
direction of asymptotic group theory is flourishing at present time and the author has no doubt that
the situation will not change in the next few decades.

\section{Growth and Amenability}\label{amenab}

\begin{defn} [John Von Neumann (1929)]\index{amenable group}
A  group $G$ is called {\em amenable}, if there is a finitely additive measure $\mu$ defined on the algebra
of all subsets of $G$ and  such that:
\begin{itemize}
\item $\mu(G)=1$, $0\leq\mu(E)\leq1$, $\forall E\subset G$
\item $\mu$ is left invariant, i.e.\ $\forall E\subset G$:
\[\mu(E)=\mu(gE), \qquad \forall g\in G, \forall E\subset G.\]
\end{itemize}
\end{defn}

The measure $\mu$   determines an invariant, positive, normalized functional $m$ on Banach space
$l^\infty(G)$ defined by
\[m(f)=\int_Gf\,d\mu\]
and  is called a \emph{left invariant mean} ($LIM$).  And vice versa, any invariant, positive,
normalized functional $m$ determines a measure $\mu$ by restriction of its values on the characteristic
functions of subsets.  In a similar way, von Neumann defined amenability of action of a group on a set
\cite{vonNeumann:1929}.

\medskip

The simplest example of a non-amenable group is the free group $F_2$ on two generators. The simplest
examples of amenable groups are finite groups and commutative groups.

 In the above definition the
group $G$ is assumed to be a group with discrete topology.  There is a version of this definition due to
 N.N.~Bogolyubov \cite{bogolyubov:amenable39}, M.~Day and others for general topological groups \cite{hewitt_ross:book} but it depends on the choice
 of the space of functions (bounded continuous or bounded uniformly continuous functions) as was observed by P. de la Harpe \cite{harpe:amenab73}.  In the case
 of locally compact groups one can use any of the mentioned spaces to define amenability \cite{hewitt_ross:book,greenleaf:means}. Growth of groups also can
be defined for locally compact compactly generated groups but we will not consider  topological groups case
here. We recommend to the reader the survey \cite{harpe_cg:paradoxical} on amenability and literature there.


\medskip

Following M.~Day \cite{day:amenable}, we use $AG$ to denote the class of amenable groups.
This class is extremely
important for various topics in mathematics. In this section we discuss briefly  this notion because of its
relation   with growth. By the theorem of Adelson-Velskii and Schreider \cite{adelson:banach57},   each finitely generated group
of subexponential growth belongs to the class $AG$. This class contains finite groups and commutative groups
and
 is closed under the following operations:
\begin{enumerate}
\item taking   a \emph{subgroup},
\item taking  a  \emph{quotient group},
\item  \emph{extensions} (i.e.\ an extension of amenable group by amenable group is amenable),
\item  \emph{directed unions}\\ $G_\alpha \in AG, G_\alpha\subset G_{\beta}\ \text{if}\ \alpha<\beta
    \Rightarrow\cup_{\alpha}G_\alpha\in AG$ 
(here  we  assume  that  the  family of  groups   $\{G_{\alpha}\}$  is  a directed family,  i.e.\  for  any
 $\alpha$ and $\beta$  there
 is  $\gamma$  with $G_{\alpha} \leq G_{\gamma}$  and  $G_{\beta}  \leq G_{\gamma}$).

\end{enumerate}
\medskip

Let $EG$ be the class of \emph{elementary} amenable groups. This is the smallest class of groups containing
finite groups, commutative groups which is closed with respect to the  operations (1)-(4). For such groups a kind
of complexity can be defined in the following way, as was suggested by C.~Chou \cite{chou:eg}. For each ordinal
$\alpha$, define a subclass $EG_{\alpha}$ of $EG$. $EG_0$ consists of finite groups and commutative groups. If
 $\alpha$ is a limit ordinal then
\[ EG_{\alpha}=\bigcup_{\beta \preceq \alpha} EG_{\beta}.\]
Further, $EG_{\alpha+1}$  is defined as the set of groups which are extensions of groups from the set $EG_{\alpha}$
by groups from  the same set. The \emph{elementary complexity} of a group $G\in EG$ is the smallest $\alpha$
such that $G \in EG_{\alpha}$.

One more class of subexponentially amenable groups, denoted here $SG$,  was defined
in~\cite{grigorch:example} as the smallest class which contains all groups of subexponential growth and is
closed with respect to the operations (1)-(4). It is contained in the class of ``good'' groups, in the terminology
of M.~Freedman and P.Teichner~\cite{freedman_teich:subexp95}. Observe that the inclusions
$EG\subset SG \subset AG$
 hold. A general concept of elementary classes is developed by D.~Osin~\cite{osin:elementary02}, along with an  indication of
 their relation to the Kurosh-Chernikov classes.  We will mention the class $SG$  again in section
\ref{probabil}.

\smallskip
The question of Day \cite{day:amenable} concerning  coincidence of the classes $AG$ and $EG$ was answered in the
negative by the author as a result of the solution of the first part of  Milnor's problem
\cite{grigorch:degrees}. In fact, each group of intermediate growth belongs to the complement $AG\setminus EG$
as was showed in \cite{chou:eg}. And as there are uncountably many 2-generated groups of intermediate growth,
the cardinality of the set $AG\setminus EG$ is the cardinality of the continuum.

\smallskip
The second question of Day about  coincidence of classes $AG$  and $NF$ (the latter is the class of groups
that do not contain a free subgroup on two generators) was answered in the negative by A.~Olshanskii
\cite{olshanskii:means81} and S.~Adian \cite{adian:nanf}. Sometimes this problem is formulated in the form of
a conjecture, with attribution to von Neumann. A counterexample to a stronger version of the von Neumann
conjecture was constructed in \cite{grigorch:invariant78}, where a subgroup $H$  of a free group $F_2$ with
the property that the action of $G$ on $G/H$ is nonamenable and  some nonzero power of each element in $F_2$
belongs to some conjugate of $H$. All three articles
\cite{grigorch:invariant78,olshanskii:means81,adian:nanf} use the cogrowth criterion of amenability from
\cite{grigorchuk:symmetrical00}.

It is important to mention that both problems of Day had  negative solutions   not only for the class of
finitely generated groups but also  for the class of finitely presented groups, as was established by the
author and Olshanskii and Sapir respectively \cite{grigorch:example,olshanskii-sapir03}. In other words, there
are examples of finitely presented amenable but non-elementary amenable groups, and there are examples of
finitely presented non-amenable groups without a free subgroup on two generators. Therefore the situation with
amenable/non-amenable groups is in a sense better than the situation with groups of intermediate growth where
the existence of finitely presented groups is unknown (see problem  \ref{finitepresent}).

The variety of different (but equivalent) definitions of the notion of amenable group is tremendous. Perhaps
 there is no any other notion in mathematics which  may compete with the notion of amenable
group in the number of different definitions. The criteria of amenability of  Tarski, F\"{o}lner, Kesten,
Reiter and  many others were found during the  eight decades of studies of amenability. Kesten's criterion
expresses amenability in terms of  spectral radius $r$ of  random walk on a group  in the following way: a
group $G$ is amenable if and only if $r=1$ (here we should assume that the measure defining the random walk
is symmetric and its support generate group, more on this in  section~\ref{probabil}).

\medskip
In the final part of his article \cite{milnor:note68}, Milnor raised the following question.

\begin{enumerate}
\item[]{\em Consider  a random walk on the fundamental group of
compact manifold of negative curvature. Is the spectral radius $r$ necessarily less than 1?
}
\end{enumerate}

In view of Kesten's criterion, this is equivalent to the question about non-amenability of the fundamental
group  of a compact manifold of negative curvature. This problem was positively solved by P.~Eberlain
\cite{eberlein:someprop73}. Related results are obtained by  Yau~\cite{yau:harmonic75} (see also
Proposition~3 on page 98 of the book by Gromov, Pansu and Lafontaine~\cite{gromov_panlaf:book81}) and Chen
\cite{chen:negatcurv78}. As observed by T.~Delzant,  the solution of Milnor's problem~\ref{miln2} can also
be deduced from the results and techniques of A.~Avez paper~\cite{avez:riemann70}. A theorem, due to
M.~Anderson, considers the case of non positive curvature (i.e.\ the curvature can be 0)
\cite{anderson:fundgroups87}.  It is proved that amenable subgroup of such a group must be virtually abelian
and must be the fundamental group of a flat totally geodesic submanifold.  In fact, after the work of Gromov
\cite{gromov:hyperbolic} on hyperbolic groups, followed by the books
\cite{ghys_haefharpe:book90,cdp90},  it became clear  that fundamental groups of compact
manifolds of negative curvature are non-elementary hyperbolic. Such groups contain a free subgroup on two
generators and hence are nonamenable.

In 1974,  J.~Rosenblatt introduced the interesting notion of a \emph{superamenable} group and proved that a group
of subexponential growth is superamenable \cite{rosenblatt:supermen74}. A group $G$ is called superamenable
if for any $G$-set $X$ and any nonempty subset $E\subset X$ there is a $G$-invariant finitely additive
measure on the algebra of all subsets of $X$ taking  values in $[0,+\infty)$ and normalized by condition
$\mu(E)=1$. This property implies there is an invariant Radon measure for any cocompact
continuous action on a locally compact Hausdorff space (private communication by N.~Monod).  It is worth
mentioning that amenability of a group can be characterized by the property of a group having an invariant
probability measure for any continuous action on a compact space, because of the  theorem of Bogolyubov-Day
(which generalizes a theorem of Bogolyubov-Krylov) \cite{greenleaf:means}. Superamenability is discussed in
 the book \cite{wagon:banach93}. In the next section we will mention it again.

\section{Polynomial growth}\label{polynom0}
Recall that a group $G$ has polynomial growth if there are constants $C>0$ and $d> 0$ such that
$\gamma(n) < Cn^d$ for all $ n\geq 1$,  i.e.\ (\ref{polynom}) holds. As was discussed   in section~\ref{1981}
a nilpotent group has  polynomial growth and therefore  a virtually nilpotent group also has
polynomial growth. In his remarkable paper \cite{gromov:poly_growth}, Gromov established the converse:
 polynomial growth implies the virtual nilpotence of the group.

\begin{thm}\label{gromov} (Gromov 1981) If a finitely generated group  $G$ has polynomial growth then $G$  contains a
nilpotent subgroup of finite index.
\end{thm}

This  theorem, together with  known information about growth of nilpotent groups, gives a wonderful
characterization of groups of polynomial growth in algebraic terms: a finitely generated group has polynomial
growth if and only if it is virtually nilpotent.

\smallskip

The proof of Gromov's theorem is  very geometric by its nature, but it also uses some fundamental facts from
algebra and analysis.   Roughly speaking, the idea of Gromov consists of considering  the group as a metric
space, and at looking  at this space via a ``macroscope''. The implementation of this idea goes through the
 development of the technique of limits of sequences of metric spaces.
Given a group of polynomial growth $G$, this method assiociates with it
a Lie group $\mathcal{L_G}$  which has finitely many connected components,
and a nontrivial homomorphism $\phi\: G \rightarrow \mathcal{L_G}$. The rest of the proof is a matter of
mathematical culture modulo  known facts from algebra and analysis including
the Jordan theorem, Tits alternative and the Milnor-Wolf theorem on growth of solvable groups.
   The skeleton  of Gromov's proof has been adapted to many other proofs. In particular, it initiated
the study of \emph{asymptotic cones}  of groups
\cite{dries_wilkie:gromov84,gromov:asymptotic_invariant93}.

Combining his theorem with the results of Shub and Franks   \cite{shub:expend70} on  properties of
fundamental group of  a  compact manifold admitting an expanding self-covering, Gromov deduced the following
statement.
\begin{cor} An expanding self-map of a compact manifold is topologically conjugate to an
ultra-nil-endomorphism.
\end{cor}

In fact Gromov obtained a stronger result about  polynomial growth.
\begin{thm}[{Gromov's effective polynomial growth theorem}]\label{gromoveffect}
For any positive integers $d$ and $k$, there exist positive integers $R, N$ and $q$ with the
following property. If a group $G$ with a fixed system of generators satisfies the inequality $\gamma(n)\leq
kn^d$ for $n=1,2,\dots,R$ then $G$ contains a nilpotent subgroup $H$ of index at most $q$ and whose degree of
nilpotence is at most $N$.
\end{thm}

In \cite{mann:ggd07}, A.~Mann made the first step in the direction of getting a concrete effective bound
on nilpotency class in Gromov's Theorem~\ref{gromov} by  showing that if $G$  is a  finitely generated
group of polynomial growth of degree $d$, then $G$  contains a finite index nilpotent subgroup of nilpotence
class at most $\sqrt{2d}$. This is a corollary of the Bass-Guivarch formula (\ref{polyn5}).
Mann  also showed  that $G$ contains a finite-by-nilpotent subgroup of
index at most $g(d)$, the latter function being the maximal order of a finite subgroup of $GL_d(\mathbb Z)$,
which is known to be $2^nn!$ in most cases (e.g. for $d=1,3,5$ and $d \ge 11$), \cite[Theorem
9.8]{mann:growthbook}.

Another fact obtained by Gromov in~\cite{gromov:poly_growth} is the following: ``If an infinitely generated
group $G$ has no torsion and each finitely generated subgroup $G_1 <G$ has polynomial growth of degree at
most $d$, then $G$ contains a nilpotent subgroup  of finite index''.  On the other hand, A.~Mann showed  in
~\cite{mann:ggd07} that the property of an infinitely generated group $G$ being virtually nilpotent and of
finite rank is equivalent to the property that there are some (fixed) constants $C$ and $d$ such that
$\gamma^A_{H}(n)\leq C n^d, n=1,2,\dots$  for every finitely generated subgroup $H < G$ with  finite system
of generators $A$.

\medskip
Gromov in \cite{gromov:poly_growth} raised the following question.

\theoremstyle{plain}
\newtheorem*{GromovProblemI}{Gromov's Problem on Growth (I)}
\begin{GromovProblemI}
What is the dependence of the numbers $R,N$ and $q$ on $d$ and $k$? In particular, does there
exist an effective estimate of these numbers in terms of $d$ and $k$?
\end{GromovProblemI}

The second part of this question was answered by Y.~Shalom and T.~Tao  in their ``Quantitative Gromov
Theorem'' \cite[Theorem 1.8]{shalom_tao:polynom10} (see below for more on results of Shalom and Tao).


\medskip
The theorem of Gromov was  improved  by van den Dries and Wilkie \cite{dries_wilkie:gromov84}. They showed
that polynomial growth takes place under a weaker assumption: that the inequality  $\gamma(n) \leq Cn^d$
holds for infinite number of values of the argument $n$ (i.e.\ that the group has weakly polynomial growth
in our terminology, see (\ref{polynom1})). Van den Dries and Wilkie applied techniques
 of nonstandard analysis and generalized the notion of asymptotic cone using ultrafilters.

The original proof of Gromov (as well as its modification by van den Dries and Wilkie, see also the
paper of Tits \cite{tits:gromov81})  is based on the use of the Gleason, Montgomery, Zippin, Yamabe
structural theory of
locally compact groups \cite{montgomery_zip:topolog55}. For a long time there was a hope to find a proof of
Gromov's theorem which would avoid the use of this remarkable (but technically very complicated) machinery
which describes a relation between
 locally compact groups and Lie groups.
 This goal was achieved
by B.~Kleiner \cite{kleiner:plynom10}  who proved Gromov's theorem using a completely different approach that
involves harmonic functions. The core of  Klener's  arguments is in the new proof of (a slight modification
of) the theorem of Colding-Minicozzi \cite{colding_minic:harmonic97}.

\begin{thm}(Kleiner 2010) Let $G$ be a group of
weakly polynomial growth, and $l \in (0,\infty)$. Then the space of harmonic
functions on $\Gamma$ with polynomial growth at most $l$ is finite dimensional.
\end{thm}

Recall that a function $f$  on a group $G$ is called $\mu$-harmonic if it satisfies $Mf=f$, where $M$ is the
Markov operator of the random walk on $G$ determined by the measure $\mu$.  A function $f$  has at most
polynomial growth of degree $l$ if for some positive constant $C$ the inequality
 $|f(g)|\leq C|g|^l$ holds, for all $g \in G$.

{Gromov's effective polynomial growth theorem}
was improved by Y.~Shalom and T.~Tao in different directions
\cite{shalom_tao:polynom10}. They  showed that for some absolute (and explicit) constant $C$, the following
holds for every finitely generated group $G$, and all $d>0$: If there is some $R_0 > \exp(\exp(Cd^C))$ for
which $\gamma_G(R_0)<R_0^d$, then $G$ has a finite index subgroup which is nilpotent of nilpotency degree
less than $C^d$. In addition, an effective upper  bound on the
 index is provided in \cite{shalom_tao:polynom10} if the word ``nilpotent'' is replaced by ``polycyclic''.

In \cite{shalom_tao:polynom10},  a pair $(G,A)$ (a group $G$ with a finite generating set $A$) is called an
$(R_0,d)$-growth group if
\begin{equation} \gamma_G^A(R_0) \leq R_0^d.
\end{equation}


{\StretchyLines
\begin{thm}[Fully quantitative weak Gromov theorem] Let $d > 0$ and $R_0 > 0$,
and assume that
$R_0 >\exp(\exp(cd^c)\!)$
for some sufficiently large absolute constant $c$. Then every $(R_0, d)$-growth group
has a normal subgroup of index at most $\exp(R_0\exp(\exp(d^c)\!)\!)$
 which is polycyclic.
\end{thm}
}

\begin{cor}[Slightly super-polynomial growth implies virtual nilpotency]
Let $(G, A)$ be a finitely generated group such that
$\gamma_G^A(n) < n^{c(\log \log n)^c}$
for some $n > 1/c$, where $c > 0$ is a sufficiently small absolute constant. Then $G$ is
virtually nilpotent.
\end{cor}

This result  shows that if the growth of a group is less than that of $n^{(\log \log n)^c}$, then it is
polynomial.

\medskip

 Gromov's original proof implies the existence of a function $\upsilon(n)$ growing faster than any
polynomial and such that if $\gamma_G(n)\prec \upsilon(n)$, then the growth of  $G$ is polynomial. Therefore
there is a \emph{gap} in the scale of growth degrees of finitely generated groups.  In fact, as there are
groups with incomparable growth functions, it may happen that there can be many gaps between polynomial and
intermediate growth.

The function $ n^{(\log \log n)^c}$  is the first concrete example of the superpolynomial bound separating
the polynomial growth case from the intermediate one. Finding the border(s) between polynomial and
intermediate growth is one of the main open problems about the growth of groups. We call this the  \emph{Gap
Problem}   and there is  the associated \emph{Gap Conjecture}  which we discuss in more detail in
section~\ref{inter2}.

\smallskip

One more result in the direction of study of polynomial growth is the characterization of finitely generated
cancellative semigroups of polynomial growth  given in \cite{grigorch:semigroups_with_cancel88}. Recall that
a semigroup $S$ is cancellative if the left and right cancellative laws hold: $\forall a,b,c \in S,  ab=ac
\Rightarrow b=c, ba=ca \Rightarrow b=c$.  Cancelation is a necessary condition for  embedding of a semigroup
into a group. In 1957, A.I.~Malcev introduced the notion of a nilpotent semigroup  \cite{malcev:nilpotent53}.
A semigroup is nilpotent if for some $n\geq 1$, it satisfies the identity $X_n=Y_n$, where $X_n,Y_n,
n=0,1,\dots$ are words over the alphabet $x,y,z_1,z_2, \dots,z_n,\dots$ defined inductively:  $X_0=x,Y_0=y$
and
\[X_{n+1}=X_nz_nY_n, \qquad Y_{n+1}=Y_nz_nX_n.\]
A semigroup  is of nilpotency degree $n$ if it satisfies the identity $X_n=Y_n$ but  not the identity
$X_{n-1}=Y_{n-1}$. Malcev proved  that a group $G$ is nilpotent of degree $n$ if and only if it is nilpotent
of degree $n$ as a semigroup.

 We  say that a
subsemigroup $L \leq S$ has a   finite (left)  index  in $S$ if there is a finite subset
$E=\{e_1,\dots,e_k\}\subset S$ such that
\[S= \bigcup_{i=1}^ke_iL.\]

\begin{thm}[\cite{grigorch:semigroups_with_cancel88}] Let $S$ be a cancellative semigroup of polynomial growth of degree $d$ i.e.

\begin{equation}\label{polynom3}
d=\overline{\lim}_{n \to \infty}\frac{\log \gamma_S(n)}{\log n} < \infty.
\end{equation}

Then

(i) $S$ has a group $G=S^{-1}S $ of left quotients which also has  polynomial growth of degree $d$ and
therefore is virtually nilpotent. In particular, $d$ is a nonnegative integer.

(ii)  The semigroup $S$ is virtually nilpotent in Malcev sense.  More precisely, $G$ contains a nilpotent
subgroup $H$ of finite index such that $S_1=H\cap S$ is a nilpotent semigroup of the same degree of
nilpotency as $H$,  and $S_1$ is of finite (left or right) index in $S$.
\end{thm}

The proof of this theorem relies  on the fact that a cancellative semigroup $S$ of subexponential growth
satisfies the Ore condition and hence has a group of left quotients $G$ (this also holds for cancellative
right amenable semigroups). The next step consists in getting of an upper bound of polynomial type for the
growth of $G$, and an application of Gromov's theorem finishes the argument.

\smallskip

 Observe that   much  less is known about growth of nilpotent cancellative semigroups than nilpotent groups.  For instance, it is
unknown if there is a constant $c>0$ such that
\[cn^d\leq \gamma_S(n), n=1,2,\dots\]
in the case when $S$ is a cancellative semigroup of polynomial growth of degree $d$. In the case of a
semigroup without cancellation,  the growth  can be equivalent to  the growth of an arbitrary function
$\gamma(n)$ satisfying reasonable restrictions, as was shown by V.~Trofimov in \cite{trofimov:semigr80}.
Also a
semigroup of polynomial growth need not be virtually nilpotent in Malcev sense.

\medskip

We conclude this section with the following questions.

\begin{prob} Let $S$ be a cancellative semigroup of polynomial growth of degree $d$. Is it true that the
limit

\begin{equation}\label{polynom4} d=\lim_{n \to \infty} \frac{\gamma_S(n)}{ n^d}
\end{equation}
exists?

\end{prob}

\begin{prob} Does there exist a finitely generated cancellative semigroup $S$ of subexponential growth such
that the group $G=S^{-1}S$ of left quotients has exponential growth?

\end{prob}

A positive answer to the previous question would provide an example of a superamenable group of exponential
growth and an answer to the J.~Rosenblatt's question from \cite{rosenblatt:supermen74}. Also observe that in
the case   $\gamma_S(n)\prec e^{\sqrt{n}}$,  the group of quotients $G=S^{-1}S$ has subexponential growth
which justifies the following   question.

\begin{prob}\quad\\[-\baselineskip]
\begin{itemize}
\item[(a)] Does there exist a finitely generated cancellative semigroup $S$ of superpolynomial growth
which is strictly less that the growth of $ e^{\sqrt{n}}$?
\item[(b)] Does there exist a finitely generated cancellative semigroup $S$ of  growth equivalent to the growth of
$e^{\sqrt{n}}$?
\end{itemize}
\end{prob}

 This problem is related to  problem~\ref{gap}  and the  Gap Conjecture discussed in
section~\ref{gapconj}.

\section{Intermediate growth: the construction}\label{inter1}

The first part of Milnor's problem was solved by the author in 1984.

\begin{thm} (\cite{grigorch:degrees})\label{growthmain}\quad\\[-\baselineskip]
\begin{itemize}
    \item There are  finitely generated groups of intermediate growth.

   \item The partially ordered set $\mathcal{W}$  of growth degrees of 2-generated groups contains a chain of
   the  cardinality of the continuum, and contains an antichain of the cardinality of the continuum.

    \item  The previous statement holds for the class of 2-generated  $p$-groups for any prime $p$.
    \end{itemize}
\end{thm}

\begin{cor}
There are uncountably many 2-generated groups, up to quasi-isometry.
\end{cor}

The main example of a group of intermediate growth, which will be  denoted through the whole  paper
by $\mathcal{G}$, is the group defined by the next figure:
\begin{center}\label{pict2}
\vspace{1cm} \includegraphics[width=200pt]{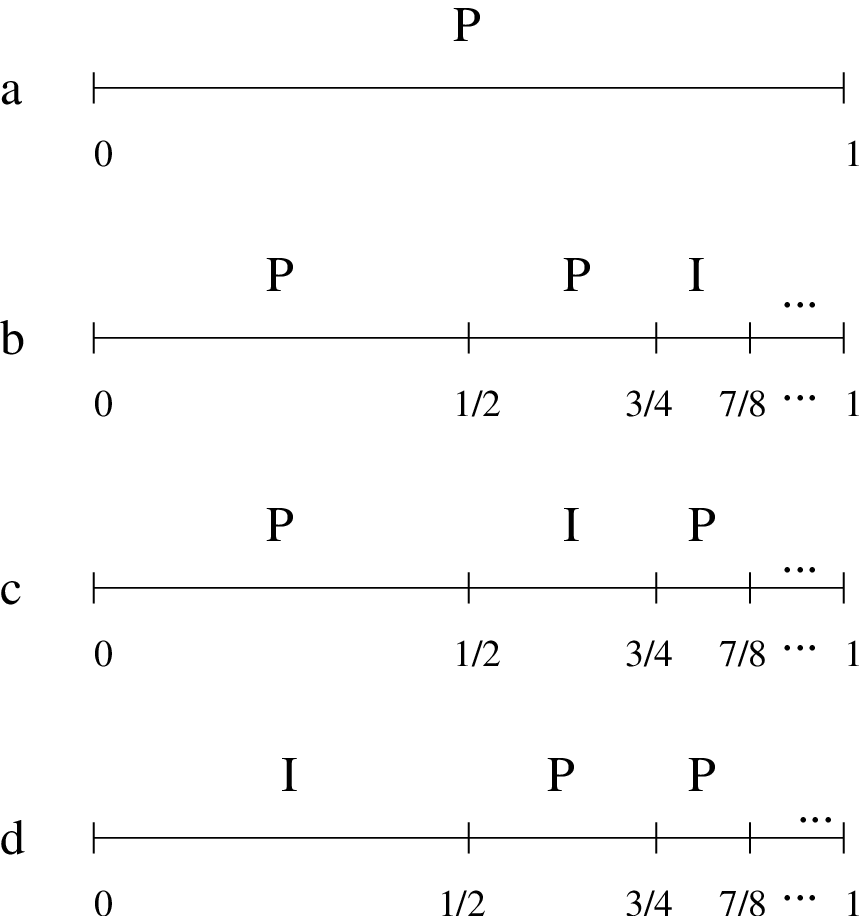}\\
Generators of   $\mathcal{G}$\\[0.2cm]
\end{center}

This group was constructed in 1980 in the author's  note~\cite{grigorch:burnside} as a simple example of a
finitely generated infinite torsion group. Recall that the question on the existence of such groups was the
subject of the \emph{General Burnside Problem} raised by Burnside in 1904. The problem was solved by E.~Golod
in 1964 on the basis of the Golod-Shafarevich Theorem~\cite{golod_s:class_field_tower,golod:p-groups}. The bounded
version of the Burnside problem (when  orders of all elements are assumed to be uniformly bounded) was solved
by S.~Novikov and
S.~Adian~\cite{novikov_a:infinite_periodic1,novikov_a:infinite_periodic2,novikov_a:infinite_periodic3,adian:b-burnside}).
A third version of the Burnside Problem, known as the Restricted Burnside Problem, was solved by E.~Zelmanov
\cite{zelmanov:burnside_2-groups,zelmanov:burnside_odd,zelmanov:icm-restricted}.

The note
\cite{grigorch:burnside} contains two examples of Burnside groups,  and the second group $\mathcal B$
presented there was, in fact, historically constructed before  group $\mathcal G$. Therefore, it can be
considered as the first example  of a nonelementary self-similar group (self-similar groups are discussed in
the next section). The group $\mathcal B$ shares with $\mathcal G$ the property of being an infinite torsion
group, but it is unknown whether it has intermediate growth. Also it is unknown if the Gupta-Sidki $p$-groups
\cite{gupta_s:burnside}, whose construction and properties are similar to $\mathcal B$   in many aspects,
have intermediate growth.

\smallskip
The  group  $\mathcal G$ acts on the unit interval $[0,1]$ with the dyadic rational points removed, i.e.\
on the set $\Delta=[0,1]\setminus\{m/2^n,\  0\leq m \leq 2^n, n=1,2,\dots \}$ (if needed, the action of
$\mathcal G$ can be extended to the entire interval $[0,1)$ ). The generator $a$ is a permutation of the two halves
of the interval, and $b,c,d$ are also interval exchange type transformations (but involve a partition of
$\Delta$ into infinitely many subintervals). They preserve the partition of the interval indicated in the
figure, and their action is described by the labeling of the atoms of  partition by symbols $\{I,P\}$.
The letter $I$ written over an interval means the action is the identity transformation there, and
the letter $P$ indicates the two halves of the
interval are permuted as indicated in the figure below.
\begin{center}\label{pict1}
\includegraphics[width=200pt]{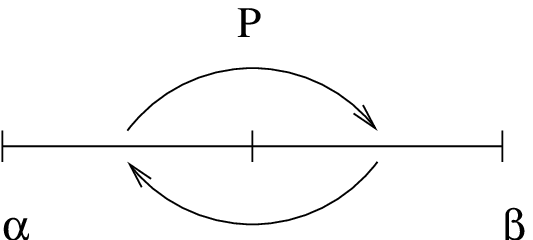}
\end{center}

\medskip

Now we are going   to describe a general construction of uncountably many groups of intermediate growth.
These groups are of the form
\[\mathcal{G}_{\omega}=\langle a,b_{\omega},c_{\omega},d_{\omega}\rangle,\]
\noindent where the generator $a$ is defined as before, and $b_{\omega},c_{\omega},d_{\omega}$ are
transformations of the set $\Delta$  defined using an \emph{oracle} $\omega \in
\Omega=\{0,1,2\}^{\mathbb{N}}$ in the same manner as $b,c,d$. To describe the construction of groups
$\mathcal{G}_{\omega}$ completely,  let us start with the space $\Omega = \{0,1,2\}^{\bb N}$ of infinite
sequences over  alphabet $\{0,1,2\}$. The bijection
\begin{align*}
0 &\longleftrightarrow \left(\begin{matrix} P\\ P\\ I\end{matrix}\right)\\
1 &\longleftrightarrow \left(\begin{matrix} P\\ I\\ P\end{matrix}\right)\\
2 &\longleftrightarrow \left(\begin{matrix} I\\ P\\ P\end{matrix}\right)
\end{align*}
between symbols  of the alphabet $\{0,1,2\}$  and corresponding columns is  used in the construction. Namely
given a sequence $\omega = \omega_1\omega_2\ldots \in\Omega$ replace each $\omega_i$, $i=1,2,\ldots$, by the
corresponding column  to get a vector
\[
\left(\begin{matrix}
U_\omega\\ V_\omega\\ W_\omega\end{matrix}\right)
\]
consisting of three infinite words $U_\omega, V_\omega, W_\omega$ over the alphabet $\{I,P\}$. For instance,
the triple of words
\begin{align*}
U_\xi &= PPI~PPI\ldots\\
V_\xi &= PIP~PIP\ldots\\
W_\xi &= IPP~IPP\ldots,
\end{align*}
corresponds  to the sequence $\xi = 012~012\ldots$, while the triple
\begin{align*}
U_\eta &= PPP~PPP\ldots\\
V_\eta &= PIP~IPI\ldots\\
W_\eta &= IPI~PIP\ldots.
\end{align*}
 corresponds to the sequence $\eta = 01~01\ldots$.

Using  the words $U_\omega, V_\omega, W_\omega$, define transformations $b_\omega$, $c_\omega$, $d_\omega$ of the
set $\Delta$ as in the case of the sequence $\xi$; the sequence $\xi$  corresponds to the group $\mathcal G$. Then
$\mathcal{G}_{\omega}=\langle a,b_{\omega},c_{\omega},d_{\omega}\rangle$. The most interesting groups from
this family correspond to sequences $\xi$ and   $\eta$, which are respectively the groups $\mathcal G$ and
$\mathcal E = \mathcal{G}_{(01)^{\infty}}$ studied by A.~Erschler \cite{erschler:boundary04}.

The generators of $\mathcal{G}_{\omega}$ satisfy the relations
\begin{align*}
&a^2 = b^2_\omega = c^2_\omega = d^2_\omega = 1\\
&b_\omega c_\omega = c_\omega b_\omega = d_\omega\\
&b_\omega d_\omega  = d_\omega b_\omega = c_\omega\\
&c_\omega d_\omega = d_\omega c_\omega = b_\omega
\end{align*}
(this is only a partial list of the relations), from which it follows that the groups $\mathcal{G}_{\omega}$ are
3-generated (and even 2-generated in some degenerate cases, for example, the case of the sequence $0~0\ldots
0\ldots$). Nevertheless, we prefer  to consider them as 4-generated groups, with the system
$\{a,b_\omega,c_\omega,d_\omega\}$  viewed as a canonical system of generators.

While the groups $\mathcal{G}_{\omega}$ are 3-generated, a simple trick used in \cite{grigorch:degrees} produces
a family  $\mathcal{H}_{\omega},\ \omega \in \Omega$ of 2-generated groups with the property that for $\omega
\in \Omega_1$ (the set $\Omega_1$ is defined below), the growth of $\mathcal{H}_{\omega}$ is the same as
growth of $\mathcal{G}_{\omega}$ raised to the fourth power.

\smallskip
Alternatively, the groups  $\mathcal{G}_{\omega}$ can be defined as groups acting by automorphism  on a binary
rooted tree $T$ depicted in the following figure.

\begin{center}\label{binary}
\includegraphics[width=\hsize]{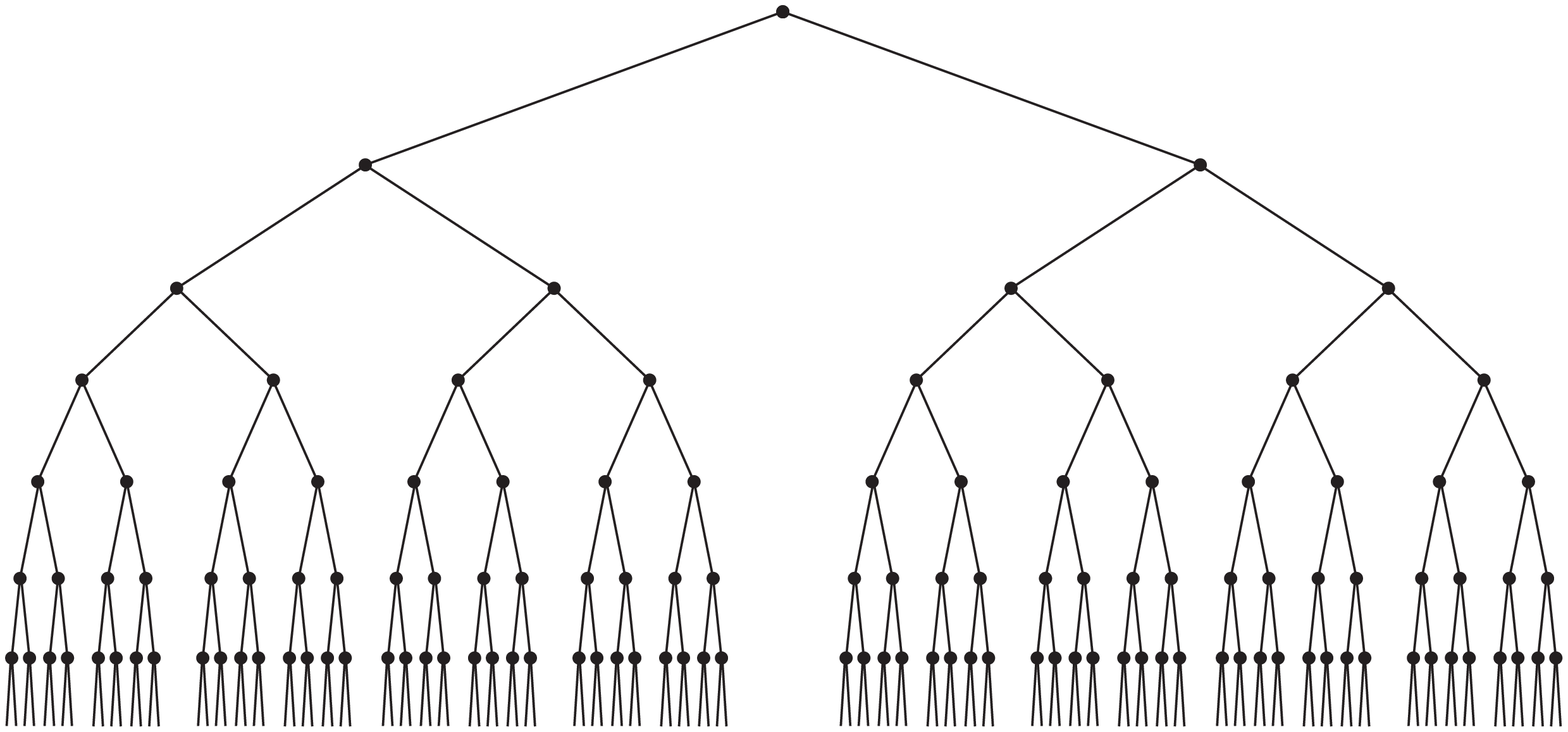}
\end{center}

\medskip
For instance, let us explain how generators of $\mathcal G$ act on $T$. The generator $a$ acts  by
permuting the two rooted subtrees $T_0$ and $T_1$  which have roots at the first level.  The generator $b$  acts on
the left subtree $T_0$ as $a$ acts on $T$ (here we use  the self-similarity property of a binary tree), and acts on
the right subtree $T_1$  as  generator $c$ (whose action on the whole tree is shifted to the right subtree).
Similarly, $c$ acts on $T_0$ as $a$ on $T$ and on $T_1$ as $d$ acts on $T$. And finally, $d$ acts on $T_0$ as
the identity automorphism and on $T_1$ as $b$ on $T$.
This gives us a recursive definition of automorphisms $b,c,d$ as in the figure
below, where dotted arrows show the action of a generator on the corresponding vertex.

\medskip
\begin{center}\label{self-sim}
\includegraphics [width=.95\textwidth]{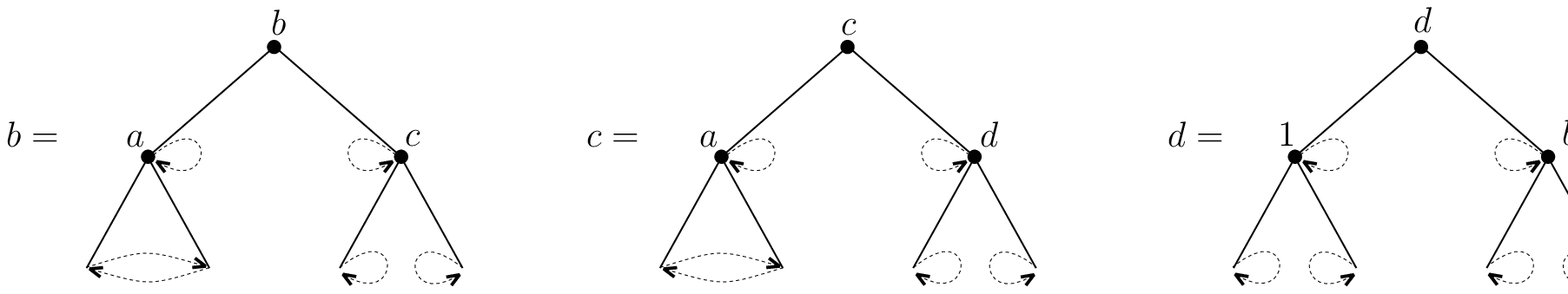}
\end{center}
\medskip

Making the identification of vertices of $T$ with  finite binary sequences, identifying the
boundary $\partial T$ of the tree with  the set of infinite sequences from $ \{0,1\}^{\mathbb{N}}$,  and
representing dyadic irrational points by the  corresponding binary expansions,  we obtain an isomorphism of the
dynamical system $(\mathcal{G},[0,1],m)$ (where $m$ is Lebesgue measure) with the system
$(\mathcal{G},\partial T, \nu)$, where $\nu$ is a uniform measure on the boundary $\partial T$  (i.e.\ a
$\{\frac{1}{2},\frac{1}{2}\}$ Bernoulli measure if the boundary is identified with the space of binary
sequences). A similar description of the action on a binary tree can be given for all groups
$\mathcal{G}_{\omega}$,  although self-similarity should be used in a more general sense.

\smallskip
There are many other ways to define the groups $\mathcal{G}_{\omega}$. For instance, a very nice approach was
suggested by Z.~Sunic which allows one to define the groups $\mathcal{G}$  and $\mathcal{E}$  by use of
irreducible (in the ring $\mathbb{F}_2[x]$) polynomials $x^2+x+1$ and $x^2+1$ respectively.  This approach
allows one to associate to any irreducible polynomial a self-similar group \cite{bar_gs:branch}.

One more
definition of the group $\mathcal{G}$ can be given by the following  presentation found by I.~Lysenok
\cite{lysionok:presentation}:

\begin{equation}\label{lysenok}\mathcal{G}=\langle a,b,c,d \ |\ a^2,b^2,c^2,d^2,bcd, \sigma^k((ad)^4), \sigma^k((adacac)^4), k\geq 0
\rangle,
\end{equation} where

\[\sigma:  \left\{\begin{array}{l}a \rightarrow aca,\\b \rightarrow d,\\c \rightarrow b, \\ d \rightarrow c.
\end{array}\right.\]
 is a  substitution. A similar presentation (called an $L$-presentation) can be found for $\mathcal E$ and
 many other self-similar groups.

\begin{thm}\label{GRI83} (\cite{grigorch:burnside,grigorch:degrees})
\item 1) $\mathcal{G}$ is a residually finite torsion 2-group, which is not finitely presentable.

\item  2) $\mathcal{G}$ is of intermediate growth with bounds
\[e^{n^{1/2}} \preceq \gamma_{\mathcal{G}}(n)\preceq e^{n^{\beta}},\]
where $\beta=\log_{32}31 < 1$.

\item 3) (A.~Erschler \cite{erschler:boundary04}) For any $\epsilon > 0$ the following upper and lower bounds
    hold:
\[ \exp{\frac{n}{\log^{2+\epsilon}n}} \prec \gamma_{\mathcal{E}}(n) \prec
\exp{\frac{n}{\log^{1-\epsilon}n}}.\]
\item  4) (\cite{grigorch:degrees85}) There is a torsion free extension $\hat{\mathcal G}$ of abelian group
    of
    infinite rank by  $\mathcal G$ which has intermediate growth.

\end{thm}

The first part of this theorem shows that there are  examples of  groups of Burnside type already in the
class of groups of intermediate growth.  The last part shows that  there are torsion free groups of
intermediate growth.

The question of finding the precise (in Milnor-Schwarz sense) growth degree of any group of intermediate
growth was open until recently, when in their excellent paper \cite{bartholdi_Ers:growth10},  Bartholdi and
Erschler constructed groups with  growth of the type $\exp n^\alpha$ for infinitely many $\alpha \in
[0.7674,1)$ (more on their result in section~\ref{inter3}). Moreover, they were able to find the growth
degree of the torsion free group  $\hat{\mathcal G}$  by showing that
\[  \gamma_{\hat{\mathcal G} }(n) \sim n^{\log(n)n^\alpha_0},\]
where $\alpha_0\approx 0.7674$,  is defined below~(\ref{barth}).
In addition to the  properties of $\hat{\mathcal G}$ listed above, we have the following fact, which was
established in
the paper of A.~Machi and author \cite{grigorch_machi:93}.

\begin{thm}
The group $\hat{\mathcal G}$  is left orderable and hence acts by orientation preserving  homeomorphisms of
the line $\mathbb R$ (or, equivalently, of the interval $(0, 1)$ ).
\end{thm}

It is unknown if there are \emph{orderable} groups (i.e.\ groups possessing a two-sided invariant order) of
intermediate growth.

{\StretchyLines
\begin{prob}
 Does there exist a finitely generated, orderable group of intermediate growth?
\end{prob}}

A possible approach to this problem could be the following. By a theorem of
M.I.~Zaiceva~\cite{kokorin_kopyt:fully74}, if a factor group $F/A$ of a free group $F$ has an infrainvariant
system  of subgroups (for the definition, see~\cite{kokorin_kopyt:fully74}) with torsion free abelian
quotients, then the group $G/A$ is left orderable. Moreover, by D.M.~Smirnov's theorem
\cite{smirnov:generalized63,kokorin_kopyt:fully74} the group $F/[A,A]$ is totally orderable in this case. It
is known that the group $\hat{\mathcal G}$ has an infrainvariant system of subgroups with torsion free
abelian quotients. Representing the group $\hat{\mathcal G}$ in the form $F/A$ we obtain an interesting
totally orderable group $F/[A,A]$.  The question is whether it is of intermediate growth.

The last theorem shows that $\mathcal G$ embeds into the group $Homeo^+([0, 1])$ of homeomorphisms of the
interval. The notion of growth is useful for study of codimension one foliations, and group growth as well as
growth of Schreier graphs play important role. The question concerning the smoothness of the action  of a
group on a manifold arises in many situations. Many topics related to the subject of actions on the
interval and on the circle can be found in
\cite{ghys:survay01,navas:book,beklaryan:umn04,beklaryan:survay08}. The following result of A.~Navas
\cite{navas:growth08} gives important information about what one can expect in the one-dimensional case.

\begin{thm}
The group $\mathcal G$  embeds into $\text{\itshape Diff\/}^1([0, 1])$. However, for every $\epsilon > 0$,
every subgroup of
$\text{\itshape Diff\/}^{1+\epsilon}([0, 1])$ without a free subsemigroup on two generators is virtually nilpotent.
\end{thm}

\medskip

Until 2004, all known examples of groups of intermediate growth were residually finite groups. The question
about  existence of non residually finite groups of intermediate growth was open for a while. In
\cite{erschler:notresid04}, Erschler constructed an uncountable family of non residually finite groups of
intermediate growth. These groups  are extensions of an elementary 2-group of infinite rank by $\mathcal G$.

Let $\Omega_0\subset \Omega$ be the subset consisting of
sequences in which all three symbols $0,1,2$ occur infinitely often,
let $\Omega_1\subset \Omega$ be the subset consisting of sequences in which at least two of
the symbols $0,1,2$ occur infinitely often, and let $\Theta \subset \Omega$ be the subset consisting of
sequences $\omega$ with the
property that there is a natural number $C=C(\omega)$ such that for every $n$, each of symbols $0,1,2$ appears
among any $C$ consecutive symbols $\omega_n\omega_{n +1}\dots\omega_{n+C-1}$. Observe that the inclusions
$\Theta \subset \Omega_0\subset \Omega_1$ hold.

Recall that a group is called \emph{just-infinite} if it is infinite but every proper quotient is finite.  We
postpone the definition of \emph{branch} groups (as well as the notions of \emph{self-similar} group and
\emph{contracting} group) until section~\ref{inter1}.

\begin{thm}\label{growth2}(\cite{grigorch:degrees})
\begin{itemize}
\item The word problem for the group $\mathcal{G}_{\omega}$ is solvable if and only if the sequence
    $\omega$ is recursive.
 \item For each $\omega \in \Omega_0 $ the group $\mathcal{G}_{\omega}$ is a just-infinite torsion 2-group.
 \item For each $\omega \in \Omega_1 $ the group $\mathcal{G}_{\omega}$ is infinitely presented and  branch.
\item For each  $\omega \in \Omega_1$ the group $\mathcal{G}_{\omega}$ has intermediate growth and the lower
bound

\[e^{n^{1/2}} \preceq \gamma_{G_{\omega}}(n)\]
holds.

\item For each $\theta \in \Theta$ there is a $\beta =\beta(\theta)<1$ such that
\[ \gamma_{G_{\omega}}(n) \preceq  e^{n^{\beta}}\]
holds.
\end{itemize}

\end{thm}

The first statement of this theorem shows that there are many groups of intermediate growth for which the word problem is  decidable, and therefore there are recursively presented groups of intermediate growth.  At the same time
there are recursively presented residually finite groups of intermediate growth for which the  word problem is
unsolvable, as is shown in \cite{grigorch:habil}.

Formally speaking  \cite{grigorch:degrees}   does not contain the proof of the
last part of the last theorem. But it can be easily obtained following the
same line as the proof of theorem 3.2 from \cite{grigorch:degrees}.

\smallskip
The lower bound by the function $e^{\sqrt n}$ for the growth of $\mathcal{G}$ was improved by L.~Barthodi
\cite{bartholdi:lower} and Y.~Leonov \cite{leonov:bound}, who showed that $\exp{n^{0.5157}}\preceq
\gamma_{\mathcal{G}}(n)$.

The upper bound for $\mathcal{G}$  given by  theorem~\ref{GRI83} was  improved  by L.~Bartholdi
\cite{bartholdi:growth} by showing that
\begin{equation}\label{barth}\gamma_{\mathcal{G}}(n) \preceq e^{n^{\alpha_0}},
\end{equation}
with $\alpha_0=\log 2/\log (2/\rho) \approx 0.7674$, where $\rho$ is the real root of the polynomial
 $x^3+x^2+x-2$.  The method used consists of consideration of a more general type of length function, arising from
 prescribing positive weights to generators, and counting the length using the weighted contribution of
each generator.
  The method   relies  on the idea of giving  the optimal
 weights  to  generators that lead to the best upper bound. It happened to be also useful for the
 construction of groups of intermediate growth with  explicitly
 computed growth, as recently was demonstrated by Bartholdi and Erschler in \cite{bartholdi_Ers:growth10}
(more on this in section~\ref{inter3}).

Another technical tool  was  explored   by  R.~Muchnik and I.~Pak \cite{muchnik_p:growth}  to get an upper bound
on growth for the whole family of groups $\{\mathcal{G}_{\omega}\}$. Surprisingly,  in the case of
$\mathcal{G}$ their approach give the same upper bound  as (\ref{barth}), so the question of improving it is
quite intriguing (see problem~\ref{pak}).

 Unfortunately, even after three decades of study of  the group $\mathcal{G}$ and other groups
 $\mathcal{G}_{\omega}$,
 we still do not know the precise growth rate of any group of intermediate growth from the
 family  $\{\mathcal{G}_{\omega}\}$.

\medskip
\begin{prob}\quad\\[-\baselineskip]
\begin{itemize}
\item[(1)] Does there exist $\alpha$ such that $\gamma_{\mathcal{G}}(n)\sim e^{n^{\alpha}}$?

\item[(2)] If the answer to the previous question is ``\emph{yes}'', what is the value of $\alpha$?
\end{itemize}
\end{prob}

The question of whether the upper bound obtained in  ~\cite{bartholdi_Ers:growth10,muchnik_p:growth}
is optimal (i.e.\ it coincides with the growth rate of the group $\mathcal{G}$) is very intriguing
and we formulate it as the
next problem. The point is that all known groups whose growth is explicitly  computed up to equivalence by
$\sim$ (i.e.\ groups considered in~\cite{bartholdi_Ers:growth10,bartholdi_ersch:given(2)11}) have growth  not
smaller than $e^{n^{\alpha_0}}$. And it looks that if the growth of $\mathcal{G}$ is less than
$e^{n^{\alpha_1}}$ with $\alpha_1 < \alpha_0$ then, using the results of the cited papers and of the paper by
Kassabov  and Pak~\cite{kassabov_pak:11}, one can extend the range of possible growth rates from the
``interval'' $[e^{n^{\alpha_0}},e^n]$ to the ``interval'' $[e^{n^{\alpha_1}},e^n]$. But so far
$e^{n^{\alpha_0}}$ is a kind of a ``mountain'' which ``closes the sky'' for people working in the area of
group growth.

In view of recent results from \cite{bartholdi_Ers:growth10,brieussel:growth11}, even obtaining a weaker
result that would answer  the next question seems to be interesting in its own right.

\begin{prob}\label{pak}\quad\\[-\baselineskip]
\begin{itemize}
\item[(1)]  Is the upper bound (\ref{barth})   the best possible for $\mathcal{G}$?

\item[(2)] Does there exist a group of intermediate growth whose growth is less than $e^{n^{\alpha_1}}$, where
$\alpha_1 < \alpha_0$ and  $\alpha_0=\log 2/\log (2/\rho)$ is the constant defined above?
\end{itemize}
\end{prob}

In section~\ref{inter2}, we will give a brief account of methods that can be used to obtain upper and lower
bounds for intermediate growth.

\medskip
The groups $\mathcal{G}$ and  $\mathcal{E}$ belong to the class of self-similar groups, that is, groups
generated by automata of Mealy type, which are discussed a bit in the next section. An important quantative
characteristic of such groups is the pair $(m,n)$ which is a rough indication of the \emph{complexity} of the
group, where $m$ is the cardinality of the alphabet and $n$ is the cardinality of the set of the states of
automaton. From this point of view, $\mathcal{G}_{\xi}$ and $\mathcal{G}_{\eta}$  are $(2,5)$-groups (the
Moore diagrams of corresponding automata are presented in the next section).   And there is even a group of
intermediate growth of complexity $(2,4)$  (namely the iterated monodromy group $IMG(z^2+i)$ in the sence of
Nekrashevych~\cite{nekrash:self-similar} of quadratic polynomial $z^2+i$) as was showed by K.~Bux and
R.~Peres \cite{bux_p:iter_monodromy}.

\begin{prob}\quad\\[-\baselineskip]
\begin{itemize}
\item[(i)] Find the growth degree of each of the groups $\mathcal{G}_{\xi}$, $\mathcal{G}_{\eta}$, and
  $IMG(z^2+i)$.
 \item[(ii)]  Are there groups of intermediate growth of complexity $(2,3)$?

 \item[(iii)]  Determine all automata of complexity  $(2,3)$, $(2,4)$ and $(2.5)$ which generate groups of
     intermediate growth.

 \item[(iv)] Determine all possible types of growth of self-similar groups generated by finite automata.
\end{itemize}
\end{prob}


It is a kind of  a miracle that  an automaton with a small number of states can generate a group with very
complicated algebraic structure and asymptotic behavior. Therefore, it is not surprising that some of the
automata groups studied prior to 1983 (when the the first examples of groups of intermediate growth were
found) are also of intermediate growth. Specifically, the 2-group of Aleshin \cite{aleshin:burnside}
(generated by two automata with 3 and 7 states) and the $p$-groups of V.~Sushchanskii (generated by automata
with a number of states growing as a quadratic function of $p$ \cite{sushch:burnside}) were shown to have
intermediate growth by the author in \cite[page 280]{grigorch:degrees} and \cite[page
197]{grigorch:degrees85}, respectively  (see also \cite{bondar_sav:susch07} where the Sushchanskii group is
treated in detail). It is worth mentioning that the papers of Aleshin and Sushchanskii deal exclusively with
the question of construction of  finitely generated infinite torsion groups (contributing to the
\emph{General Burnside Problem}) and Milnor's problem is not considered at these articles at all.

\medskip

As all the groups of intermediate growth from theorem~\ref{growth2} have only finite quotients (and consequently,
at most a countable set of quotients), in 1983 it was reasonable to ask if there are groups of intermediate
growth with
uncountably many homomorphic images, one of the properties that a finitely generated virtually nilpotent
group does not have.
 This was affirmatively answered in \cite{grigorch:continuum84}.
The next theorem gives a hint to the main result of \cite{grigorch:continuum84} and its proof.
\smallskip

 Let $\Lambda \subset \Omega$ be the subset
 consisting of sequences that are products of blocks
$012, 120, 201$,
and let
 $\mathcal{G}_{\lambda}$ be presented as the quotient  $F_4/\mathcal{N}_{\lambda}$ of the free group $F_4$ of rank 4 by a normal subgroup $\mathcal{N}_{\lambda}$. Call the group
 $\mathcal{U}_{\Lambda}=F_4/\bigcap_{\lambda \in \Lambda}\mathcal{N}_{\lambda} $
 $\Lambda$-universal.

 \begin{thm} The  $\Lambda$-universal group $\mathcal{U}_{\Lambda}$ has
 intermediate growth and
 has uncountably many quotients which are pairwise non-isomorphic.
 \end{thm}

The first part of the theorem follows from \cite[theorem 3.2]{grigorch:degrees}.

 As each of the groups $\mathcal{G}_{\lambda}$ is a homomorphic image of $\mathcal{U}_{\Lambda}$, the second
 part of the theorem is obvious modulo the fact that the classes of isomorphisms of groups from the family
 $\mathcal{G}_{\omega}$, $\omega \in \Omega$ are at most countable. It was shown by Nekrashevych in \cite{nekrash:minimal07}  that
 $\mathcal{G}_{\omega}  \simeq \mathcal{G}_{\zeta}$ if and only if the sequence $\zeta$ can be obtained from
 $\omega$ by the diagonal action at all coordinates of  an element from  the symmetric group $Sym(3)$ acting on the  set
 $\{0,1,2\}$. (In fact, it is enough to quote  theorem~5.1 from \cite{grigorch:degrees}, which states that for each $\omega$ there are at most countably many groups $\mathcal G _{\eta}$ isomorphic to
$\mathcal G _{\omega}$).

{\StretchyLines
 In contrast with  the $\Lambda$-universal group
$\mathcal{U}_{\Lambda}$, the growth of
the $\Omega$-universal group $\mathcal{U}_{\Omega}=F_4/\bigcap_{\omega \in
\Omega}\mathcal{N}_{\omega}$ is exponential \cite{muchnik:amenab_of_universal}. It is known that
$\mathcal{U}_{\Omega}$ does not contain a free subgroup on two generators, is self-similar }
 (of complexity $(6,5)$ \cite{grigorch:solved}), weakly branch and contracting.  However, the following question is still open
 (unfortunately the article \cite{muchnik:amenab_of_universal}, which contains the claim about
 amenability of $\mathcal{U}_{\Omega}$, has  a mistake).

\begin{prob}  Is $\mathcal{U}_{\Omega}$ amenable or not?
\end{prob}

\medskip
In \cite{grigorch:degrees} the author proved that for any monotone function $\rho(n)$ growing slower than
exponential functions, there is a group with  growth not slower than $\rho(n)$ (so either $\rho(n)\preceq
\gamma_G(n)$ or  $\rho(n)$ and $\gamma_G(n)$ are incomparable with respect to the preoder $\preceq$). This result
was improved by Erschler .

\begin{thm} \label{approachexp}(Erschler \cite{erschler:degrees05}) For any increasing function $\rho(n)$ growing
slower than an exponential function, there is a finitely generated
group $G$  of intermediate growth with $\rho(n)\preceq \gamma_G(n).$
\end{thm}

The last result shows that there is no upper bound for intermediate growth, in contrast with the lower bound
given by the Shalom-Tao function $ n^{c(\log \log n)^c}$ ($c$ a constant), as discussed in section~\ref{polynom0}.

\section{The Gap Conjecture}\label{gapconj}

The history of the Gap Conjecture is as follows. 
While reading Gromov's paper on polynomial growth in 1982 (soon after its publication),
the author realized that the effective version of Gromov's polynomial growth
theorem (Theorem~\ref{gromov})  implies the existence of a function $\upsilon(n)$ growing faster than any
polynomial such that, if $\gamma_G(n)\prec \upsilon(n)$, then the growth of  $G$ is polynomial. Indeed, taking a
sequence $\{k_i,d_i\}_{i=1}^{\infty}$ with $k_i \to \infty $ and $d_i \to \infty$ as $i \to \infty$ along with the
corresponding sequence $\{R_i\}_{i=1}^{\infty}$ (whose existence follows from Theorem~\ref{gromoveffect}), one
can build a function $\upsilon(n)$ which coincides with the polynomial $k_in^{d_i}$ on the interval
$[R_{i-1}+1,R_i]$.  The constructed function $\upsilon(n)$ grows faster than any polynomial and separates
polynomial growth from intermediate growth. In fact, as was already mentioned in section~\ref{polynom0} and at the
end of previous section,  the function $ n^{c(\log \log n)^c}$ ($c$ a constant) separates polynomial growth
from exponential~\cite{shalom_tao:polynom10}.  As, according to Theorem~\ref{growthmain}, the set of growth
degrees is not linearly ordered, it may happen that there is more than one ``gap'' between polynomial growth
and intermediate growth. But in any case, it would be nice to obtain the best possible estimate of the
asymptotics of a function which ``uniformly'' separates the polynomial and intermediate growth.

Approximately at the same time, while reading Gromov's paper (thus around 1982), the author was establishing
his results on groups of intermediate growth discussed in the previous section and in 1983-1985 published
\cite{grigorch:milnor,grigorch:degrees,grigorch:degrees85}. The lower bound of the type $e^{\sqrt n}$  for
all groups $\mathcal G_{\omega}$ of intermediate growth established in those papers and in his Habilitation
\cite{grigorch:habil} allowed to author to guess that the equivalence class of the function $e^{\sqrt n}$ could be
a good candidate for a ``border'' between polynomial and exponential growth.  This guess became stronger in
1988 when the author obtained the results published in~\cite{grigorch:hilbert} (see Theorem~\ref{gap1}).
In the ICM Kyoto paper~\cite{grigorch:ICM90}, the author raised a question of whether the function  $e^{\sqrt n}$
gives a universal lower bound for all groups of intermediate growth. Moreover, at approximately the same time,
he conjectured that indeed this is the case, and stated this later at numerous talks.

\begin{conj}  \label{c:gap1}(Gap Conjecture) If the growth function $\gamma_G(n)$ of a finitely
generated group $G$ is strictly bounded from
above by   $e^{\sqrt{n}}$ (i.e.\ if $\gamma_G(n) \prec e^{\sqrt n}$), then the growth of $G$ is polynomial.
\end{conj}

We are also interested to know whether there is a group, or more generally a cancellative semigroup, with growth
equivalent to $e^{\sqrt{n}}$.

\smallskip
In this section we formulate several results in the direction of confirmation of the Gap Conjecture, and
suggest slightly different versions of it. Later in section~\ref{probabil} we will formulate  analogous
conjectures about some other asymptotic characteristics of groups.

The next few results, together with the results about lower bounds on growth discussed in the next section,
are the main  source of  support of the Gap Conjecture. Recall that a group $G$ is said to be a residually
finite-$p$ group (sometimes also called residually finite $p$-group) if it is approximated by finite
$p$-groups, i.e.\ for any $g\in G$ there is a finite $p$-group $H$ and a homomorphism $\phi\:G\rightarrow H$
with $\phi(g)\neq 1$. This class is, of course, smaller than the class of residually finite groups, but it is
pretty large.  For instance, it contains Golod-Shafarevich groups, the $p$-groups $\mathcal G_{\omega}$ from
\cite{grigorch:degrees,grigorch:degrees85}, and many other groups.

\begin{thm} \label{gap1} (\cite{grigorch:hilbert})  Let $G$ be a finitely generated residually finite-$p$  group.
If $\gamma_G(n)\prec e^{\sqrt{n}}$ then $G$ has polynomial growth.
\end{thm}

As was established by the author in a discussion with A.~Lubotzky and A. Mann during the conference on
profinite groups in Oberwolfach in 1990,
  the same arguments as given in \cite{grigorch:degrees} combined with
Lemma~1.7 from \cite{lubotzky_mann:polyn91}  allow one to prove a stronger version of the above theorem (see
also the remark after Theorem~1.8 in \cite{lubotzky_mann:polyn91}, but be aware that capital~$O$ has to
be replaced by small~$o$).

\begin{thm}\label{lubmann}
Let $G$ be a residually nilpotent finitely generated group. If $\gamma_G(n)\prec e^{\sqrt{n}}$ then $G$ has
polynomial growth.
\end{thm}

The main goal of the  paper~\cite{lubotzky_mann:polyn91}, followed by the article~\cite{lubotzky_man_seg:93},
was to give a complete description of finitely generated groups with polynomial subgroup growth (the growth
of the function which counts the number of subgroups of given finite index). The remarkable result achieved
in~\cite{lubotzky_man_seg:93} shows that such groups are precisely the \emph{solvable groups of finite rank}.

As was already mentioned in the introduction,  the original proof of Gromov's polynomial growth theorem is
based on the use of the solution of Hilbert's 5th problem by Montgomery and Zippin, concerning the
isometric actions of locally compact groups and their relation to Lie groups. Surprisingly, in the proofs of
the results stated in Theorems~\ref{gap1} and~\ref{lubmann}, as well as in the results from
\cite{lubotzky_mann:polyn91,lubotzky_man_seg:93} about polynomial subgroup growth, M. Lazard's solution of the $p$-adic analog of Hilbert's 5th problem~\cite{lazard:analytic65} plays an important role.
The result of Lazard gives a characterization of analytic pro-$p$-groups. After a long period of search, a proof
of Gromov's Theorem which avoids the use of the 5th Hilbert problem was found by
B.~Kleiner~\cite{kleiner:plynom10}.  Now we will formulate a theorem (Theorem~\ref{polyciclic}), which
generalizes theorems~\ref{gap1} and~\ref{lubmann}, and whose proof  is based on the techniques of J.S.~Wilson
from~\cite{wilson:growthsolv05,wilson:gap11} and some other results.  Wilson's arguments how to handle with
growth of resudually finite groups are quite original (the techniques of ultraproducts  is used at some
point), but eventully they reduce the arguments  to the case of residually nilpotent groups (i.e.\ to the
previous theorem). It would be interesting to find a proof of the theorem~\ref{gap1} which avoids the use
of the $p$-adic version of Hilbert's 5th problem.

Recall that a group is called supesolvable if it has  a finite  normal descending  chain of subgroups with
cyclic quotients. Every finitely generated nilpotent group is supersolvable~\cite{robinson:book96}, therefore
the next theorem improves Theorem~\ref{lubmann}.

\begin{thm} \label{polyciclic}
The Gap Conjecture holds for residually supersolvable  groups.
\end{thm}

The proof of this theorem  is based on the techniques of J.S.~Wilson developed
in~\cite{wilson:growthsolv05,wilson:gap11} for studies around Gap Conjecture. Currently, it is not known if
the main Conjecture holds for residually polycyclic and, more generally, residually solvable groups.  However, it
is quite plausible that it does.

\begin{conj}  \label{c:gap2}(Gap Conjecture with parameter $\beta$) There exists $\beta$, $0<\beta<1$,
such that if the growth function
$\gamma_G(n)$ of a finitely generated group $G$ is strictly bounded from above by $e^{n^\beta}$
(i.e.\ if $\gamma(n) \prec e^{n^\beta}$), then the growth of $G$ is polynomial.
\end{conj}

Thus the Gap Conjecture with parameter $1/2$ is just the Gap Conjecture~(Conjecture~\ref{c:gap1}).
If  $\beta<1/2$ then
the Gap Conjecture with parameter $\beta$ is weaker than the Gap Conjecture,  and if  $\beta>1/2$ then it is
stronger than the  Gap Conjecture.

\begin{conj} \label{c:gap3} (Weak Gap Conjecture)  There is a $\beta, \beta <1$ such that if  $\gamma_G(n) \prec
e^{n^\beta}$ then the Gap Conjecture with parameter $\beta$ holds.
\end{conj}

As was already mentioned, there are some results of J.S.~Wilson  in the direction of confirmation of the
Gap Conjecture. He showed that, if $G$ is a residually solvable group whose growth is strictly less than
$e^{n^{1/6}}$, then it has polynomial growth~\cite{wilson:growthsolv05,wilson:gap11}. Therefore the Gap
Conjecture with parameter $1/6$ holds for residually solvable groups. The proof of Wilson's result is based
on the estimate of the rank of chief factors of finite solvable quotients of $G$. The methods
of~\cite{wilson:growthsolv05,wilson:gap11} combined with the theorem of
Morris~\cite{morris:amenable_gps_acting_on_line} and theorem of Rosset~\cite{rosset:76} can be used to prove
the following statement.

\begin{thm}\label{orderable}\quad\\[-\baselineskip]
\begin{itemize}
\item[(i)] The Gap Conjecture  with parameter $1/6$ holds for left orderable groups.
\item[(ii)] The Gap Conjecture holds for left orderable groups if it holds for residually polycyclic groups.
\end{itemize}
\end{thm}

Let us also formulate an open problem which is related to the  above discussion.

\begin{prob} \label{gap}\quad\\[-\baselineskip]
\begin{itemize}
\item[(i)]  Does there exist  $\alpha$, $0 < \alpha <1$ such that if the growth of a finitely generated group is
strictly less than the growth of $e^{n^{\alpha}}$ then it is polynomial?

\item[(ii)]  If such $\alpha$ exists, what is its maximal value?  Is it $< 1/2, = 1/2$ or $> 1/2$?

\item[(iii)] In the case $\alpha$ exists (and is chosen to be maximal), is there a group (or a cancellative
semigroup) with growth equivalent to $e^{n^{\alpha}}$?

\item[(iv)]  Is there a finitely generated group approximated by nilpotent groups with  growth equivalent to
$e^{\sqrt{n}}$?  Is there a residually  finite-$p$ group with growth  equivalent to $e^{\sqrt{n}}$? ($p$ -
fixed prime).
\end{itemize}
\end{prob}

There is some evidence  based on considerations presented in  the last section of this article and some
additional arguments that the
above  conjectures and  problem (parts (i), (ii), (iii)) can be reduced to consideration of  three classes of
groups: \emph{simple} groups, \emph{branch} groups and  \emph{hereditary just-infinite}  groups. These three
types of groups appear in a natural division of the class of just-infinite groups
into three subclasses described in Theorem~\ref{just-inf}.
Branch groups are defined in section~\ref{inter2} and a hereditary just-infinite group is a residually
finite group with  the property that every proper quotient
    of every subgroup of finite index (including the group itself) is just-infinite.
In any case, the following theorem holds (observe that branch groups and hereditary just-infinite groups are
residually finite groups).

\begin{thm}\label{resid}
If the Gap Conjecture holds for the classes of  residually finite groups and simple groups, then it holds for
the class of all groups.
\end{thm}

This theorem is a corollary of the main result of \cite{bajorska_maked:note07}.  A different  proof is suggested
in \cite{grigorch:gapconj12}. It is adapted to the needs of the proof of Theorem \ref{wilson},
which is discussed at the end of the article.

As it is quite plausible that the Gap Conjecture could be proved for residually finite groups (the
classification of finite simple groups may help), we suspect that the validity of the Gap Conjecture depends on
its validity for the class of simple groups. We will return to just-infinite groups at the end of the article
and state one more reduction of the Gap Conjecture.

It is unknown if there are simple groups of intermediate growth (problem~\ref{medyn}), but a recent article of
K.~Medynets and the author~\cite{grigorch_medyn:simple11} shows that there are infinite finitely generated
 simple groups  which belongs to the class $LEF$ (locally embeddeble into
 finite groups); the authors
 conjectured that these groups are amenable.
This was recently confirmed by K. Juschenko and N. Monod \cite{juschenko_monod:12}.

 This gives some hope that   groups of intermediate growth may exist within the
subgroups of the groups considered in~\cite{grigorch_medyn:simple11} (namely among subgroups of \emph{full
topological groups} $[[T]]$ associated with   minimal homeomorphisms $T$ of a Cantor set). Algebraic
properties of $[[T]]$ were studied by H.~Matui~\cite{matui:simple06}, who showed that their commutator
subgroup $[[T]]'$ is simple, and is finitely generated in the case when the homeomorphism $T$ is a minimal
subshift over a finite alphabet. Observe that $[[T]]'$ always has exponential growth, as  was recently
shown in \cite{matui:growth11}. Therefore the only hope is that groups of intermediate growth may exist
among subgroups of $[[T]]$.

There are several other Gap Conjectures related to various asymptotic characteristics of groups.  We list
some of them in section~\ref{probabil} and discuss briefly their relation to the problems and conjectures
considered in this section.

\section{Intermediate growth: the upper and lower bounds} \label{inter2}

In this section we give an overview of the main methods of getting  upper and  lower bounds of  growth in the
intermediate growth  case.  We begin with upper bounds. For establishing if a group is of intermediate growth,
it is more important to have tools to obtain upper bounds because as soon as it is known that a group is not
virtually nilpotent (and usually this is not difficult to check),  its growth is known to be superpolynomial
(by Gromov's theorem).
Therefore, if a finitely generated infinite group possesses any property such as being  simple, torsion, not
residually finite, nonhopfian etc.,  one immediately knows that the growth is superpolynomial.

Finding a lower bound for the growth is of interest not only because of its connection with the
Gap Problem~\ref{gap} discussed above,
 but also because of the  connection with various topics in the theory of random walks on groups and spectral
 theory of the discrete Laplace operator.

\medskip
The   method  for getting upper bounds for  the growth in the intermediate case that we are going to
describe was used in \cite{grigorch:degrees,grigorch:degrees85}. Roughly, the idea consists  in encoding of
each element $g$  of a group $G$ by a set of $d$ elements $g_1,\dots, g_d$, $g_i \in G$ ( $d\geq 2$ and fixed)
in such  a  way that for some fixed constants $C$ and $\lambda$, $0< \lambda <1$ (independent of $g$) the
inequality
\begin{equation}\label{growthst1}
 \sum_{i=1}^d|g_i| \leq \lambda |g| + C
\end{equation}
holds.  The meaning of this inequality is that an element of  (large) length $n$ is coded by a set of $d$
elements of the total length  strictly less than $n$ (with coefficient of the reduction $\lambda$).  In
dynamics, this situation corresponds to the case when the entropy of the system is  zero, while in the context
of growth it corresponds to the case when the constant $\kappa$ given by (\ref{growthexp})
 is equal to 1 (i.e.\ when the growth is subexponential).

\smallskip

  There are some variations of condition~(\ref{growthst1}). For instance,  in some cases  the coefficient
$\lambda $ can be taken to be equal to $1$, but some  additional conditions on
 the group have to be
 satisfied in order to claim that the growth is subexponential. For instance this happens in the case of the group $\mathcal{E}$,
 and more generally in the case of the
 groups $\mathcal{G}_{\omega}$, $\omega \in \Omega_1$.

At the moment, the above  idea  is realized only in the case of certain groups acting on \emph{regular rooted
trees} that are \emph{self-similar} or have  certain self-similarity features. The simplest property that
leads to intermediate growth is the  \emph{strong contracting property} that we are going to describe.

The set of vertices of a $d$-regular rooted tree $T=T_d$
 is in a natural bijection with the set of finite words  over an alphabet $X$ of cardinality $d$
  (usually one of the alphabets $\{0,1,\dots,d-1\}$ or
$\{1,2,\dots,d\}$ is used).   The set of vertices of the tree is  graded by levels $n=0,1,2,\dots$, and
vertices of the $n$-th level are in natural bijection with the words of length $n$ (i.e.\ with the
elements of the  set $X^n$) listed in lexicographical order. Let $V=X^{\ast}=\bigcup_{n=1}^{\infty}X^n$
be the set of all vertices.

For the full group $Aut (T)$ of automorphisms of a $d$-regular rooted tree $T$, a  natural decomposition into a
semidirect product
\begin{equation}\label{semidir1}
 (Aut (T) \times  \dots \times Aut (T))\rtimes Sym(d)
\end{equation}
($d$ factors) holds, as well as  a corresponding  decomposition of any element $g\in Aut (T)$
\begin{equation}\label{semidir2}
g=(g_1,g_2,\dots,g_d)\sigma.
\end{equation}
  Here the element $\sigma$ of  the symmetric group $Sym(d)$  shows how $g$ acts on the vertices of the first
level, while the \emph{projections} $g_i$ show how $g$ acts on the  subtree $T_i$ which has its root at vertex $i$
of the first level (we  identify the subtree $T_i$ with $T$ using the canonical  self-similarity of the
regular rooted tree). In a similar  way, the projection  $g_v$ of an element $g \in Aut (T)$ can be defined
for an arbitrary vertex $v\in V$.

 \begin{defn} The action of a  group $G$  on a  regular rooted tree $T$ is called \emph{self-similar} if for any
 vertex $v$, the projection $g_v$ is again an element of $G$ modulo the canonical identification of the
 subtree $T_v$
  with the original tree $T$. A group  is called self-similar if it has a faithful self-similar action.
\end{defn}

For instance, the groups $\mathcal{G}$ and $\mathcal{E}$ are self-similar. For  generators $a,b,c,d$ of
$\mathcal{G}$ the following relations of the type
  (\ref{semidir2})
  hold
  \begin{equation}\label{selfsim1} a=(1,1)\sigma, b=(a,c)e, c=(a,d)e, d=(1,b)e,
  \end{equation}
  where  $e$ is the identity element  and  $\sigma$ is a permutation (both are elements of $Sym(2)$ acting on the
  alphabet $\{0,1\}$). Observe that  usually the above equalities are used in the simplified form
  \begin{equation}\label{selfsim2}
  a=\sigma, b=(a,c), c=(a,d), d=(1,b).
  \end{equation}

An equivalent definition of a self-similar group is that it is a group generated by states of non-initial
Mealy type automaton with the operation of composition of automata. We are not going to explain here what
the Mealy  automaton is, nor the group defined by it. For this approach to self-similarity, see
\cite{gns00:automata,nekrash:self-similar}). For instance, the groups $\mathcal{G}$  and $\mathcal{E}$ are
groups generated by  the following automata




\begin{center}\label{pict3}
\includegraphics[width=\hsize]{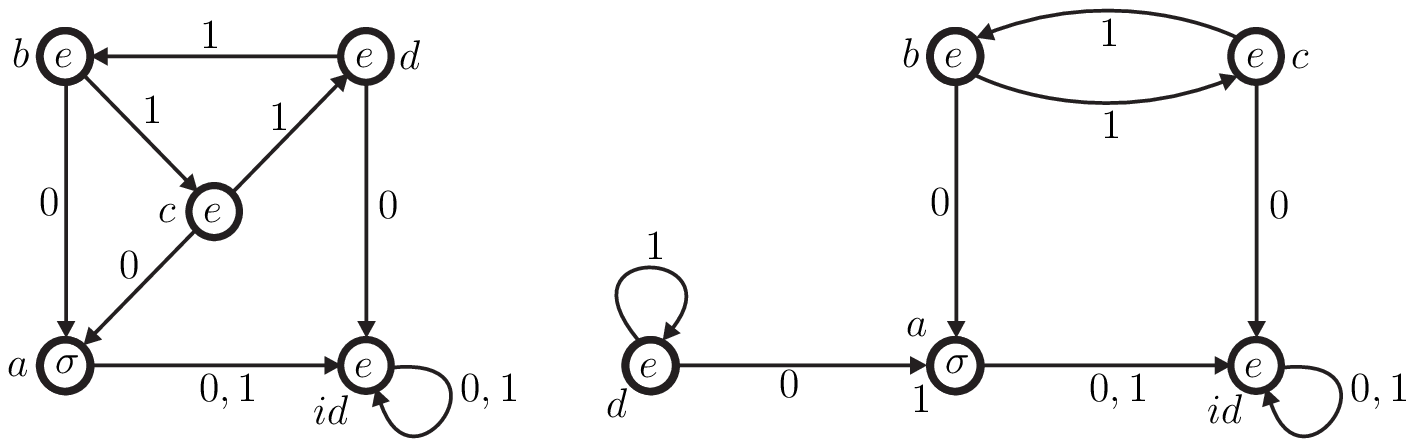}

Automata generating $\mathcal G$  and $\mathcal E$
\end{center}

A challenging problem is to understandi what the class self-similar groups is, and especially what constitutes
the subclass consisting of  groups that are generated by finite automata.  There is a description of groups
generated by 2-state automata over the alphabet on two letters  (there are 6 such groups), see
\cite{gns00:automata}.  There are not more that 115 groups generated by 3-state automata over an alphabet on
two letters (see \cite{bondarenko_gkmnss:full_clas32} where an \emph{Atlas} of self-similar groups is
started).

\medskip

 Another important notion useful when studying growth  is the notion of a \emph{branch} group, introduced in
 \cite{grigorch:branch,grigorch:jibranch}.
 With each sequence $\bar{m}=\{m_n\}_{n=1}^{\infty}$ consisting of integers $m_n \geq 2$, one can associate a
 spherically  homogeneous rooted tree $T_{\bar{m}}$  as done in
 \cite{bass_or:cyclic,grigorch:jibranch,bar_gs:branch}. Given a group $G$ acting faithfully by automorphisms
 on such a tree $T=T_{\bar{m}}$ and a vertex
 $v\in V(T)$, the rigid stabilizer $rist_G(v)$ is defined as the subgroup in $G$ consisting of elements which act
 trivially outside the subtree $T_v$ with its root at $v$.  The rigid stabilizer $rist_G(n)$ of level $n$
 is the subgroup generated by the rigid stabilizers of vertices of $n$th level (it is their internal direct
 product).

 \begin{defn}  \label{branch} A group $G$ is called  branch  if it has a faithful, level-transitive action on some
 spherically homogeneous rooted tree $T_{\bar{m}}$  with the property that
 \[ [G:rist_G(n)] < \infty\]
 for each $n=1,2,\dots$.

 \end{defn}

 This is a geometric definition of branch groups.  There is also an algebraic definition discussed in
 \cite{grigorch:jibranch,bar_gs:branch} which gives a slightly larger class of groups.

 As was already mentioned in the introduction,   branch groups of self-similar type are basically the only
 source for the constructions of groups of intermediate growth and, as is stated in the theorem
 \ref{growth2}, the
 groups $\mathcal{G}_{\omega}$, $\omega \in \Omega_1$, are branch. Branch groups constitute one of three classes
 in which the class of just-infinite groups naturally splits. The role of this class in the theory of growth of groups will
 be emphasized  once more in the last section of the paper.

\begin{defn} A finitely generated self-similar group is called \emph{contracting}  with parameters $\lambda <1 $  and $C$ if the inequality
\begin{equation}\label{project}
|g_i| < \lambda|g| +C
\end{equation}
holds for every element $g\in G$ and any of its section $g_i$  on the first level.
\end{defn}

It follows from (\ref{project}) that
\[|g_i| < |g| \]
if $|g|>C/(1-\lambda)$.  Therefore for elements of sufficiently large length the projections on the vertices
of the first level are shorter.

\medskip
The contracting property and self-similarity allow one to study algebraical properties of the group.  This was first
used in \cite{grigorch:burnside}  to show that the groups $\mathcal{G}$  and $\mathcal{B}$ are torsion, and
later many other interesting and sometimes unusual properties of the group $\mathcal{G}$ and similar type
groups (for instance Gupta-Sidki $p$-groups \cite{gupta_s:burnside}) were established.

To be more precise, the groups $\mathcal{G}_{\omega}$ are not self-similar for a typical $\omega$. But in
fact, the whole family $\{\mathcal{G}_{\omega}\}_{\omega \in \Omega}$ of groups is self-similar, because the
projections of any element $g \in \mathcal{G}_{\omega}$ on vertices of arbitrary level $n$ belong to the
group $\mathcal{G}_{\tau^n(\omega)}$ where $\tau$ is the shift in the space of sequences. Moreover the
contracting property holds for this family with respect to the canonical generating sets
$\{a,b_{\omega},c_{\omega},d_{\omega}\}$. This allows one to get various results (including the estimates of the
growth) using a simultaneous induction on the length of elements for infinite shift invariant families of
groups from the class  $\{\mathcal{G}_{\omega}\}_{\omega \in \Omega}$.

\begin{defn} A  finitely generated self-similar group acting on a $d$-regular rooted tree $T$ is called
\emph{strictly contracting}  (or \emph{sum contracting}) with parameters $k\geq 1$, $\lambda <1 $ and $C$ if
there is a level $k$ of the tree such that the inequality
\begin{equation}\label{growth9}
 \sum_{v \in V_k}|g_v| \leq \lambda|g| +C
\end{equation}
holds for an arbitrary element $g \in G$ ($V_k$ denotes the set of vertices of level $k$).
\end{defn}

The idea of encoding of the elements of length $n$ by tuples of elements of total length less than $n$
discussed above  has a clear implementation for strictly contracting groups.  Namely, every element $g\in G$
can be encoded by the set $\{\bar{g},g_1,g_2,\dots, g_{d^k}\},$ where $\bar{g}$ is a shortest representative
of the coset $gst_G(k)$ of the stabilizer $st_G(k)$ of $k$th level, and $g_1,g_2,\dots,g_{d^k}$ are
projections on the vertices of $k$-th level
 (observe that the length of $\bar{g}$ is uniformly bounded by some constant independent of the choice of the element).
\begin{prop} (\cite{grigorch:degrees})  Let $G$ be a strictly contracting group with parameters $k, \lambda$ and $C$.  Then there is a constant
$\beta =\beta(\lambda,k) <1$ such that
\[\gamma_G(n) \preceq  e^{n^{\beta}}.\]

\end{prop}

The strict contracting property  holds, for instance, for the groups $\mathcal{G}_{\omega}, \omega \in \Theta$. For
the groups $\mathcal{G}_{\omega}$, $\omega \in \Omega_1 \setminus \Omega_0$, condition~(\ref{growth9}) does not
hold in general. Regardless, a  modification of the arguments allows one show that the growth is intermediate in this
case as well, and even to get an upper bound of the type shown in the third part of the theorem~\ref{GRI83}.

\medskip

Now let us turn  to lower bounds. There is a much wider variety of methods used to obtain lower bounds than
for upper bounds. We are going to mention the \emph{branch commensurability} property, and give a short
account (partly here, partly in the next section) of other methods: the \emph{anti-contracting property}
 (which is a kind of opposite to the contracting property), the \emph{Lie method}, the \emph{boundary of random
walks} method \cite{erschler:boundary04} and the \emph{Munchhausen trick} method
\cite{bartholdi_v:amenab,kaiman:munchhausen,grigorchuk-n:schur}.

 The first observation is that the group $\mathcal{G}$ and many
other self-similar groups of branch type, for instance the so called \emph{regularly branch} groups
\cite{bar_gs:branch}, have the property of being \emph{abstractly commensurable} to some power $d\geq 2$  of
itself (i.e.\ to a direct product of $d$ copies of the group, in the case of $\mathcal{G}$ the power   $d=2$).
Recall that two groups are said to be abstractly commensurable if they have isomorphic subgroups of finite
index. Observe that there are even finitely generated groups $G$ which are isomorphic to their proper powers
$G^d$, $d \geq2$ \cite{jones:hopf74,baumslag_miller:odd88}. Such groups are non-hopfian (and hence not
residually finite) in contrast to regularly branch groups, but among them, groups of intermediate growth have
not yet been detected.

\begin{prob}  Does there exist a finitely generated group $G$ of intermediate growth isomorphic to some power
$G^d$, $d\geq 2$?
\end{prob}

\smallskip

For a group $G$ commensurable with $G^d$  there is a lower bound on  growth given by the inequality
\begin{equation}\label{growth3}
e^{n^{\alpha}}\preceq \gamma_G(n)
\end{equation}
for some positive $\alpha$.
 Verifying this  is an easy exercise if one keeps in mind that the growth of a group coincides with that of a
subgroup of finite index, and that
\[ \gamma_{G^d}(n) \sim \gamma^d_G(n).\]
Observe that in general we do not have  control over $\alpha$ in this type of argument.


Another condition  that allows us to get a lower bound of type (\ref{growth9}) is the following. Let $G$ be a
self-similar group acting on a $d$-regular tree.

\begin{defn} A group $G$ satisfies the \emph{anti-contracting} property with parameters $k$, $\mu$ and $C$   if   for arbitrary $g\in G$ the inequality
\begin{equation}\label{growth6}
|g|\leq \mu \sum_{v \in V_k} |g_v| +C
\end{equation}
holds, where $V_k$ is the set of vertices of level $k$ and $g_v$ is the projection of $g$ on $v$.
\end{defn}

Inequality (\ref{growth6}) is a kind of opposite to inequality (\ref{growth9}).  Using the methods
of~\cite{grigorch:degrees} one can prove the following fact.

\begin{thm}  \label{growth4} Let $G$ be a group satisfying  the anti-contracting property with parameters $\mu$ and $C$. Then there is $\alpha=\alpha(k,\mu)>0$ such that inequality (\ref{growth3}) holds.
\end{thm}

The lower bound for $\alpha$ in terms of $k$ and $\mu$ ($C$ is not important) can be explicitly written. In
the case of the group $\mathcal{G}$, the \emph{anti-contracting} property holds  with parameters $k=1,\mu=2$,
$C=1$ and $\alpha=1/2$. In fact, a lower bound  of type  $e^{\sqrt{n}}$ holds for all groups from
the family $\mathcal{G}_{\omega}, \omega \in \Omega_1$.

Another approach to obtain lower bounds is via an idea of  W.~Magnus: using Lie algebras and associative algebras
associated with  a group.  Given a  descending central series $\{G_n\}_{n=0}^{\infty}$ of a group $G$, one can
construct the graded Lie ring $\mathcal{L}=\oplus_{n=0}^{\infty}G_n/G_{n+1}$ and the graded associative
algebra $\mathcal{A}=\oplus_{n=0}^{\infty}A_n$, where $A_n=\Delta^n/\Delta^{n+1}$, and $\Delta$ is the
fundamental ideal of the group algebra $\mathbb{F}[G]$   ($\mathbb{F}$ a  field).  The Lie operation in
$\mathcal{L}$ is induced by the operation of taking commutators of pairs of elements in the group (defined
first on the abelian quotients $G_n/G_{n+1}$ and then extended by linearity to the whole ring).  The
important cases are given by the lower central series and the Jennings-Lazard-Zassenhaus
 lower $p$-central series, which in the case of a simple field of characteristic $p$ can be defined
as $G_n=\{g\in G : 1-g \in \Delta^n\}$.

By Quillen's theorem, the algebra $\mathcal{A}$ is the universal enveloping algebra of $\mathcal{L}$
(or $p$-universal in the case $char~{\mathbb{F}}=p$).
There is a close relationship between the growth of the algebras
$\mathcal{A}$ and $\mathcal{L}$, i.e., the growth of the dimensions of the homogeneous components of these
algebras. Namely, $\mathcal{L}$ has exponential growth if and only if $\mathcal{A}$ has exponential growth,
and if $\mathcal{L}$ has polynomial growth of degree $d$ then $\mathcal{A}$ has intermediate growth of type
$e^{n^{\alpha}}$ with $\alpha=\frac{d+1}{d+2}$ (for details see e.g.,
\cite{grigorch:hilbert,bartholdi-g.Lie00}). More information about growth of algebras can be found in
\cite{krauser-lenagan:growth00}.  Observe that finitely generated Lie algebras may be of fractional (and even
irrational) power  growth as is shown in~\cite{petrograd_shestzel:10}.

The following fact shows that the growth of $\mathcal{A}$ gives a universal lower bound for the growth of a
group independently of the system of generators.

\begin{prop} (~\cite{grigorch:hilbert})  \label{growth1}
Let $G$ be a finitely generated group with a finite system $S$ of semigroup generators (i.e.\ each element of
$G$ can be expressed as a product of elements from $S$).  Let $\gamma_G^S(n) $ be the growth function of $G$
with respect to $S$ and $a_n=dim_{\mathbb{F}}\mathcal{A}_n$.  Then, for any $n \in \mathbb{N}$,
\begin{equation}\label{gap4}
\gamma_G^S(n)\geq a_n.
\end{equation}
\end{prop}

If the algebra $\mathcal{L}$ is infinite dimensional, then  the growth  of $\mathcal{A}$ is at least
$e^{\sqrt{n}}$, and therefore the growth of the group is at least $e^{\sqrt{n}}$. The Lie algebras approach
was used to show that  Gap Conjecture holds for residually nilpotent
groups~\cite{lubotzky_mann:polyn91}.

 Observe that in order to have an example of a residually-$p$ finite group  (i.e.\ a group approximated by finite
 $p$-groups) whose growth is exactly $e^{\sqrt{n}}$, the ranks of the consecutive quotients
 $G_n/G_{n+1}$ must be uniformly bounded (i.e., the group $G$ has to have \emph{bounded width}). But this condition is not enough. For instance, the group $\mathcal{G}$ has  finite width, as it was proved in~\cite{bartholdi-g.Lie00}, but its growth is bounded from below by
  $e^{n^{0.51}}$.  The growth of Lie algebra $\mathcal{L}$ associated
  with the Gupta-Sidki $3$-group $\mathcal{S}$ is linear, which implies that the group $\mathcal{S}$  has
  growth at least $e^{n^{2/3}}$ (it is not known yet if Gupta-Sidki $p$-groups have intermediate growth  or
  not).
Inequality~(\ref{gap4}) also gives a way to prove that a group has uniformly exponential growth
\cite{grigorch:hilbert,bartholdi-g.Lie00}.

\medskip
Recall that a group $G$ is said to have \emph{uniformly exponential} growth if
\begin{equation}\label{unif}
\kappa_{\ast}=\inf_A \kappa_A >1,
\end{equation}
where $\kappa_A$  denotes the base of exponential growth with respect to the system of generators $A$
($\kappa$ is defined by relation (\ref{growthexp})), and the infimum is taken over all finite systems of
generators. An immediate corollary of Proposition~\ref{growth1} is

\begin{cor}  Assume that the Lie algebra $\mathcal{L}$ associated to the group $G$ has  exponential
growth. Then $G$ has uniformly exponential growth.
\end{cor}

It follows that Golod-Shafarevich groups~\cite{golod:p-groups,golod_s:class_field_tower} have uniformly
exponential growth.

\section{Asymptotic invariants of probabilistic and analytis nature and corresponding gap type conjectures}\label{probabil}

In this section we discuss the relation between  group growth  and asymptotic behavior of random walks on a
group.  At the end of it we formulate gap type conjectures related to the asymptotic characteristics of
random walks, and discuss their relation with the growth Gap Conjecture from section~\ref{gapconj}.

The ICM paper of A.~Erschler~\cite{erschler:ICM11}, which we recommend to the reader, contains  important
material related to the subject of random walks and growth. Also we recommend the book of
W.~Woess~\cite{woess:rw}, the paper of Kaimanovich and Vershik~\cite{kaiman-ver:random83}, the
article~\cite{bendikov_pitsaur:10}, and the unfinished manuscript ``A survey on the relationships between
volume growth, isoperimetry, and the behavior of simple random walk on Cayley graphs, with examples'' of
C.~Pittet and L.~Saloff-Coste that can be downloaded from \url{http://www.math.cornell.edu/~lsc/lau.html}.

Let $G$ be a finitely generated group and $\mu$ be a probability measure on $G$ whose support $A$ generates
the group. Consider a random walk $(G,\mu)$ on $G$ which starts at the identity element $e$ and the
transitions $g\rightarrow ga$ take place with probability $\mu(a)$. Let $P(n)=P^{(n)}_{e,e}$ be the
probability of return after $n$ steps. Observe that  $P(n)=\mu^{\ast n}(e)$, where  $\mu^{\ast n}$ denotes
the $n$th convolution of $\mu$. In the case of a symmetric measure (i.e.\ when $ \mu(a)=\mu(a^{-1})$, for every
$a\in A$) the inequality
\begin{equation}\label{prob}
\frac{1}{P(2n)}\leq \gamma_G^A(2n)
\end{equation}
holds, since the maximal mass of $\mu^{\ast 2n}$ is concentrated at the identity element~$e$ \cite{woess:rw}.
In (\ref{prob}) the probability $P(n)$ is evaluated only for even values of $n$  because for odd values it
can be zero (this happens when the identity element of the group can not be expressed as a product of an odd
number of generators). From now on, when we discuss the rate of decay  of the probabilities  $P(n)$ as
$n\to \infty$, we assume that the argument $n$ takes only even values. We will use the comparison  $\preceq$
of the rate of decay of $P(n)$, or the rate of growth when $n \to \infty$ of some other functions that will
be introduced later, in the sense of the definition given in section~\ref{prelim}.

The rate of decay  of the probabilities  $P(n)$ can obey a power law (of the type $n^{\alpha}, \alpha < 0$),
be exponential (of type $\lambda^n, 0< \lambda <1$), or can be  intermediate between the two. A
power low holds if and only if the group  has polynomial growth. This follows from  a combination of Gromov's
theorem and results of Varopoulos on random walks  (\cite{Varopoulos:icm91} and
\cite[Theorem VI.5.1 on p. 84]{varopoulos_salcoul:groups92}).

An important characteristic  of random walks, introduced by Kesten~\cite{kesten:symmetric}, is the spectral
radius defined by  relation~(\ref{spectral}), which in the case of symmetric measures coincides with the norm of
the Markov operator
\[Mf(x)=\sum_{g\in G} \mu(g)f(xg).\]
Observe that this is also the operator given by right convolution with a measure $\mu$ acting on $l^2(G)$.
By Kesten's criterion a group is amenable if and only if for some (``some'' can be replaced by ``every'')
symmetric measure $\mu$ whose support generates a group, the spectral radius takes its maximal possible value
$r=1$~\cite{kesten:amenab}. Therefore amenable groups have subexponential rate of decay of return
probabilities, and the rate of decay  is exponential in the case of nonamenable groups.

For groups of exponential growth, the rate of decay  is not slower than $e^{-\sqrt[3]{n}}$. In other words, the
 upper bound
\begin{equation}P(n)\preceq \label{prob1} e^{-\sqrt[3]{n}}
\end{equation}
holds~\cite{Varopoulos:icm91}. This result cannot be improved, as there are groups of exponential growth
 for which the upper bound  (\ref{prob1})  is sharp (for
instance, the lamplighter group $\mathcal{L}= \mathbb{Z}_2\wr \mathbb{Z}$ or the Baumslag-Solitar solvable
groups $BS(1,n), n\geq 2$).

Because of inequality (\ref{prob}), if a bound of the type
\begin{equation}\label{prob2}P(n)\preceq e^{-n^{\alpha}}
\end{equation}
holds, then
\begin{equation}\label{growth5}\gamma_G(n)\succeq  e^{n^{\alpha}}.
\end{equation}
On the other hand if (\ref{growth5}) holds then
\begin{equation}\label{prob3}P(n)\preceq e^{-n^{\frac{\alpha}{\alpha+2}}}
\end{equation}
\cite{Varopoulos:icm91,woess:rw,bendikov_pitsaur:10}.

This leads  to the following natural questions. What
are the slowest and  fastest  rates of decay   of probabilities $P(n)$ for groups of intermediate growth? Is
it of type $e^{-n^\alpha}$ (respectively  $e^{-n^\beta}$) for some positive $\alpha$ and $\beta$? What are
the values of $\alpha$ and $\beta$?  The values $1/2$, $1/3$ and $1/5$  are the first candidates for these
numbers. Is there a group of intermediate growth with rate of decay  $P(n)\succ e^{-n^{1/3}}$? If the rate of
decay  of $P(n)$ for a group of intermediate growth cannot be slower than $e^{-n^{1/3}}$ then the weak
version of Gap Conjecture~\ref{c:gap1} holds with  parameter $1/3$. Later we will formulate a conjecture
related to the above discussion.
\medskip

A new  approach for obtaining lower bounds for    growth based on  techniques of random walks is developed by
A.~Erschler  in~\cite{erschler:boundary04,erschler:critical05}.  Without getting into details, let us briefly
outline some features of her approach. For a  random walk  given by  pair $(G,\mu)$  the \emph{entropy}
$h=h(G,\mu)$ and the \emph{drift} (or the \emph{rate of escape}) $l=l(G,\mu)$ are defined as
\begin{equation}
h=\lim_{n\to \infty} \frac{H(n)}{n},
\end{equation}
where $H(n)=H(\mu^{\ast n})$ and $H(\mu)=-\sum_{g\in G}\mu(g)\log \mu(g)$ is Shannon entropy, and
\begin{equation}
l=\lim_{n\to \infty}\frac{L(n)}{n},
\end{equation}
where $L(n)=\sum_{g\in G}|g|\mu^{\ast n} (g)$ is the expectation of the length $|g|$ of a random element at
the $n$th moment of the random walk (the length $|g|$ is considered  with respect to the system of generators
given by the support of $\mu$). By the Guivarc'h  inequality~\cite{guivarc'h:inequl80}, in the case of symmetric
measure with finite support (or more generally with finite first moment  $\sum_{g\in G_{\omega}}|g|\mu(g)$),
the numbers $h$, $l$ and $\kappa$ (the base of exponential growth defined by~(\ref{growthexp}) ) are related as
\begin{equation}\label{guivarc'h}
h\leq l\kappa.
\end{equation}
Therefore the equality $l=0$  implies $h=0$.

An important notion due to Furstenberg is  the Poisson boundary (which we will call the Poisson-Furstenberg
 boundary). It
is a triple $(G,\mathcal{B},\nu)$   consisting of a $G$-space $(\mathcal{B},\nu)$ with a $\mu$-stationary
probability measure $\nu$ (i.e.\ $\mu\ast\nu=\nu$). This boundary describes  space of \emph{bounded}
$\mu$-harmonic functions:
\begin{equation}
f(g)=\int_{\mathcal{B}}\phi(gx) d\nu(x)
\end{equation}
(the Poisson integral). The left hand side of the last equality takes values in the space of  bounded
$\mu$-harmonic functions while $\phi$  belongs to the space $L^{\infty}(\mathcal{B})$ (see
\cite{kaiman-ver:random83} for details).  The Liouville property of a group (more precisely, of a pair
$(G,\mu)$) is that every bounded $\mu$-harmonic function is constant; this property is equivalent to the
triviality of the Poisson-Furstenberg boundary.

The entropy criterion due to  Avez-Derriennic-Kaimanovich-Vershik ~\cite{kaiman-ver:random83} states the
following. Let $G$ be a countable group and  $\mu$ a probability measure on $G$ with finite entropy $H(\mu)$.
Under this assumption the Poisson-Furstenberg boundary is trivial if and only if the entropy $h$ of the
random walk is equal to zero.

If $\mu$ is symmetric and has  finite first moment  with respect to some (and hence with respect to every)
word metric on $G$, then the entropy $h$ is positive if and only if the rate of escape  $l$ of the random walk
determined by  $(G, \mu)$ is positive. In one direction this follows from the Guivarc'h inequality
(\ref{guivarc'h}). The converse was proved by Varopoulos ~\cite{Varopoulos:icm91} for finitely supported
measures and then extended by Karlsson and Ledrappier  to  the case of a measure with finite first moment
\cite{karlsson_ledr:drift07}.  It is known that for a group of intermediate growth and a measure with finite
first moment, the entropy  is zero; therefore the drift is also zero and the Poisson-Furstenberg boundary is
trivial. The vanishing of the entropy, and hence the triviality of the boundary, can easely be deduced for
instance from inequality~(18) in~\cite{kaiman-ver:random83}. In the case of a nonsymmetric measure, the drift $l$
can be nonzero, even on a group of polynomial growth (for instance, for a $(p,1-p), p<1/2$ random walk on
$\mathbb Z$). But, in the case of a symmetric measure, it is zero if and only if $h=0$. Therefore, in the case of
groups of subexponential growth  and symmetric measures with finite first moment $h=l=0$ and the functions
$H(n)$ and $L(n)$ grow sublinearly. The Poisson-Furstenberg boundary is also trivial for each group of
polynomial growth and any measure $\mu$; this follows from Gromov's result on groups of polynomial growth
and the theorem of Dynkin and Malyutov concerning Martin boundaries of nilpotent groups
(from which the triviality of
the Poisson-Furstenberg boundary follows in this case)~\cite{dynkin_malyut:61} (see also the work of
G.~Margulis~\cite{margulis:nilpotent66}).

\medskip
To obtain a lower bound for growth for some  groups from the family  $G_{\omega}$ and  to obtain  new results
about the Poisson-Furstenberg boundary of random walks on groups of intermediate growth,
in \cite{erschler:boundary04} Erschler  introduced   the so called    ``strong condition''
  $(\ast)$  for some type actions on the interval $(0,1]$. She  proved
 that if a group $G$  satisfies condition $(\ast)$  and  the group $germ(G)$ of germs of
 $G$ (also defined in~\cite{erschler:boundary04}) satisfies some extra condition, then $G$
 admits a symmetric probabilistic measure $\mu$  with finite entropy $H(\mu)$ and nontrivial
Poisson-Furstenberg boundary.
These conditions are satisfied for all groups $G_{\omega}, \omega \in \Omega_4$, where the set
$\Omega_4=\Omega_1 \setminus \Omega_0$ consists of the sequences containing only two symbols from
$\{0,1,2\}$, with each of them occuring infinitely many times  in the sequence.  Recall that by
Theorem~\ref{growth2} all such groups have intermediate growth. The groups $G_{\omega}, \omega \in \Omega_4$
are  first examples of groups of intermediate growth possessing a symmetric measure with nontrivial
Poisson-Furstenberg boundary. Also for all these groups  and any $ \epsilon >0$, a lower bound on growth of
type
\begin{equation}
\exp{\frac{n}{\log^{2+\epsilon}n}}\preceq \gamma_{\mathcal G_{\omega}}(n)
\end{equation}
holds~\cite{erschler:boundary04}. The proof of this result uses the existence of a special element $g\in G$
of infinite order. It is based on the combination of facts of existence of  measure with non-trivial Poisson
boundary and the analogue of inequality~(\ref{guivarc'h}) for measures with infinite first moment.

Interesting results  concerning  growth and triviality of the Poisson-Furstenberg boundary were obtained by
Karlsson, Ledrappier and Erschler~\cite{karlsson_ledr:drift07,erschler_karl:homom10}. These results led to
upper bounds on growth of $H(n)$ and $L(n)$ for the groups $G_{\omega}$, and show  that, under  certain
conditions, non vanishing of the drift implies that the group is indicable (i.e.\  existence of surjective
homomorphism onto $\mathbb Z$).

Let us mention now a new method of studying of asymptotic properties of self-similar groups discovered by
L.~Barthlodi and B.~Virag  in~\cite{bartholdi_v:amenab}.  It received further development  in the paper of
V.~Kaimanovich ~\cite{kaiman:munchhausen}, where it was called the \emph{ Munchhausen trick},  and in the papers
of Bartholdi, Kaimanovich and Nekrashevych~\cite{bartholdi_kainekr:10}, and of Amir, Angel and
Virag~\cite{amir_anvir:linear09}. In~\cite{kaiman:munchhausen} the entropy arguments were used to prove
amenability, the notion of  a self-similar measure was introduced, and the map $\psi$ in the space of
probabilistic measures on a self-similar group was defined. It  allows one to describe self-similar measures
as  fixed points of this map.  The relation between   $\psi$  and the classical tool of  linear algebra known
as  the Schur complement was established in~\cite{ grigorchuk-n:schur}.

The Munchhausen trick has been  used to prove amenability of certain self-similar groups of exponential
growth~\cite{bartholdi_v:amenab,bartholdi_kainekr:10,amir_anvir:linear09}. For the first time, this method was
applied to prove the amenability of the group $B$, named ``Basilica'', which was introduced in the paper of
A.~Zuk and the author~\cite{grigorch_z:basilica}, and can be
 defined as the \emph{iterated monodromy group} of the polynomial $z^2-1$, or alternatively as the group
generated by the automaton $\mathcal{A}_{852}$ from the \emph{Atlas} of self-similar complexity $(2,3)$ groups
~\cite{bondarenko_gkmnss:full_clas32}. Observe that  $B$ has exponential growth. The amenability of $B$
  allows us to separate the class $AG$ from the class $SG$ of subexponentially amenable groups that was
mentioned in section~\ref{amenab}. It was  originally defined  in~\cite{grigorch:example}
 (where the question about the possible coincidence of classes $AG$ and $SG$ was raised). It is currently unknown
whether the cardinality of the set $AG\setminus SG$ is the cardinality of the continuum or not.

Unfortunately, the Munchausen trick has not been used so far to obtain new information on growth of groups.
But it
 was used in~\cite{bartholdi_v:amenab,kaiman:munchhausen,bartholdi_kainekr:10} to obtain interesting results
 about the rate of growth of the functions  $H(n)$,  $L(n)$, and the rate of decay  of $P(n)$.
For instance, in the case of
 the Basilica group, $P(n)\succeq e^{-n^{2/3}}$~\cite{bartholdi_v:amenab}, while for the group of intermediate
growth $\mathcal{G}$, the lower bound is $P(n)\succeq e^{-n^{1/2-\epsilon}}$ for any  positive number $\epsilon$
(this follows for instance from the results of the paper by Bartholdi, Kaimanovich and
Nekrasevych~\cite{bartholdi_kainekr:10}).
 Using inequalities (\ref{prob1}) ) and (\ref{prob3}) we obtain the following estimates:
 \[e^{-n^{2/3}}\preceq P_{B}(n) \preceq e^{-n^{1/3}},\]
\[e^{-n^{1/2-\epsilon}}\preceq P_{\mathcal{G}}(n) \preceq e^{-n^{1/5}}.\]
It would be interesting to find the asymptotics of the rate of decay  of $P(n)$ for each of the groups
$\mathcal{B}$ and $\mathcal{G}$ (as well as the rate of growth of the functions $H(n)$  and $L(n)$).

\medskip
Behind the idea of the Munchausen trick is the conversion of the self-similarity of the group into
self-similarity of the random walk on the group. Let $G$ be a self-similar group acting level transitively on
a $d$-regular rooted tree $T_d$, and let $\mu$ be a probability measure on $G$. Denote by $H=st_G(x)$ the
stabilizer of a vertex $x$ on the first level. Then $[G:H]=d$. Let $p_x\:H\rightarrow G$ be the projection
homomorphism of $H$ on the subtree $T_x$ with root $x$  ($p_x(g)=g_x$, where $g_x$ is the section of $g$ at
vertex $x$), and let $\mu_H$ be the probability distribution on $H$ given by the probability of the first hit
of $H$ by a random walk  on $G$ determined by $\mu$. Denote by $\mu_x = ({p_x})_{\ast}(\mu_H)$  the image
of the measure $\mu_H$ under the  projection $p_x$.

 \begin{defn} A measure $\mu$ is called \emph{self-similar}  if for some vertex $x$ of the first level
 \[\mu_x=(1-\lambda)\delta_e+\lambda \mu\]
 for some $\lambda, 0<\lambda<1$.
\end{defn}

\begin{thm} (\cite{bartholdi_v:amenab,kaiman:munchhausen}) If  a self-similar group $G$ possess  a self-similar
symmetric probability measure $\mu$ with finite entropy and contracting coefficient $\lambda, 0 < \lambda
<1$, then the entropy $h$ of the corresponding random walk on $G$ is zero and therefore the group $G$ is
amenable.
\end{thm}

For instance, for  $\mathcal G$ the measure $\mu=\frac{4}{7}a+\frac{1}{7}(b+c+d)$ is self-similar with
contracting coefficient $\lambda=1/2$ ~\cite{kaiman:munchhausen}.

\medskip
Let $G$ be an amenable group  with a finite generating system $A$. Then  one can associate with $(G,A)$ a
function
\[
 F_G(n)=\min \{n: \text{ there is a finite subset\\
  } F \subset G
 \text{ s.t. } \frac{|F \bigtriangleup aF|}{|F|} < \frac{1}{n}, \forall a \in A\},
\]
which is  called the F\"{o}lner function of $G$ with respect to $A$ (because of F\"{o}lner's criterion of
amenability~\cite{folner:mean}). This function was introduced by Vershik in the Appendix to the Russian
edition of Greenleaf's book on amenability~\cite{greenleaf:means}. The growth type of this function does not
depend on the generating set. By a result of Coulhon and Saloff-Coste~\cite{coulhon_saloff:isoperim93}, the
growth function and the F\"{o}lner function are related by the inequality
\begin{equation}\label{coulhon}
\gamma(n)\preceq F(n).
\end{equation}
On the other hand, for groups with $\gamma(n)\preceq e^{n^{\alpha}}$,  the F\"{o}lner function can be
estimated as in~\cite[Lemma 3.1]{erschler:piecewise06}
\begin{equation}\label{felner}
F(n)\preceq e^{n^{\frac{\alpha}{1-\alpha}}}.
\end{equation}

It is proved by Erschler that for an arbitrary function $f\:\mathbb N\rightarrow \mathbb N$  there is a group
of intermediate growth with F\"{o}lner function $F(n)\succeq f(n)$.  The method of construction of such
groups is based on the ``oscillation'' type techniques that we discuss briefly in the last section. Very
interesting results about the asymptotics of random walks on Schreier graphs associated with finitely
generated groups, and in particular with  $\mathcal G$, are obtained by Erschler
in~\cite{erschler:critical05}. 
We shall reformulate the following question soon in the form of a conjecture.

\begin{prob}  \label{erschler}
Can a F\"{o}lner function  grow strictly slower than the exponential function but faster than any polynomial?
\end{prob}

In view of the inequality (\ref{felner}), the Gap Conjecture holds if the answer to the last problem is
negative.

\bigskip
The next conjecture (consisting of three subconjectures) was formulated by P.~Pansu and the author in 2000 in
an unpublished note. Let $G$ be an  amenable finitely generated group, $\mu$ a  symmetric probability
measure with finite support that generates $G$, $M$ a Markov operator of the associated  random walk on $G$
given by $\mu$, and $P(n)$ the probability of return after $n$ steps for this random walk. Let  $F(n)$ be the
F\"{o}lner function,
and $\mathcal N(\lambda)$ be the spectral density  defined by the relation
\[\mathcal N (\lambda)= tr_{vN}(\chi_{(-\infty,\lambda]}(\Delta),\]
where $tr_{vN}$ is the von Neumann trace defined on the von Neumann algebra $\mathcal N(G)$ of $G$, generated
by the right regular representation,  $\Delta=I-M \in \mathcal N(G)$ is the discrete Laplace operator on
$G$, and $\chi_{(-\infty,\lambda]}(\Delta)$ is the projection obtained by application of the characteristic
function $\chi_{(-\infty,\lambda]}$ to  $\Delta$. We are interested in the asymptotic behavior of $P(n)$ and
$F(n)$ when $n \to \infty$ and of $\mathcal N(\lambda)$ when $\lambda \to 0+$.  Their asymptotic behavior
does not depend on the choice of the measure $\mu$ ~\cite{bendikov_pitsaur:10}.

\begin{conj}\label{pansu}\quad\\[-\baselineskip]
\begin{itemize}
\item[(i)] (Gap conjecture for the heat kernel). The function  $P(n)$ is either of power rate of decay or satisfies
\[P(n)\preceq e^{-\sqrt[3]{n}}. \]

\item[(ii)] (Gap Conjecture for the F\"{o}lner function)  The F\"{o}lner function $F(n)$ has either  polynomial
growth or the growth is at least exponential.

\item[(iii)]  (Gap Conjecture for the spectral density) The spectral density $\mathcal N (\lambda)$ either has power
decay of type $\lambda^{d/2}$ for some $d \in \mathbb N$ when $\lambda \to 0$ or
\[ \mathcal{N}(\lambda) \preceq e^{-\frac{1}{\sqrt[4]{\lambda}}}. \]
\end{itemize}
\end{conj}

Let us also  formulate a  modified version of the previous conjecture. We guess that for each of three
conjectures stated below there is a number $\beta>0$ for which it holds.

\begin{conj} (Gap Conjectures with parameter $\beta,\beta >0$) \label{pansu1}\\[-\baselineskip]
\begin{itemize}
\item[(i)]  The function $P(n)$ is either of power rate of decay or satisfies
\[P(n)\preceq e^{-n^{\beta}}. \]

\item[(ii)]   The F\"{o}lner function $F(n)$ has either  polynomial growth or the growth is not less than
$e^{n^{\beta}}$.

\item[(iii)]  The spectral density $\mathcal N (\lambda)$ either has power decay of type $\lambda^{d/2}$ for some $d
\in \mathbb N$ or
\[ \mathcal N(\lambda) \preceq e^{-\lambda^{-\beta}}. \]
\end{itemize}
\end{conj}

Each of the above alternatives separates the case of polynomial growth from the intermediate growth case.

Perhaps  gap type conjectures can  also be formulated in a reasonable way for some other asymptotic
characteristics of groups such as the isoperimetric profile, $L^2$-isoperimetric profile,  entropy function
$H(n)$, drift function $L(n)$,  etc.

There are  relations between all of the conjectures stated here to the growth Gap Conjecture from
section~\ref{gapconj} (which we will call here Growth Gap Conjecture), and to its generalization,
Conjecture~\ref{c:gap2}). For instance, using  inequality (\ref{felner})  we conclude that the Gap Conjecture
for the F\"{o}lner function
 implies  the Growth Gap Conjecture. On the other hand,  inequality~\ref{coulhon} shows that the
Growth Gap Conjecture with parameter $\beta$
 implies the Gap Conjecture with parameter $\beta$ for the F\"{o}lner function.
Therefore the Weak Gap Conjecture
 for growth is equivalent to the Weak Gap Conjecture for the F\"{o}lner function
(the latter is formulated similarly to the conjecture  \ref{c:gap3}).

  The Growth Gap Conjecture 
 with parameter $\alpha$
 implies the Gap Conjecture with parameter $\frac{\alpha}{\alpha+2}$ for the return probabilities $P(n)$
(because of inequality (\ref{prob3})) etc.
 It would be interesting to find the relation between all stated conjectures.

The  results of~\cite{grigorch:hilbert,lubotzky_mann:polyn91,wilson:gap11} and some statements from this
article  provide the first  classes of groups for which the Gap Conjectures of the given type hold. For
instance, the Gap Conjecture with parameter $\frac{1}{5}$ for return probabilities $P(n)$, and  with parameter
$\frac{1}{2}$ for the F\"{o}lner function $F(n)$  hold  for residually supersolvable groups.

\section{Inverse orbit growth and examples with explicit growth} \label{inter3}

Until recently there was no exact computation of the intermediate growth in the sense of the Schwarz-Milnor
equivalence. The first such examples were produced  recently by L.~Bartholdi and A.~Erschler
\cite{bartholdi_Ers:growth10}.  The idea is very nice and we shall explain it briefly.

Observe that the notion of growth can be defined for transitive group actions. Namely, if a finitely
generated group $G$ with a system of generators $A$ acts transitively  on a set $X$ and a base point $x\in X$
is selected, then the growth function $\gamma^A_{X,x}(n)$ counts the number of points in $X$ that can be
reached from $x$ by consecutive applications of at most $n$ elements from set $A\cup A^{-1}$. The growth type
of this function (in the sense of the equivalence $\sim$) does not depend on the choice of $x$. The triple
$(G,X,x)$ can be encoded by the \emph{Schreier graph} (or the \emph{graph of the action}) $\Gamma$, with set
of vertices $X$ and a set of oriented edges (labeled by the elements of $A$) consisting of  pairs $(x,a(x)),
x\in X, a\in A$. The growth function of the action is the same as the growth function of the graph $\Gamma$,
which is a $2d$-regular graph (viewed as non-oriented graph).

\begin{figure}[!htb]
\ifColorSection\includegraphics[width=.95\textwidth]{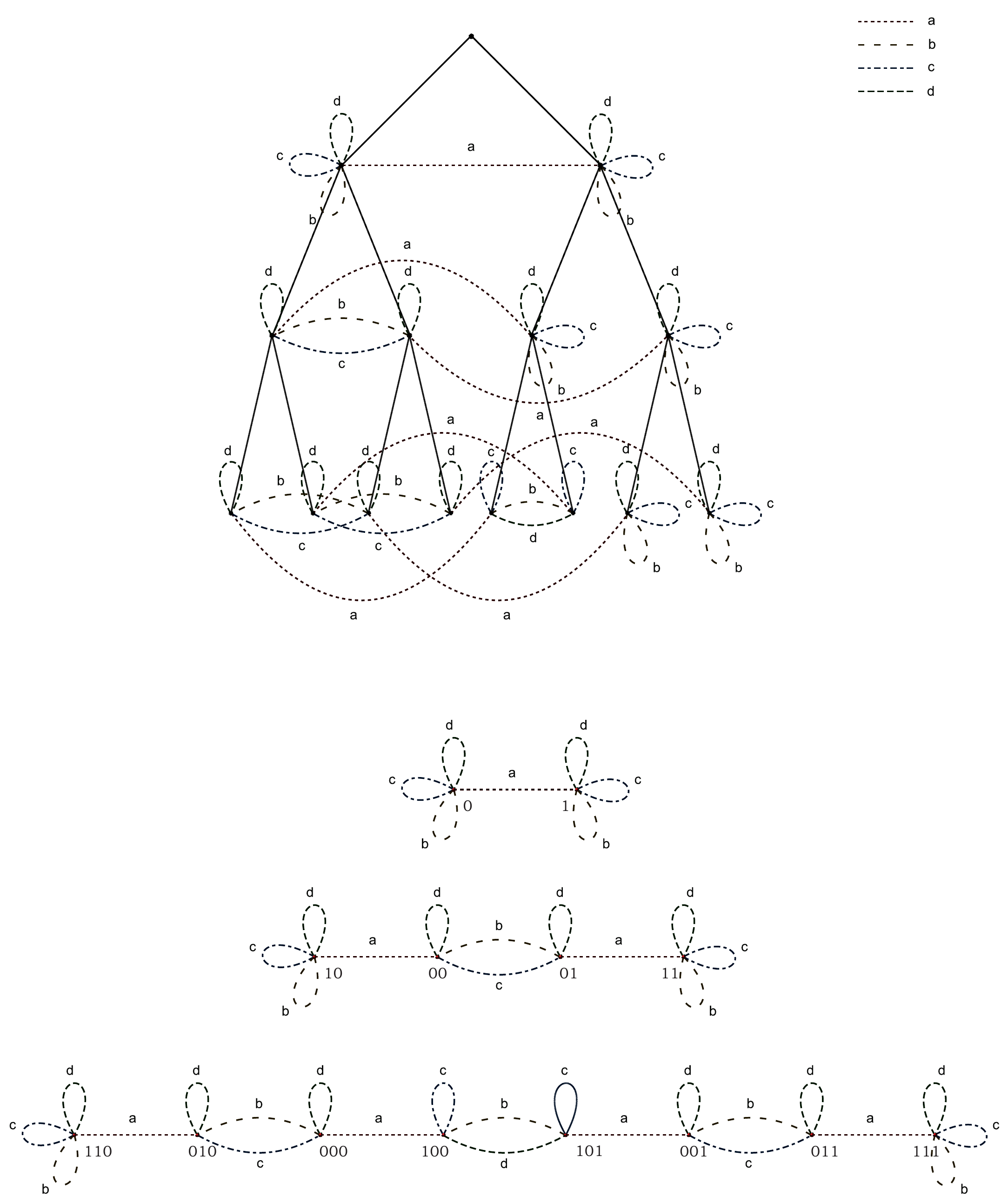}
\else\includegraphics[width=.95\textwidth]{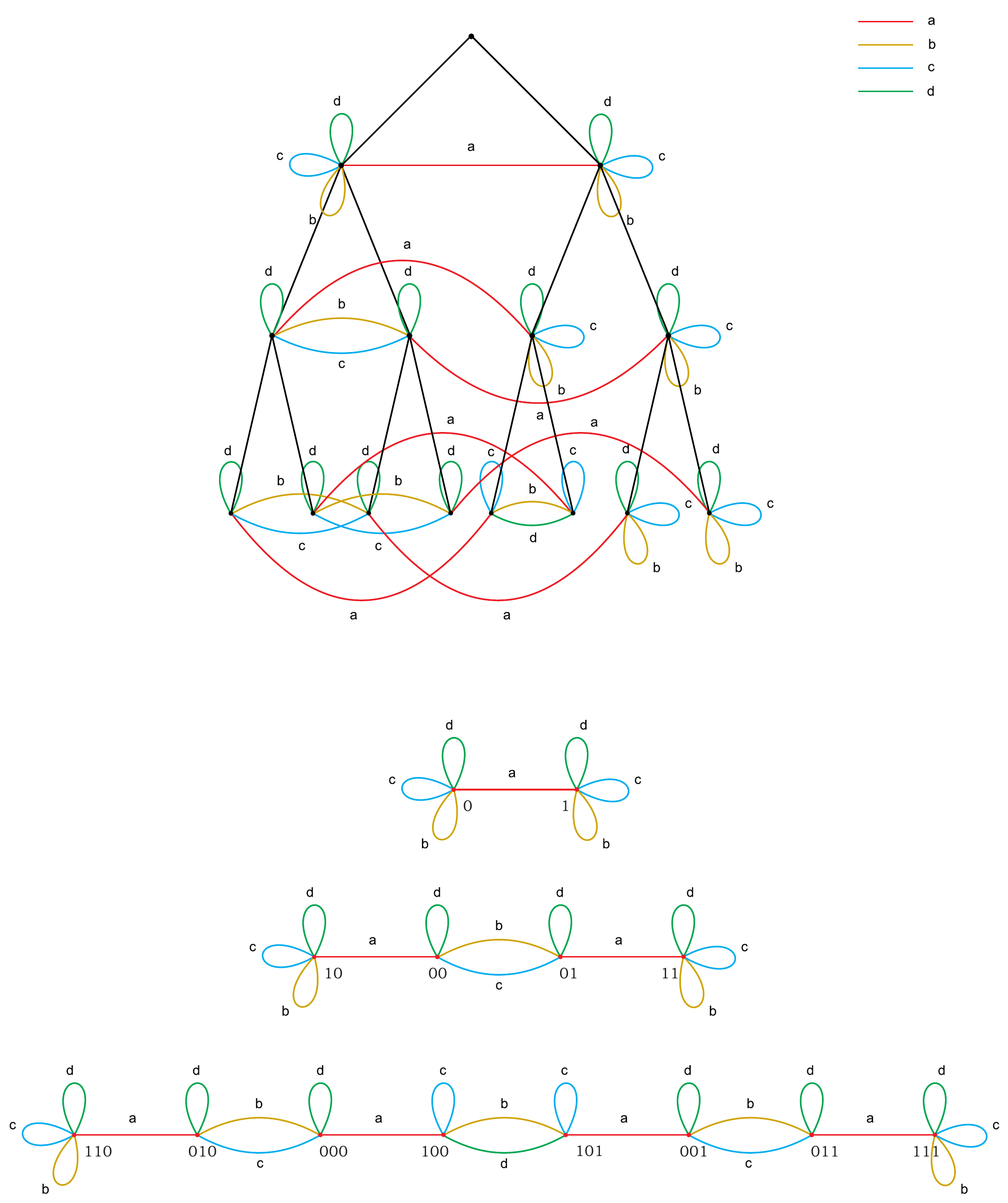}\fi
\caption{The action of $\mathcal G$ on the first three levels of the tree.
\SeeColorFig{C:gaction3}}
\label{f:gaction3}
\end{figure}

\begin{figure}[!ht]
\ifColorSection\includegraphics[width=.95\textwidth]{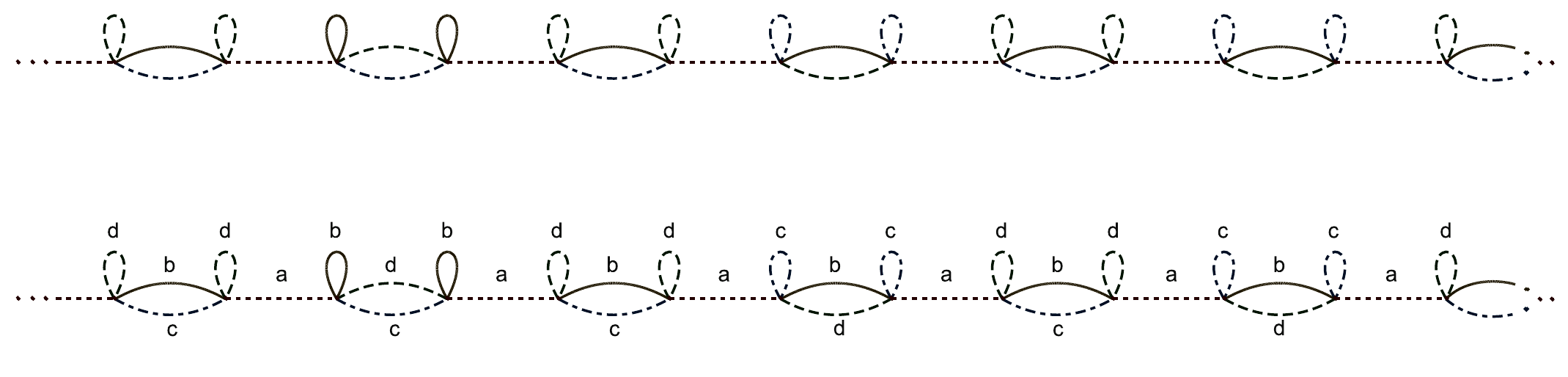}
\else\includegraphics[width=.95\textwidth]{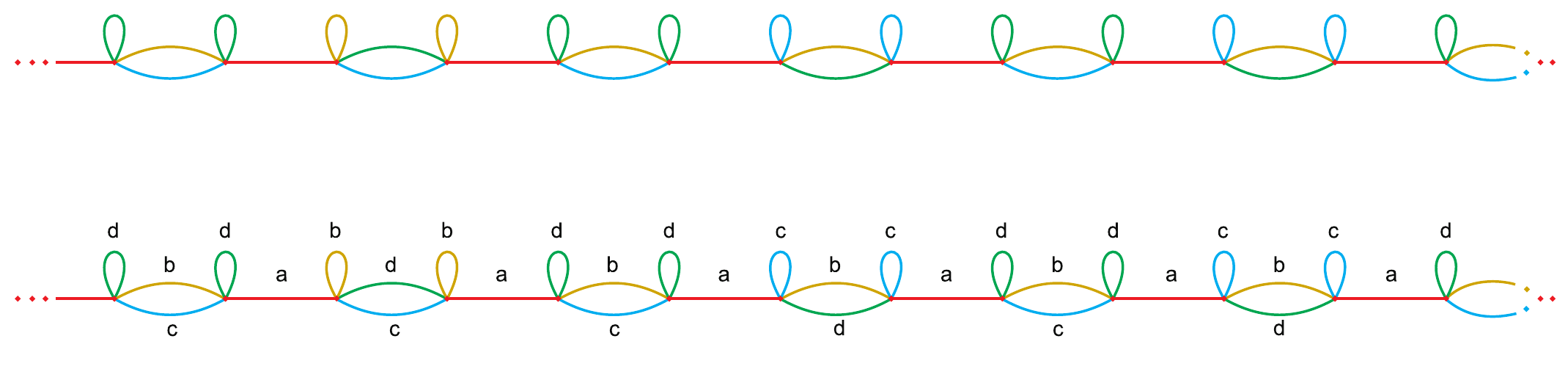}\fi
 \caption{Typical Schreier graph of the boundary action.
\SeeColorFig{C:grigfig2}}
 \label{fig2}
\end{figure}

It is easy  to construct actions with intermediate growth between polynomial and exponential. An interesting
topic is  the study of growth of Schreier graphs associated with actions of self-similar groups on corresponding
rooted trees and their boundaries. The graphs we allude to here are of the form $\Gamma= \Gamma(G,H,A)$, where
$G$ is a finitely generated self-similar group with generating set $A$, and $H=St_G(\xi),$ where $\xi$ is a
point on the boundary of the tree. Such graphs are isomorphic to  the corresponding graphs of the action of
the group on the orbit of the base  point (namely $\xi$).  For contracting groups they have  polynomial
growth, which can be of  fractional and even irrational degree \cite{bartholdi_g:spectrum,bondarenko:PhD}.
There are examples of actions with quite exotic intermediate behavior like $n^{\log^m n}$ for some $m>0$
\cite{benjamini_h:omega_per_graphs,grigorchuk-s:hanoi-cr,bondarenko:growth11,bondarenko-3:family11}.  At the
same time there are examples with quite regular polynomial type orbit growth.  For instance, for the group
$\mathcal G$ the orbit growth is linear. Schreier graphs of the action of $\mathcal G$ on the first three
levels of the binary tree are shown in Figure~\ref{f:gaction3}
and the infinite graph of the action on the orbit of a typical point of the boundary is given in
Figure~\ref{fig2} (drawn in two versions: with labels and without).

In~\cite{bartholdi_Ers:growth10,bartholdi_ersch:given(2)11} the notion of \emph{inverted orbit growth
function} has been introduced. Let us briefly  explain the idea and  formulate some of the results of these
articles. Let $G$ be a group  acting on the right on a set $X$,  $S$ be a generating set for  $G$ (viewed as
monoid), and let $x\in X$ be a base point. Denote by $S^{\ast}$ the set of (finite) words over the alphabet
$S$. For a word $w = w_1 . . .w_l \in S^{\ast}$, its inverted orbit is
\[\mathcal O(w) = \{x, xw_l,xw_{l-1}w_l, xw_1\dots w_{l-1}w_l\},
\]
and the  inverted orbit growth of $w$ is $\delta(w) = |\mathcal O(w)|$. The inverted orbit growth function of
$G$ is the function
\[\Delta(n)=\Delta_{(G,X,x)}(n) = \max_{w}\{\delta(w) | |w|=n \}.\]
Clearly $\Delta_{(G,X,x)}(n)\preceq \gamma_{(G,X,x)}(n)$.

The notions of \emph{wreath product} and \emph{permutational wreath product} are standard in group theory.
Consider groups $A, G$ and a $G$-set $X$, such that $G$ acts on $X$ on the right. The wreath product $W = A
\wr_X G$ is a semidirect product of
 $\sum_X  A$  (a direct sum of copies of $A$ indexed by $X$) with $G$  acting on $\sum_X  A$ by automorphisms induced by the corresponding permutations  of  $X$. In other words, view  elements of $\sum_X  A$  as finitely supported functions $X\rightarrow A$. A left action of $G$ on $\sum_X  A$ by automorphisms is then defined by $(g f)(x) =
f(xg)$. There are two versions of (permutational) wreath products, restricted and unrestricted, and here we
use the \emph{restricted} version (because of the assumption of finiteness of the support of the elements
from the base group $\sum_X A$).

\medskip
In the next two theorems, proved in~\cite{bartholdi_Ers:growth10}, $X$ is a $\mathcal{G}$-orbit of an
arbitrary point of the boundary $\partial T$ of the binary rooted tree $T$ that does not belong to the orbit
of the point $1^{\infty}$. Observe that the orbit of $1^{\infty}$ consists of the sequences cofinal to
$1^{\infty}$, where two sequences are cofinal if they coincide starting with some coordinate, and more
generally, two points of the boundary of $T$ are in the same orbit if they are cofinal
\cite{bartholdi_g:spectrum,bartholdi_g:parabolic,grigorch:solved}. The graph of the action of $\mathcal{G}$
(after deletion of the labels) looks similar to the one shown in the last figure. Recall that the number
$\alpha_0 \approx 0.7674$  was defined above and is included in~(\ref{barth}).

\begin{thm} Consider the following sequence of groups: $K_0 = \mathbb{Z}/2\mathbb{Z}$, and
$K_{k+1} = K_k \wr_X \mathcal{G}$. Then every $K_k$ is a finitely generated infinite torsion group, with
growth function $$\gamma_{K_k} (n) \sim exp(n^{[1-(1-\alpha_0)^k)]}).$$
\end{thm}

Recall that the torsion free group of intermediate growth $\hat{\mathcal G}$, constructed in
\cite{grigorch:degrees85}, was mentioned in section~\ref{inter1} (theorem~\ref{GRI83}).

\begin{thm} Consider the following sequence of groups: $H_0 = \mathbb{Z}$, $H_1= \hat{\mathcal G}$ and $H_{k+1} =
H_k \wr_X  \mathcal{G}, k\geq 1$. Their growth functions satisfy
$$\gamma_{H_k} (n) \sim
\exp(log(n)n^{1-(1-\alpha_0)^k}).$$
\end{thm}

It is remarkable that not only two infinite series of groups of intermediate growth  with precisely computed
growth have been constructed, but also that the precise growth of the torsion free group $\hat{\mathcal G}$,
constructed more than 25 years ago, has finally been evaluated.

The proof of these two theorems is based on the calculation of the inverse orbit growth of the action of
$\mathcal{G}$ on $X$, which is $n^{\alpha_0}$ in this case. In the realization of  the first step of this
program the technique used previously by Bartholdi~\cite{bartholdi:growth} for improving the upper bound in
the case of the group $\mathcal{G}$ from the value $\log_{32}31$ used in~\cite{grigorch:degrees} to the value
$\alpha_0$ established in~\cite{bartholdi:growth,muchnik_p:growth} is explored again. The technique is based
on assigning positive weights to the canonical generators $a,b,c,d$  and finding the values that give the
best possible (for this approach) upper bound for growth.

The above ideas and some other tools (in particular, the dynamics of partially continuous self-maps of
simplices) has been used by Bartholdi and Erschler in~\cite{bartholdi_ersch:given(2)11} to construct a
 family of groups with growth of the type $e^{n^{\alpha}}$  with $\alpha$ belonging to
 the interval  $(\alpha_0, 1)$.  Moreover, they presented the following impressive result. Let $\eta_+\approx 2.4675$ be the positive root of $x^3 - x^2-2x-4$.

 \begin{thm}\label{ersch}
 Let  $f\:\mathbb R\rightarrow \mathbb R$ be a function satisfying
 \[ f(2x)\leq f(x)^2\leq f(\eta_{+}x)\]
 for all $x$ large enough. Then there exists a finitely generated group with growth equivalent to the growth
 of $f$.
 \end{thm}

This theorem   provides a large class of growth functions of finitely generated groups that ``fill'' the
 ``interval'' $[e^{n^{\alpha_0}}, e^n]$.


\section{Miscellaneous}\label{miscellaneous}

Recall that a group has uniformly exponential growth if it has exponential growth and moreover the number
$\kappa_{\ast}$ (the base of exponential growth)  defined by (\ref{unif}) is $> 1$.

In~\cite{gromov_panlaf:book81} Gromov raised the following  question.

\theoremstyle{plain}
\newtheorem*{GromovProblemII}{Gromov's Problem on Growth (II)}
\begin{GromovProblemII}
 Are there groups of exponential but not uniformly exponential growth?
\end{GromovProblemII}

Some preliminary results concerning uniformly exponential growth were obtained by P. de la Harpe
\cite{harpe:uniform02}. The problem  of Gromov (II) was solved by  J.S.~Wilson
~\cite{wilson:nonuniform,wilson:further04} by providing an example of such a group. A shorter
  solution was found later by Bartholdi  in~\cite{bartholdi:nonuniform}.
 Tree-like constructions and  techniques of self-similar groups naturally leading to such examples were explored by V.~Nekrashevych~\cite{nekrash:non-unifom10}. Similar to the intermediate growth case, all  known
  examples of groups of exponential but not uniformly exponential growth are based on the use of self-similar groups of branch type.

There are results that show that groups  within certain classes of groups of exponential growth are of
uniformly exponential growth. This holds, for instance, for hyperbolic groups
(M.~Koubi~\cite{koubi:unifom98}), one-relator groups, solvable groups, linear groups over fields of
characteristic 0 (A.~Eskin, S.~Mozes and H.~Oh~\cite{eskin_mozoh:uniform05}), subgroups of the mapping class
group which have exponential growth (J.~Mangahas~\cite{mangahas:unif08}), and some other groups and classes
of groups.

The fact that solvable groups of exponential growth have uniformly exponential growth was proved by D.~Osin
\cite{osin:entropy03} and he generalized this result to elementary amenable groups \cite{osin:elementary04}.
E.~Breuillard  gave a ``ping-pong'' type proof of Osin's result \cite{breuillard:uniform07}. A connection
between the so called ``slow growth'' and the Lehmer conjecture is another result in Breuillard's paper. The
fact that one relator groups of exponential growth have a uniformly exponential growth is proved by P. de la
Harpe and the author~\cite{grigorch_harpe:onerelator01}, and that a one relator group has either polynomial
growth or exponential growth is shown in the paper of T.~Ceccherini-Silberstein and the
author~\cite{ceccherini_grigorch:97}.

A second topic of this section is the \emph{oscillation} phenomenon that exists in the world of  growth of
finitely generated groups. The meaning of this is that there are groups whose growth in a certain sense may
oscillate between two types of growth. This was first discovered by the author in~\cite{grigorch:degrees} and
was further developed in his Habilitation~\cite{grigorch:habil}.   The goal achieved in that work was the
construction of a chain and an anti-chain of cardinality of the continuum in the space of growth degrees of
finitely generated groups. The ``trick'' used for this purpose can be described briefly as follows.

The groups $\mathcal G_{\omega}, \omega \in \Omega \setminus \Omega_1$ (recall that  $\Omega \setminus
\Omega_1$ consists of sequences that are constant at infinity) are virtually abelian,  while the rest of the
groups $\mathcal G_{\omega}, \omega \in \Omega_1$  have intermediate growth and the set $\{\mathcal
G_{\omega}, \omega \in \Omega_1\}$   has the property that if two sequences $\lambda,\mu \in  \omega \in
\Omega_1$ have the same prefix of length $n$ then the subgraphs with vertices in the balls $B_{\mathcal
G_{\lambda}}(2^{n-1})$ and $B_{\mathcal G_{\nu}}(2^{n-1})$ of radius $2^{n-1}$ with centers at the identity
elements in the Cayley graphs of these groups are isomorphic. It was suggested in~\cite{grigorch:degrees} to
replace the groups from the set $\{\mathcal G_{\omega}, \omega \in \Omega \setminus \Omega_1\}$ by the set of
accumulation points of $\{\mathcal G_{\omega}, \omega \in \Omega_1\}$ in the space of the Cayley graphs of
4-generated groups, supplied with a natural topology introduced in the article. Then each of deleted groups
$\{\mathcal G_{\omega}, \omega \in \Omega \setminus \Omega_1\}$ is replaced by some virtually metabelian
group (for which we will keep the same notation) of exponential growth and the modified set of groups
$\{\mathcal G_{\omega}, \omega \in \Omega\}$ becomes a closed subset in the space of groups homeomorphic to a
Cantor set.

Consider a sequence $\theta$ of the type

\begin{equation}\label{sequence}
\theta=(012)^{m_1}0^{k_1}(012)^{m_2}0^{k_2}\dots
\end{equation}
with a very fast growing sequence $m_1,k_1,m_,k_2,\dots,m_i,k_i,\dots$ of parameters. In view of the property
that the balls of radius $2^{n-1}$ coincide for the considered groups when they are determined by sequences
with the same prefix of length $n$, and the fact that the group $\mathcal G _{\omega}, \omega \in \Omega_1$
is (abstractly) commensurable to the group $\mathcal G_{\tau^n(\omega)}^{2^n}$ ($\tau$ the shift in space of
sequences), we  conclude that for initial values of $n$ when  $n \in [1,R(m_1)]$ ($R(m_1)$ determined by
$m_1$) the group $\mathcal G_{\theta}$ behaves as a group of intermediate growth (similar to $\mathcal
G=\mathcal G_{(012)^{\infty}}$), but then for larger values of $n, n \in [R(m_1)+1,R(m_1,k_1)]$ ($R(m_1,k_1)$
determined by $m_1,k_1$) it starts to grow as the group $\mathcal G_{0^{\infty}}^{2^{3m_1-1}}$ (i.e.
exponentially), but then again the growth slows down and behaves in intermediate fashion, etc. Taking
sequences of the type \eqref{sequence} but determined by various sequences of parameters, one can construct
a chain and an anti-chain of cardinality of the continuum  in the space of growth degrees of 3-generated
groups; this was done in~\cite{grigorch:degrees,grigorch:habil}.

The possibility of application of the oscillation technique is based  on the use of a topology in the space
of Cayley graphs (or what is the same, in the space of marked groups) introduced in \cite{grigorch:degrees}
(this topology is a relative to Chabauty topology known in the theory of locally compact groups and geometric
topology \cite{harpe:topics}).
 The space of marked groups is a compact totally disconnected space and one of the major
problems is to find its Cantor-Bendixson rank  (for more on this and related problems see
\cite{grigorch:solved}). Although the fact that this space has a nontrivial perfect core follows from the
result established by B.Neumann in 1937 (a construction of uncountably many 2-generated groups, up to
isomorphism),  there is a considerable interest in  finding Cantor subsets in the space of marked groups
consisting of interesting families of groups. The first such subset was identified in
\cite{grigorch:degrees}, which, with exception of a countable set, consists of groups of intermediate growth
$\mathcal G _{\omega}, \omega \in \Omega_1$ (its construction was described  in section \ref{inter1}).

Oscillation techniques received further development in the paper of Erschler \cite{erschler:piecewise06}
(which was already mentioned at the end of section~\ref{probabil}), where the notion of piecewise automatic
group was introduced. Using  this notion, Erschler constructed groups of intermediate growth with arbitrary
fast growth of the F\"{o}lner function. Moreover, the asymptotic entropy of a random walk on Erschler
groups can be arbitrarily close to a linear function, while at the same time, the Poisson boundary can be
trivial. Oscillation techniques, in combination with ideas from~\cite{bartholdi_Ers:growth10} were used also
by J.~Brieussel~\cite{brieussel:growth11}. The most recent result of ``oscillation character'' is due to
M.~Kassabov and I.~Pak~\cite{kassabov_pak:11}. They demonstrated a very unusual oscillation phenomenon for
groups of intermediate growth and their result is based on a new idea which we are going to explain briefly.

The groups $G_{\omega}$, as well as groups of branch type, act on spherically homogeneous rooted trees, that
is, trees $T_{\bar{k}}$ defined by sequences of integers $\bar{k}=k_1,k_2,\dots$, with $k_i\geq2$; these sequences
$\bar{k}$ are called the \emph{branch index} ($k_i$ is the branching number for the $i$th level of the tree).
Kassabov and Pak
suggested to modify this approach by considering actions on \emph{decorated} trees by attaching to some
levels of $T_{\bar{k}}$ finite subtrees with actions of suitably chosen finite groups  (each vertex of the
corresponding  level is decorated by the same structure). The sequence $\{F_i\}$ of attached groups has to
satisfy certain properties (in particular the groups need to be generated by four involutions, three of which
commute  as in  case of $\mathcal G_{\omega}$), but the main property is that the groups
$F_i$ must behave as expanders, in a certain sense. Namely, their Cayley graphs have to have
diameters $d_i$ growing as $i \to \infty $ as the logarithm of the size of the group, and for values of $n$
in the range $1\leq n\leq Cd_i$ ($C, 0 < C < 1$ some constant independent of $i$), the growth
functions $\gamma_i(n)$ have to  behave as the exponential function.

\bigskip
As it was already mentioned in the introduction, there are various types of asymptotic characteristics that
can be associated with algebraic objects. In addition to the  group  growth, and other characteristics
considered in section~\ref{probabil}, the \emph{subgroup} growth, the \emph{conjugacy} growth, the
\emph{geodesic} growth, and many other types of growth have been studied.  The subgroup growth was already
discussed a little bit and we refer the reader to the comprehensive book on this and other subjects,
by Lubotzky and Segal~\cite{lubotzky_segalbook} and the literature cited therein.

\medskip
The conjugacy growth counts the number of conjugacy classes of length $\leq n, n=1,2,\dots$, where the length
of a conjugacy class is the length of the shortest representative of this class. Interesting results on this
subject are obtained by I.~Babenko~\cite{babenko:closed88} (who was perhaps the first who introduced this
notion), M.~Coornaert and G.~Knieper, who studied the hyperbolic groups
case~\cite{Coornaert_Knieper:conjugacy02},  Breuillard and Coornaert~\cite{breuillard_corn:conjugacy10} (the
case of solvable groups), and M.~Hull and D.~Osin~\cite{hull_osin:conjugacy11} who showed that basically any
monotone function growing not faster than an exponential function is equivalent to the conjugacy growth
function. This is a far from complete list of papers and results on this subject (for a more complete list see
the literature in the cited papers).

\medskip
The geodesic growth of a pair $(G,A)$  (a group and a finite system of generators) is a growth of the
language of geodesic words over  alphabet $A\cup A^{-1}$   (i.e.\ words over an alphabet of generators that
represent geodesics in the Cayley graph $\Gamma(G,A)$ with the origin at the identity element). It can be
polynomial with respect to some system of generators but exponential with respect to other systems of
generators, and it is unclear if it can be intermediate between polynomial and exponential. This notion was
studied
in~\cite{cannon:combinatorial84,neumann_Shapiro:95,grigorch_nagnib:complete97,bridson_Bur_Eld_Sun:conjugacy11}.
The following question was discussed by M.~Shapiro and  the author  around 1993.

\begin{prob}
Are there pairs $(G,A)$ consisting of a group of intermediate  growth $G$ and a finite system of generators
$A$ with intermediate geodesic growth?
\end{prob}

All known groups of intermediate growth have exponential geodesic growth.  For instance for $\mathcal G$ this
follows from the fact that  the Schreier graph presented in figure~\ref{fig2}  has exponential geodesic growth
which is obvious.

The study of geodesic growth is a particular case of study of growth of formal languages. Such questions
originated in the work of Sch\"{u}tzenberger in the 1950s. The literature on this subject
related to group theory can be found
in\cite{bridson_gil:context02,Ceccher_Woess:context02,Ceccher_Woess:sensitivity03,gilman:languages05}.
Observe that for regular and for context-free languages, growth can be only polynomial or exponential
(R.~Incitti~\cite{incitti:growth01}, M.~Bridson and R.~Gilman \cite{bridson_gil:context02}), while the
intermediate type behavior is possible for indexed  languages (the next in the language hierarchy type of
languages after the context-free languages), as is shown in the note of A.~Machi and the
author~\cite{grigorch_Machi:indexed99}.

\medskip
There are interesting studies about growth of regular graphs. One of the first publications on this subject
is the article of V.~Trofimov~\cite{trofimov:graphs84} where, under certain conditions on the group of
automorphisms of the graph, the case of polynomial growth is studied.  In the last decade, the study of Schreier
graphs of finitely generated groups, and in particular of their growth and amenability properties has been
intensified. Some results were already mentioned in previous sections. We mention here the paper of
Bartholdi and the author~\cite{bartholdi_g:spectrum}, where it is observed that Schreier graphs of
self-similar contracting groups are of polynomial growth, and that in this case the degree of the polynomial
(more precisely of \emph{power}) growth can be non-integer and even a transcendental number. Interesting
results on the growth of Schreier graphs  are obtained by I.~Bondarenko~\cite{bondarenko:PhD} (see also
\cite{benjamini_h:omega_per_graphs,bondarenko-3:family11}). The amenability of Schreier graphs associated
with
 actions of almost finitary groups on the boundary of rooted tree is proven in~\cite{grigorchuk-n:amenable05}.

\medskip
Finally, let us return to the discussion on  the role of just-infinite groups in the study of  growth. Recall
that a group $G$ is called just-infinite if it is infinite but every proper quotient is finite. Such groups
are on the border between finite groups and infinite groups, and surely they should play an important role in
investigations around various gap type conjectures considered in this article. The following statement is an
easy application of Zorn's lemma.

\begin{prop} \label{prop}
Let $G$ be a finitely generated infinite group.  Then $G$ has a just-infinite quotient.
\end{prop}

\begin{cor}  Let $\mathcal{P}$ be a group theoretical property preserved under taking quotients. If there
is a finitely generated  group satisfying the property  $\mathcal{P}$  then there is a just-infinite group
satisfying this property.
\end{cor}

Although the property of a group to have intermediate growth is not preserved when passing to a quotient
group (the image may have polynomial growth), by  theorems of Gromov~\cite{gromov:poly_growth} and Rosset
\cite{rosset:76}, if the quotient $G/H$ of a group $G$ of intermediate growth is a virtually nilpotent group
then $H$ is a finitely generated group of intermediate growth and one may look for a just-infinite quotient
of $H$ and iterate this process in order to represent $G$ as a consecutive extension of a chain of groups
that are virtually nilpotent or just-infinite groups.  This observation is the base of the arguments for
statements given by theorems~\ref{orderable}, \ref{polyciclic}, \ref{resid} and \ref{wilson}.

The next theorem  was  derived by the author from a result of J.S.~Wilson~\cite{wilson:ji71}.

\begin{thm}\label{just-inf} (\cite{grigorch:jibranch})
The class of just-infinite groups naturally splits into three subclasses:
\begin{itemize}
\item[(B)] Algebraically branch just-infinite groups,
\item[(H)] Hereditary just-infinite groups, and
\item[(S)] near-simple just-infinite groups.
\end{itemize}
\end{thm}

Recall that  branch groups  were already defined in section~\ref{inter1}.  The definition of algebraically
branch groups can be found in~\cite{grigorch:branch,bar_gs:branch}. Every geometrically branch group is
algebraically branch but not vice versa.  The difference between the two versions of the definitions is not
large but still there is no complete understanding how much the two classes differ.  Not every branch group
is just-infinite but every proper quotient of a branch group is virtually abelian. Therefore branch groups
are ``almost just-infinite''  and most of the known finitely generated branch groups are just-infinite.

\begin{defn} A group $G$ is \emph{hereditary
 just-infinite} if it is infinite, residually finite,  and every subgroup $H< G$
of finite index is just-infinite.
\end{defn}

For instance $\mathbb{Z}, D_{\infty}$, and  $PSL(n,\mathbb{Z}),~n\geq 3$  (by the result of G.Margulis
\cite{margulis:book91}) are hereditary just-infinite groups.

\begin{defn}
We call a  group $G$ \emph{near-simple} if it contains a subgroup of finite index $H$ which is a direct
product
$$ H=P\times P \dots \times P$$
of finitely many copies of  a \emph{simple} group $P$.
\end{defn}

We already know that  there  are  finitely generated  branch  groups of intermediate growth (for
 instance  groups $G_{\omega}, \omega \in \Omega_1$).
The question on the existence of non-elementary amenable hereditary just-infinite groups is still  open
(observe that the only elementary amenable hereditary just-infinite groups  are $\mathbb{Z}$ and
$D_{\infty}$).

\begin{prob}\label{medyn}
Are there  finitely generated hereditary just-infinite groups of intermediate growth?
\end{prob}

\begin{prob}
Are there  finitely generated  simple groups of intermediate growth?
\end{prob}

As it was already mentioned in section~\ref{gapconj}, we believe that there is a reduction of the Gap
Conjecture to the class of just-infinite groups, that is to the classes of (just-infinite) branch groups,
hereditary just-infinite groups and simple groups. The corresponding result would hold if the Gap Conjecture
holds for residually solvable groups. Using the results of Wilson~\cite{wilson:gap11}, one can prove the
following result.

\begin{thm}(\cite{grigorch:gapconj12}) \label{wilson}\quad\\[-\baselineskip]
\begin{itemize}
\item[(i)] If the Gap Conjecture with parameter $1/6$ holds for just-infinite groups,
 then it holds for all groups.

\item[(ii)] If the Gap Conjecture holds for residually polycyclic groups and for just-infinite groups,
 then it holds for all groups.
\end{itemize}
\end{thm}

Therefore to obtain a complete reduction of the Gap Conjecture to just-infinite groups it is enough to prove
it for residually polyciclic groups which is quite plausible.  Similar reductions hold for some other gap
type conjectures stated in this article.

As was already mentioned in section~\ref{gapconj}, uncountably many finitely generated simple groups that
belong to the class $LEF$ were recently constructed by K.~Medynets and the
author~\cite{grigorch_medyn:simple11}. It may happen that among the subgroups considered
in~\cite{grigorch_medyn:simple11} (they are commutator subgroups of \emph{topological full groups} of
\emph{subshifts of finite type}), there are groups of intermediate growth but this has to be checked.  On the
other hand, it may happen that there are no simple groups of intermediate growth and that there are no
hereditary just-infinite groups of intermediate growth at all.
In this case the Gap Conjecture would be reduced to
the case of branch groups. In~\cite{grigorch:degrees,grigorch:degrees85} the author proved that growth
functions of all $p$-groups of intermediate growth $\mathcal G_{\omega}$ discussed in section~\ref{inter1}
satisfy the lower bound $\gamma_{\mathcal G_{\omega}}(n)\succeq e^{\sqrt n}$, and this was proved by direct
computations based on the anti-contracting property given by definition~\ref{growth4}.  This gives some hope
that the Gap Conjecture can be proved for the class of branch groups by a similar method.

\bibliographystyle{amsalpha}
\newcommand{\etalchar}[1]{$^{#1}$}
\def\cprime{$'$} \def\cprime{$'$} \def\cprime{$'$} \def\cprime{$'$}
  \def\cprime{$'$} \def\cprime{$'$} \def\cprime{$'$} \def\cprime{$'$}
  \def\cprime{$'$}
\providecommand{\bysame}{\leavevmode\hbox to3em{\hrulefill}\thinspace}
\providecommand\arXiv[1]{\href{http://arxiv.org/abs/#1}{\quad arXiv:#1}}
\providecommand{\MR}{\relax\ifhmode\unskip\space\fi MR }
\providecommand{\MRhref}[2]{%
  \href{http://www.ams.org/mathscinet-getitem?mr=#1}{#2}
}
\providecommand{\href}[2]{#2}

\end{document}